\documentclass[11pt]{amsart}
\usepackage{amscd, amssymb,latexsym,amsmath, amscd, amsmath}
\usepackage[all]{xy}
\usepackage{anysize}
\usepackage[sc]{mathpazo} 
\usepackage{color}
\usepackage{tikz-cd}
\usepackage{tikz}
\usetikzlibrary{matrix, arrows}
\marginsize{2.1cm}{2.1cm}{2.1cm}{2.1cm}
\usepackage{braket}
\usepackage{enumitem}
\usepackage{appendix}
\newcommand{\fv}{{\mathfrak{v}}}
\newcommand{\Smod}{\mathcal{S}\mathrm{mod}}
\newcommand{\Span}{{\operatorname{Span}}}
\newcommand{\Hom}{\operatorname{Hom}}
\newcommand{\Ind}{\operatorname{Ind}}
\newcommand{\Rad}{\operatorname{Rad}}
\newcommand{\Supp}{\operatorname{Supp}}
\newcommand{\sign}{\operatorname{sign}}
\newcommand{\SInd}{\mathcal{S}\operatorname{Ind}}
\newcommand{\rmh}{\mathrm{H}}

\newtheorem{theorem}{Theorem}[section] 
\newtheorem{lemma}[theorem]{Lemma}
\newtheorem{corollary}[theorem]{Corollary}
\newtheorem{proposition}[theorem]{Proposition}
\theoremstyle{definition}
\newtheorem{conjecture}[theorem]{Conjecture}

\newtheorem{example}{Example}
\newtheorem{definition}{Definition}
\newtheorem{question}{Question}
\theoremstyle{remark}
\newtheorem{remark}{Remark}
\newtheorem{Case}{Case}
\newtheorem*{claim}{Claim}

\begin{document}

\title[On distinguished representations of $\mathrm{GL}_{n}(\mathbb{C})$]{Classification of $\mathrm{GL}_{n}(\mathbb{C})$-Representations Distinguished by $\mathrm{GL}_n(\mathbb{R})$}
\author[Pattanayak,\ Wu \& Zhang]{Basudev Pattanayak, Kaidi Wu and Hongfeng Zhang}

\address{(Pattanayak) Department of Mathematics, Technion-Israel Institute of Technology, Haifa, Israel.}
\email{pbasudev93@gmail.com }
\address{(Wu) Department of Mathematics and New Cornerstone Science Laboratory, The University of Hong Kong, HK.}
\email{kaidiwu24@connect.hku.hk}
\address{(Zhang) School of Mathematics and Statistics, Huazhong University of Science and Technology, Wuhan, 430074, China}
\email{zhanghongf@pku.edu.cn}


\date{}
\begin{abstract}
This paper provides a complete classification of $\mathrm{GL}_n(\mathbb{R})$-distinguished irreducible representations of $\mathrm{GL}_n(\mathbb{C})$ when the representations are either generic or unitary. Additionally, for each such $\mathrm{GL}_n(\mathbb{R})$-distinguished representation, we explicitly construct the associated period and prove its non-vanishing on the distinguished minimal $K$-type excluding certain non-generic cases. Furthermore, we offer some applications to the branching problem using theta correspondence.
\end{abstract}
\maketitle
\tableofcontents

\section{Introduction}
Let $H$ be a subgroup of $G$ and $\chi$ be a character of $H$. We say that a representation $\pi$ of $G$ is $(H, \chi)$-\textbf{distinguished} if 
\[\mathrm{Hom}_H(\pi, \chi) \neq 0.\]
 If $\chi$ is the trivial character of $H$, then such representations are called $H$-distinguished representations of $G$. One of the most important problems in the relative Langlands program is the classification of distinguished representations (see \cite{Prasad}). Over the last three decades, classifying distinguished representations of $G=\mathrm{GL}_n$ over local and global fields has been a very active area of research, see \cite{AP} for earlier work and \cite{Off, Mat, FLO, Mat14, Kem15a, ST} for more recent works. By the multiplicity-one theorem \cite{AG09}, for any irreducible representation $\pi$ of $\mathrm{GL}_n(\mathbb{C})$,
 \begin{equation}
     \dim \Hom_{\mathrm{GL}_n(\mathbb{R})}(\pi,\mathbb{C}) \leq 1.
 \end{equation} 
A more difficult problem is to figure out when this dimension is non-zero.

 In~\cite[Theorem 1.3]{Kem15a}, Kemarsky established a necessary condition. Suppose the Langlands parameter of $\pi$ is $\lambda=\oplus_{j=1}^n\lambda_j$, where $\lambda_j$ is a character of the complex Weil group $W_{\mathbb{C}}= \mathbb{C}^{\times}$. Then $\pi$ is $\mathrm{GL}_n(\mathbb{R})$-distinguished only if there exists an involution $w\in S_n$ such that $\lambda_j=\overline{\lambda_{w(j)}}^{-1}$ for $1\leq j\leq n$, and for any $j$ such that $w(j)=j$, $\lambda_j(-1)=1$. Here, the character $\overline{\lambda_j}$ is defined by $\overline{\lambda_j}(z):=\lambda_j(\overline{z})$.
 
In this paper, we provide a characterization of irreducible unitary $\mathrm{GL}_n(\mathbb{R})$-distinguished representations and generic $\mathrm{GL}_n(\mathbb{R})$-distinguished representations of $\mathrm{GL}_n(\mathbb{C})$ in terms of their Langlands parameters. 
 We note that the results for $\chi$-distinguished representations follow immediately from that for distinguished representations, since $\chi$ can be extended to a character of $\mathrm{GL}_n(\mathbb{C})$. Now, we formulate our main results.
 
\begin{theorem}[Generic case]\label{main result_generic}
Assume that $\pi$ is an irreducible generic representation of $\mathrm{GL}_n(\mathbb{C})$, which corresponds to the Langlands parameter $\lambda=\oplus_{j=1}^n\lambda_j$, where $\lambda_j$ is a character of $\mathbb{C}^{\times}$. Then $\pi$ is $\mathrm{GL}_n(\mathbb{R})$-distinguished if and only if there exists an involution $w\in S_n$ such that $\lambda_j=\overline{\lambda_{w(j)}}^{-1}$ for $1\leq j\leq n$; and for any $j$ such that $w(j)=j$, $\lambda_j(-1)=1$. 
\end{theorem}
The necessary condition of Theorem \ref{main result_generic} was already proved by Kemarsky \cite[Theorem 1.3]{Kem15a}. Therefore, we only need to prove the sufficiency direction (see Section \ref{sec:proof_gen}). To facilitate the formulation of the unitary case, we define a character of $\mathbb{C}^\times$ by
\[
\varkappa_{m,s}: z \longmapsto \left( \frac{z}{|z|} \right)^m |z|^{2s} \quad \text{for some } m \in \mathbb{Z} \text{ and } s \in \mathbb{C}.
\]
Here, $\mathrm{Re}\,s$ is called the real part of the character.
\begin{theorem}[Unitary case]\label{main result_unitary}
Assume that $\pi$ is an irreducible unitary representation of $\mathrm{GL}_n(\mathbb{C})$, which corresponds to the Langlands parameter $\lambda=\oplus_{j=1}^n\lambda_j$, where $\lambda_j=\varkappa_{m_{j},s_{j}}$ for some $m_j \in \mathbb{Z} \text{ and }s_j \in \mathbb{C}$. Then $\pi$ is $\mathrm{GL}_n(\mathbb{R})$-distinguished if and only if the parameter $\lambda$ satisfies the following two conditions:
\begin{enumerate}
\item[(i)] there exists an involution $w\in S_n$ such that  $\lambda_j=\overline{\lambda_{w(j)}}^{-1}$ for $1\leq j\leq n$; and for any $1\leq j\leq n$ such that $w(j)=j$, $\lambda_j(-1)=1$. 
\item[(ii)] for any $\lambda_j$ such that $m_j$ is odd and $2s_j\in \mathbb{Z}$, the multiplicity of $\lambda_j$ in the Langlands parameter $\lambda$ is even.
\end{enumerate}
\end{theorem}

 We can equivalently rephrase our Theorem \ref{main result_unitary}, see Theorem \ref{equi_thm}, which can be compared to the $p$-adic cases. Let $\pi$ be the irreducible unitary representation of $\mathrm{GL}_n(\mathbb{C})$. We can write $\pi$ as a product of some unitary characters and complementary series. Then $\pi$ is $\mathrm{GL}_n(\mathbb{R})$-distinguished if and only if the contragredient $\pi^{\vee}\simeq \overline{\pi}$, and when $\left(\frac{\det}{|\det|}\right)^r_{\mathrm{GL}_m(\mathbb{C})}$ has odd multiplicity in the product, $\left(\frac{\det}{|\det|}\right)^r_{\mathrm{GL}_m(\mathbb{C})}$ is $\mathrm{GL}_m(\mathbb{R})$-distinguished (i.e. $r$ is an even number). This version of our statement is comparable to the $p$-adic situation that occurred in the main theorem of \cite[Introduction]{Mat14}. For the generic case, it is also analogous to the result \cite[Theorem 5.2]{Mat} for $p$-adic general linear groups. 

 In the proof of the necessity of condition (ii), we utilize the Archimedean Bernstein-Zelevinsky derivative and filtration, which originated  from~\cite{AGS15a,AGS15b} and were further developed in~\cite{WZ25,WY26}. We expect that this method will extend to other Galois models.

As a byproduct of Theorem~\ref{main result_generic} and Theorem~\ref{main result_unitary}, we establish a period integral using the degenerate Whittaker model of irreducible unitary distinguished representations and propose a family of invariants for these representations based on this period integral, see Section~\ref{derivative section}. 

After getting a distinguished representation $(\pi,V)$, there is a more refined question about finding \textbf{test vectors} $v\in V$, such that $v$ is nonzero under the period. We call a $K$-type of $\pi$  \textbf{distinguished minimal $K$-type} if it has minimal norm among all the $K$-types $\tau$ of $\pi$ with $\mathrm{Hom}_{\mathrm{O}_n(\mathbb{R})}(\tau,\mathbb{C})\neq 0$, see Definition~\ref{distinguished minimal}.
We show that such a test vector $v$ can be found in the distinguished minimal $K$-type except some special cases (see \eqref{exception}):
\begin{theorem}\label{Non}
    Let $\pi$ be a distinguished irreducible representation of $\mathrm{GL}_n(\mathbb{C})$, and $\mathcal{P}\in \mathrm{Hom}_{\mathrm{GL}_n(\mathbb{R})}(\pi,\mathbb{C})$ be a non-zero period. When $\pi$ is either generic or unitary, except for the non-generic unitary cases specified in \eqref{exception}, $\mathcal{P}$ is non-vanishing on the distinguished minimal $K$-type of $\pi$.
\end{theorem}
In the proof, we introduce a novel approach to constructing periods, see Subsection~\ref{construction of period}. We study the parabolic induction of the period in Subsection~\ref{parabolic induced}. We first construct the period on certain principal series via intertwining operators and period integrals on the basic distinguished $\mathrm{GL}_2(\mathbb{C})$-representations. We realize the representations of interest as quotients of these principal series and then obtain the period on them by applying a multiplicity-one result (Lemma~\ref{mult<=1}) for these principal series together with a continuity argument. 

There exists another compelling reason for considering distinguished minimal $K$-types. For an irreducible representation $\pi$ with minimal $K$-type $\tau$ which is distinguished, we can consistently identify a parabolically induced representation where $\pi$ emerges as the unique subquotient containing $\tau$, see \cite{Vo79}. Given that the Schwartz homology of parabolic induced representations is more tractable to compute and their periods are more straightforward to construct, this approach provides an effective method for determining whether $\pi$ is distinguished.

From our Theorem \ref{main result_generic}, we get some branching laws of $\mathrm{GL}_{2m}(\mathbb{R}) \downarrow \mathrm{GL}_m(\mathbb{C})$ by applying seesaw identity.
\begin{corollary}\label{new-branching}
   Let $\pi$ be an irreducible generic representation of $\mathrm{GL}_n(\mathbb{C})$, with the Langlands parameter $\oplus_{j=1}^n \lambda_j$. Then it appears in the branching of $\mathrm{Ind}_{\mathrm{GL}_n(\mathbb{R})\times \mathrm{GL}_n(\mathbb{R})}^{\mathrm{GL}_{2n}(\mathbb{R})}(\mathbb{1}\boxtimes \mathbb{1})|_{\mathrm{GL}_n(\mathbb{C})}$, which means
   \begin{equation*}
       \mathrm{Hom}_{\mathrm{GL}_n(\mathbb{C})}\left(\mathrm{Ind}_{\mathrm{GL}_n(\mathbb{R})\times \mathrm{GL}_n(\mathbb{R})}^{\mathrm{GL}_{2n}(\mathbb{R})}(\mathbb{1}\boxtimes \mathbb{1})|_{\mathrm{GL}_n(\mathbb{C})},\pi\right)\neq 0,
   \end{equation*} if and only if there exists an involution $w\in S_n$ such that $\lambda_j=\overline{\lambda_{w(j)}}^{-1}$ for all $1\leq j\leq n$, and for any $j$ such that $w(j)=j$, $\lambda_j(-1)=1$. Here, $\Ind$ denotes the standard parabolic induction.
\end{corollary}
The local factors like $\varepsilon$-factor, $\gamma$-factor have arithmetic
 significance and their explicit values are of fundamental importance in the Langlands program. We have the following triviality of  the $\varepsilon$-factor for distinguished representations: 
\begin{theorem}\label{factor}
 Let $\pi$ be an irreducible $\mathrm{GL}_n(\mathbb{R})$-distinguished representation of $\mathrm{GL}_n(\mathbb{C})$,  and $\psi$ be an additive character of $\mathbb{C}$ such that $\psi\big|_\mathbb{R}=1$. Then
 \[\varepsilon\left(\frac{1}{2}, \pi , \psi\right)=1.\]  
\end{theorem}
By utilizing the strategy of the proof of Theorem \ref{factor} for distinguished representations, we reproduce the result in \cite[Theorem 6.3]{MO}.

Finally, we emphasize that our classification of distinguished representations and our non-vanishing result are closely related to the Asai $L$-function. Our findings immediately identify all the nearly tempered representations of $\mathrm{GL}_n(\mathbb{C})$ whose Asai $L$-function has an exceptional pole of level zero at $0$, see \cite[Theorem 1.1]{Yad} for details. Moreover, Theorem~\ref{Non} recovers~\cite[Theorem 4.1]{Sun} in a completely different way, thereby reproving the non-vanishing hypothesis of Grobner-Harris-Lapid in the study of non-critical values of the Asai $L$-function, see~\cite[Theorem 4.5]{Sun} and~\cite[Section 6.2]{GHL}.

\subsection{Structure of the paper}
In Section~\ref{Pre}, we introduce notation, general methods, and preliminary results that will be utilized in the subsequent sections. In Section~\ref{sec:cls}, an equivalent necessary and sufficient condition for Theorem~\ref{main result_unitary} is presented; see Theorem~\ref{equi_thm}. This equivalent theorem forms the primary focus of our proof.

In Section~\ref{Necessary}, we prove the necessary condition. Part (i) of Theorem~\ref{equi_thm} is already known from Kemarsky's work, so we reduce the necessary condition to part (ii) by applying Lemma~\ref{red_int}. We then use the Archimedean Bernstein-Zelevinsky derivative to prove part (ii); see Theorem~\ref{deri_thm}. In Section~\ref{sufficient}, we demonstrate that any irreducible representation satisfying the necessary condition is indeed distinguished. We show that the distinction is preserved under parabolic induction. By leveraging this property, we reduce the problem to Lemma~\ref{dist_comp}. Later, we prove the sufficient parts of both Theorem \ref{main result_generic} and \ref{main result_unitary} respectively in Section \ref{sec:proof_gen} and \ref{sec:proof_uni}.

In Section~\ref{non-vanishing}, we prove Theorem~\ref{Non} through the generic calculation of the distribution kernel and an explicit construction of the period; see \eqref{period def} for the period map. The non-vanishing property follows from tracing the distinguished minimal $K$-type when applying intertwining operators. Some technical details are provided in Appendix~\ref{appendix B}. In Section~\ref{local factor}, we characterize the distinguished representation in terms of the $\epsilon$-factor. Finally, in Section~\ref{application}, we present an application of our main result. We note that the see-saw identity transforms our result into the branching law of $\mathrm{GL}_{2m}(\mathbb{R}) \downarrow \mathrm{GL}_m(\mathbb{C})$.

\bigskip
\centerline{\scshape Acknowledgements}
We would like to thank Dipendra Prasad, Huajian Xue, Zhibin Geng, Nadir Matringe for their valuable suggestions along our research on this problem. We also thank Rui Chen for confirming Theorem~\ref{Theta}. We are deeply grateful to the anonymous referee for their thorough review and extremely helpful suggestions, which greatly strengthened the quality of our work.

Pattanayak and Zhang are supported by the Research Grants Council of the Hong Kong Special Administrative Region, China (Project No: 17305223) and the NSFC grant for Excellent Young Scholar (Project No: 12322120), and they also benefit from working with their supervisor Kei Yuen Chan. Wu is partially supported by the New Cornerstone Science Foundation through the New Cornerstone Investigator Program awarded to Professor Xuhua He. Wu is also supported by the National Natural Science Foundation of China (Grant No. 123B1004). He thanks Xuhua He and Binyong Sun for their continued support and encouragement. Zhang is also supported by National Natural Science Foundation of China (Grant No. 12601050, Grant No. 12271460 and  Grant No. 12341101).

\section{Preliminary}\label{Pre}
 \subsection{Notation}

Throughout, $G_n$ denotes $\mathrm{GL}_n(\mathbb{C})$, and $G_n'$ denotes $\mathrm{GL}_n(\mathbb{R})$. Unless otherwise specified, 
\begin{itemize}
    \item $G$ represents a Lie group, and $H$ represents a closed subgroup.
    \item the corresponding lowercase Gothic letter denotes the Lie algebra of a Lie group.
    \item $Z_G$ denotes the center of $G$.
\end{itemize}
When $G$ is an almost linear Nash group, we denote by $\Smod_G$ the category of smooth moderate growth Fr\'echet representations of $G$. This is the category with which we are concerned. If $G$ is reductive, we denote by $\mathrm{CW}_G$ the category of Casselman-Wallach representations, namely, the subcategory of $\Smod_G$ whose objects are finite length. For a representation $\pi$ of $G$, let $\pi^{\vee}$ denote the contragredient of $\pi$. We use $\omega_\pi$ to denote the central character of $Z_G$ on the representation $\pi$ which admits central character. For a representation $\pi$ of $G_n$, the representation $\overline{\pi}$ is defined by $\overline{\pi}(g):=\pi(\overline{g})$, where $\overline{g}$ is the complex conjugation on $G_n$. 

We have several closed subgroups of $G_n$ (respectively $G_n'$) defined as follows:
\begin{itemize}
\item $G_m$ is a subgroup of $G_n$ by the embedding $g \mapsto \begin{bmatrix}
    g & \\
     & I_{n-m}
\end{bmatrix}$ for $m\leq n$;
\item $A_m$ : the central torus of $G_m$. For $m\leq n$, it is a subgroup of $G_n$ as well. We fix the isomorphism $t_m:\mathbb{C}^{\times}\rightarrow A_m$ by sending $z$ to the diagonal matrix with each entry $z$. We use $m(a_{p},\dots,a_{q})$ to denote $t_p(a_p)\dots t_q(a_q)\in G_n$ for $1\leq p\leq q \leq n$.
\item $B_n$: the Borel subgroup of $G_n$, consisting of upper triangular matrices, and $N_n$ denotes the unipotent radical of $B_n$; 
\item $K=\mathrm{U}(n)$ ( resp. $K'=\mathrm{O}(n)$) denotes the standard maximal compact subgroup of $G_n$ (resp. $G_n'$);
\item $P_n$: the mirabolic subgroup of $G_n$, consisting of matrices with last row $(0,\dots,0,1)$;
\item $V_n$: the subgroup of $P_n$, consisting of matrices of the form $\begin{bmatrix}
        I_{n-1} & v\\
                &  1
    \end{bmatrix}$;
    \item $H_{n,d}$: the subgroup of $P_n$ of the form $h(a,u)=\begin{bmatrix}
        a & x\\
        0        &  u
    \end{bmatrix}$, with $a \in G_{n-d}$, $u \in N_{d}$, and $x$ an $(n-d,d)$-matrix. Here, $\det h(a,u)$ is defined as $\det a$. When $d=2$, we will write $H_{n,d}$ simply as $H_n$.
\end{itemize}
We also fix the following characters:
\begin{itemize}
    \item $\psi$: a character of $\mathbb{C}$ trivial on $\mathbb{R}$, defined by $\psi(x)=e^{\pi( x-\overline{x} )}$; 
    \item $\psi_n$: a character of $V_n$ defined by $\psi_n(\begin{bmatrix}
        I_{n-1} & v\\
                &  1
    \end{bmatrix}):=\psi(x_{n-1})$, for $v=[x_1,\dots,x_{n-1}]^{t}\in \mathbb{C}^{n-1}$; and also, denote by $\psi_n$ the corresponding character of the Lie algebra $\mathfrak{v}_n$ of $V_n$. If $\pi$ is a smooth moderate growth Fr\'echet representation of $P_n$, define $\pi(\psi_n):=\Span\{\alpha w - \psi_n(\alpha)w \, : \, w \in \pi, \, \alpha \in \fv_{n}\}$.
\end{itemize}
For a partition $\sigma=(\sigma_1,\dots,\sigma_r)$ of $n$, define a family of degenerate characters $$\psi^{\sigma}_j:=\psi(\sum_{i=1}^{\sigma_j-1}n_{i,i+1}+\sum_{i=\sigma_j+1}^{\sigma_j+\sigma_{j-1}-1}n_{i,i+1}+\dots +\sum_{i=\sigma_j+\dots+\sigma_2+1}^{\sigma_j+\dots+\sigma_1-1}n_{i,i+1}).$$
Here $\psi^{\sigma}_j$ is a character of $N_{\sigma_j+\dots+\sigma_1}$ for each $1\leq j\leq r$.

For a real Lie group $G$, we will always use $dg$ to denote its right-invariant Haar measure (up to some scalar), and use $\delta_G$ to denote its modular character. Given a closed subgroup $H$ of $G$, we will consider the ‘quotient measure’ on $H\backslash G$, which is, strictly speaking, not necessarily a measure but rather a continuous functional on the space of continuous left $(H, \frac{\delta_H}{\delta_G})$-equivariant functions on $G$ that are compactly supported modulo $H$. Thus,

$$ \int\limits_{H\backslash G} \int\limits_H \frac{\delta_G(h)}{\delta_H(h)} f(hx)  dh dx = \int\limits_G f(g) dg $$
for any compactly supported continuous function $f$ on $G$.

For a locally convex topological vector space $V$, we usually use $V'$ to denote its continuous dual and endow it with the strong topology if necessary. To employ the Schwartz homology, we will also consider the Nash group structure and tempered Fréchet bundles on Nash manifolds. For more details, we refer to \cite{CS}. 

Let $X$ be a Nash manifold, and let $\mathcal{V}$ be a tempered Fr\'echet bundle on $X$. Let $\Omega_X := \wedge^{\text{top}} T^* X$ be the  canonical bundle on $X$, and let $|\Omega_X|$ denote the density bundle. Then we will use $\mathcal{S}(X, \mathcal{V})$ to denote the Schwartz sections of $\mathcal{V}$, and $\mathcal{G}(X, \mathcal{V}) := \mathcal{S}(X, \mathcal{V}^{\vee} \otimes |\Omega_X|)'$ to denote the generalized sections, where $\mathcal{V}^{\vee}$ is the dual bundle. Write $\mathcal{V}^* := \mathcal{V}^{\vee} \otimes |\Omega_X|$ for simplicity. If $X$ has a right action of the Nash group $G$ and $\mathcal{V}$ is a right $G$-bundle, then $\mathcal{S}(X, \mathcal{V})$ and $\mathcal{G}(X, \mathcal{V})$ are both representations of $G$ via right translation. Here, we are only concerned with the following two cases: $\mathcal{V}$ is finite-dimensional, or $X$ is homogeneous and $\mathcal{V}$ comes from a Casselman-Wallach representation of a reductive subgroup of $G$. In the latter case, we define the dual bundle as a Fr\'echet bundle using the contragredient representation.

 All inductions in this paper are normalized inductions. Let $P$ be a Nash subgroup of an almost linear Nash group $G$, and let $(\sigma,V_{\sigma})$ be a moderate growth smooth Fr\'echet representation of $P$. Then, we denote
\begin{itemize}
    \item $\Ind_P^G(\sigma):=\{f\in \mathcal{C}^{\infty}(G,V_{\sigma})| f(p g)=\delta_{P}(p)^{1/2}\delta_{G}(p)^{-1/2}\sigma(p)f(g)\}$;
    \item normalized Schwartz induction by $\mathcal{S}\Ind_P^G(\sigma)$, see \cite{CS} for more details.
    \item normalized compact induction by $c\text{-}\Ind_P^G(\sigma)$, as the subrepresentation of $\mathcal{S}\Ind_P^G(\sigma)$ consisting of sections with compact support in $P\backslash G$.
\end{itemize}
These inductions are equipped with $G$'s actions by right translation. When $G=P_n,P=H_{n,2}=P_{n-1}\ltimes V_n$ and $\kappa$ is a representation of $P_{n-1}$, we will simply denote $\mathcal{S}\Ind_P^G(\kappa\otimes \psi_n)$ by $I(\kappa)$, which is called the Mackey induction.

Assume that $\pi_i$ is a moderate growth smooth Fr\'echet representation of $\mathrm{GL}_{n_i}(\mathbb{F})$ for $i=1,2$, where $\mathbb{F}=\mathbb{R}$ or $\mathbb{C}$. We use $\pi_1\times \pi_2$ to denote the representation of $\mathrm{GL}_{n_1+n_2}(\mathbb{F})$, which is parabolically induced from the standard block upper triangular subgroup with the Levi subgroup $\mathrm{GL}_{n_1}(\mathbb{F})\times \mathrm{GL}_{n_2}(\mathbb{F})$.

\subsection{Archimedean local Langlands correspondence}\label{sec:cls}
The building blocks for irreducible admissible representations of $G_n$ are characters of $\mathrm{GL}_1(\mathbb{C})=\mathbb{C}^\times$, which are defined by
\[\varkappa_{m,s}: ~ z \mapsto \left( \frac{z}{|z|} \right)^m  |z|^{2s}    \text{ for some } m \in \mathbb{Z} \text{ and }s \in \mathbb{C}.\]
Let $\lambda_i=\varkappa_{m_i,s_i}: z \mapsto \left( \frac{z}{|z|} \right)^{m_i}  |z|^{2s_i}$ be a character of $\mathrm{GL}_1(\mathbb{C})$ for each integer $i$ with $1 \leq i \leq n$. The tuple $(\lambda_1, \lambda_2,\dots ,\lambda_n)$ defines a one-dimensional smooth representation of the diagonal subgroup of $G_n$. Then it extends trivially on $N_n$ to a one-dimensional smooth representation of $B_n$. We denote the normalized parabolic induced representation by \[I(\lambda_1, \cdots, \lambda_n )=\mathrm{Ind}^{G_n}_{B_n} (\lambda_1 \otimes \cdots \otimes \lambda_n)\]
\begin{theorem}[Langlands classification for $\mathrm{GL}_n(\mathbb{C})$ \cite{ZN}]  \label{Lang_cls}Consider the above notation.
   \begin{enumerate}
       \item Suppose the parameters $(s_1,\dots ,s_n)$ of the character $(\lambda_1, \lambda_2,\dots ,\lambda_n)$ satisfy 
       \[\mathrm{Re}(s_1) \geq \mathrm{Re}(s_2) \geq \cdots \geq \mathrm{Re}(s_n).\] Then, the induced representation $I(\lambda_1 , \cdots , \lambda_n )$, called \textbf{standard representations}, has a unique irreducible quotient $Q(\lambda_1 , \cdots , \lambda_n )$
       \item Any irreducible admissible representation of $\mathrm{GL}_n(\mathbb{C})$ is of the form $Q(\lambda_1 , \cdots , \lambda_n )$, up to equivalence.
       \item Two such representations $Q(\lambda_1, \cdots, \lambda_n )$ and $Q(\lambda_1^\prime, \cdots, \lambda_n^\prime )$ are infinitesimally equivalent if and only if there exists a permutation $\omega \in S_n$ such that $\lambda_i= \lambda_{\omega(i)}^\prime$ for $1 \leq i \leq n$.
   \end{enumerate} 
\end{theorem}
Let $W_\mathbb{F}$ denote the Weil group of archimedean local field $\mathbb{F}$, so $W_\mathbb{C} = \mathbb{C}^\times$. We are interested in the set of equivalence classes of $n$-dimensional complex representations (denoted by $\lambda$'s) of $W_\mathbb{C}$ whose images consist of semisimple elements. Since $W_\mathbb{C}$ is abelian, every such semisimple $n$-dimensional complex representation $\lambda$ can be written as a direct sum of one-dimensional representations $\kappa_j=\varkappa_{m_j, s_j}$ for $1 \leq j \leq n$. Fix a permutation $\omega \in S_n$ such that 
\[\mathrm{Re}(s_{\omega(1)}) \geq \mathrm{Re}(s_{\omega(2)}) \geq \cdots \geq \mathrm{Re}(s_{\omega(n)}).\]
We now define the one-dimensional representation $\lambda_j$ of $\mathrm{GL}_1(\mathbb{C})$ by
\begin{equation}\label{explicit:llc}
 \lambda_j= \kappa_{\omega(j)}=\varkappa_{m_{\omega(j)},s_{\omega(j)}}   
\end{equation}
for each integer $j$ with $1 \leq j \leq n$. 
\begin{theorem}[Local Langlands correspondence, ``LLC" for short]Consider the above association (\ref{explicit:llc}). Using Theorem \ref{Lang_cls}, the map defined by
   \begin{equation}\label{llc}
     \lambda=\bigoplus\limits_{j=1}^n \lambda_j   \mapsto Q(\lambda_1,\dots ,\lambda_n)  
   \end{equation} 
   is a well-defined bijection between the set of equivalence classes of $n$-dimensional semisimple complex representations of $W_\mathbb{C}$ and the set of equivalence classes of irreducible admissible complex representations of $\mathrm{GL}_n(\mathbb{C})$. This $\lambda$ is called the \textbf{Langlands parameter} of the irreducible admissible representation $\pi=Q(\lambda_1,\dots ,\lambda_n)$.
\end{theorem}

\subsubsection{Vogan's classification} 
In his phenomenal work \cite{Vog86}, Vogan classifies all irreducible unitary representations of $\mathrm{GL}_n(\mathbb C)$ by the parabolic induction from unitary characters and Stein complementary series.

Let ${\underline n}=(n_1,\dots,n_r,2m_1,\dots,2m_s)$ be a partition of $n$, that is, $\sum_{i=1}^r n_i+\sum_{j=1}^s2m_j=n$,  and the corresponding standard Levi subgroup $L_{\underline n}=\mathrm{GL}_{n_1}(\mathbb C)\times\dots\times \mathrm{GL}_{n_r}(\mathbb C)\times \mathrm{GL}_{2m_1}(\mathbb C)\times\dots\times \mathrm{GL}_{2m_s}(\mathbb C)$. 

For $1\leq i\leq r$, define unitary characters $\chi_{i}: \mathrm{GL}_{n_i}(\mathbb C) \rightarrow \mathbb{C}^\times$ of $\mathrm{GL}_{n_i}(\mathbb C)$ by
\[\chi_{i}(g)= \left( \frac{\mathrm{det}(g)}{|\mathrm{det}(g)|}\right)^{k_i} \cdot |\mathrm{det}(g)|^{u_i}   ~~\text{ for }k_i \in \mathbb{Z}, u_i\in \sqrt{-1} \mathbb{R}.\] For $1\leq j\leq s$,  define non-unitary characters $\gamma_j, \gamma_j^{\flat}: \mathrm{GL}_{m_j}(\mathbb C) \rightarrow \mathbb{C}^\times$ of $\mathrm{GL}_{m_j}(\mathbb C)$ by
\[\gamma_j(g)= (\frac{\det}{|\det|})_{G_{m_j}}^{k_j'}\cdot |\det|_{G_{m_j}}^{u_j'+t_j},\quad  \gamma_j^{\flat}(g)= (\frac{\det}{|\det|})_{G_{m_j}}^{k_j'}\cdot |\det|_{G_{m_j}}^{u_j'-t_j}\]
for $k_j' \in \mathbb{Z}, u_j'\in \sqrt{-1} \mathbb{R}, 0 < t_j <1$. The induced representation $ \sigma_j=\gamma_j\times \gamma_j^{\flat}$ is an irreducible unitary but non-tempered representation--these are called the complementary series of $G_{2m_j}$. Moreover, given two complementary series of $G_{2m_j}$ corresponding to $(u_j',k_j',t_j)$ and $(\widetilde{u_j'},\widetilde{k_j'},\widetilde{t_j})$, they are isomorphic if and only if their parameters are the same, that is, $u_j'=\widetilde{u_j'}$, $k_j'=\widetilde{k_j'}$ and $t_j=\widetilde{t_j})$.

The representation $\chi_{1} \otimes\dots \otimes \chi_{r} \otimes \sigma_{1} \otimes\dots \otimes \sigma_{s}$ defines a unitary representation of $L_{\underline n}$. Then, the normalized parabolic induced representation 
\[I(\chi_{1},\dots ,\chi_{r},\sigma_1,\dots,\sigma_s)=\mathrm{Ind}^{\mathrm{GL}_n(\mathbb{C})}_{P_{\underline n}}(\chi_{1} \otimes\dots \otimes \chi_{r}\otimes\sigma_1\otimes\dots\otimes \sigma_s)\] is an irreducible unitary representation of  $\mathrm{GL}_n(\mathbb{C})$, where $P_{\underline n}$ is the parabolic subgroup with Levi $L_{\underline n}$.

Every irreducible unitary representation of $\mathrm{GL}_n(\mathbb{C})$ is of the form $I(\chi_{1},\dots ,\chi_{r},\sigma_1,\dots,\sigma_s)$ for some partition $\underline n$ and unitary characters $\chi_{i}$ and complementary series $\sigma_j$.

Furthermore,  given another representation  $I(\widetilde{\chi}_{1},\dots ,\widetilde{\chi}_{\widetilde{r}},\widetilde{\sigma}_1,\dots,\widetilde{\sigma}_{\widetilde{s}})$ of the above form, 
\[I(\widetilde{\chi}_{1},\dots ,\widetilde{\chi}_{\widetilde{r}},\widetilde{\sigma}_1,\dots,\widetilde{\sigma}_{\widetilde{s}})\simeq I(\chi_{1},\dots ,\chi_{r},\sigma_1,\dots,\sigma_s)\] if and only if as multisets,
\[\{\widetilde{\chi}_{i}\ |\ 1\leq i\leq \widetilde{r}\}=\{\chi_i\ |\ 1\leq i\leq r\} \text{ and }\{\widetilde{\sigma}_j\ |\ 1\leq j\leq \widetilde{s}\}=\{\sigma_j\ |\ 1\leq j\leq s\}\]

For each $G_{n_i}$ factor, the unitary character $\chi_{i}=\left( \frac{\mathrm{det}}{|\mathrm{det}|}\right)^{k_i} \cdot |\mathrm{det}|^{u_i}$ corresponds under LLC to a $n_i$-dimensional representation $\lambda_i$ of the Weil group $W_\mathbb{C}$ defined by
\[\chi_{i} \longleftrightarrow \lambda_i:=\bigoplus_{r=1}^{n_i}\left( \frac{z}{|z|}\right)^{k_i} \cdot |z|^{u_i+n_i-2r+1}.\]
For each $G_{m_j}$ factor, the non-unitary character $\gamma_j=(\frac{\det}{|\det|})_{G_{m_j}}^{k_j'}\cdot |\det|_{G_{m_j}}^{u_j'+t_j}$ corresponds under LLC to a $m_j$-dimensional representation $\kappa_j$ of the Weil group $W_\mathbb{C}$ defined by
\[\gamma_j \longleftrightarrow \kappa_j:=\bigoplus_{r=1}^{m_j}\left( \frac{z}{|z|}\right)^{k_j'} \cdot |z|^{u_j'+t_j+m_j-2r+1}.\]
Similarly, $\gamma_j^{\flat}$ corresponds to $\kappa_j^{\flat}$.
Then the full Langlands parameter for $\pi=I(\chi_{1},\dots ,\chi_{r},\sigma_1,\dots,\sigma_s)$ is given by the following $n$-dimensional representation
\[\varphi_\pi=\bigoplus_{i=1}^r\lambda_i \oplus \bigoplus_{j=1}^s(\kappa_j\oplus \kappa_j^{\flat}).\]

By the above classification of the unitary dual of $\mathrm{GL}_n(\mathbb{C})$, one can see that the condition (i) $+$ (ii) in Theorem \ref{main result_unitary} ``$\Longleftrightarrow$" condition (i) $+$  (ii) in the following Theorem \ref{equi_thm}.
\begin{theorem}\label{equi_thm}
Let $\pi$ be an irreducible unitary representation of $G_n$, and write it as 
\[\pi=\chi_1\times \dots \times \chi_r\times\sigma_1\times \dots\times \sigma_s,\] 
where $\chi_i=|\det|^{u_i}\otimes \left(\frac{\det}{|\det|}\right)^{k_i}$ ($u_i\in \sqrt{-1}\mathbb{R}$, $k_i\in \mathbb{Z}$) is a unitary character of $G_{n_i}$, and \[\sigma_j=|\det|_{G_{2m_j}}^{u_j'}\otimes \left(\frac{\det}{|\det|}\right)_{G_{2m_j}}^{k_j'}\otimes \left(|\det|_{G_{m_j}}^{t_j}\times |\det|_{G_{m_j}}^{-t_j}\right)\]
with $u_j'\in \sqrt{-1}\mathbb{R}$, $k_j'\in \mathbb{Z}$, $0<t_j<1$,
is a complementary series of $G_{2m_j}$. 

Then $\pi$ is $G'_n$-distinguished if and only if it satisfies the following conditions
\begin{itemize}
\item[(i).] \begin{itemize}
    \item[(ia)] For $\chi_i$ with $u_i\in\sqrt{-1}\mathbb{R}\backslash \{0\}$, the multiplicity of $\chi_i$  is the same as that of $\chi_{i'}$ with $n_{i'}=n_i$, $k_{i'}=k_i$ and $u_{i'}=-u_i$. 
    
    \item[(ib)] For $\sigma_j$ with $u_j'\in\sqrt{-1}\mathbb{R}\backslash \{0\}$, the multiplicity of $\sigma_j$ is the same as that of $\sigma_{j'}$ with $m_j=m_{j'}$, $u_j'=-u_{j'}'$, $k_j'=k_{j'}'$ and $t_j=t_{j'}$.  
           \end{itemize}
\item[(ii).]  For $\chi_i$ with $u_i=0$ and $k_i\in 2\mathbb{Z}+1$, the multiplicity of $\chi_i$ is even. 
\end{itemize}
\end{theorem} 

In the following sections \ref{Necessary} and \ref{sufficient}, we will prove Theorem \ref{equi_thm}. Before that, let us give an example to illustrate the above theorem.

\begin{example}
Here is an example of a unitary representation $\pi$ of $G_6$ that satisfies condition (i)
but not condition (ii) of Theorem \ref{main result_unitary}. Consider the representation \[ \pi = \left(\frac{\det}{|\det|}\right)_{G_2} \times \left(\frac{\det}{|\det|}\right)_{G_2} \times \left(\frac{\det}{|\det|}\right)_{G_2} \] of $G_6$. Its Langlands parameter is $\oplus_{i=1}^6 \lambda_i$, where $\lambda_i(z) = \frac{z}{|z|} |z|$ for $1 \leq i \leq 3$, and $\lambda_i(z) = \frac{z}{|z|} |z|^{-1}$ for $4 \leq i \leq 6$. By taking involution $w = (1~~4)(2~~5)(3~~6)$ of $S_6$, we can see that $\pi$ satisfies condition (i). However, it does not satisfy condition (ii). Our theorem implies that $\pi$ is not $G_6'$-distinguished. Although $\pi$ is not $G_6'$-distinguished, $\pi$ contains some distinguished $K$-type (for example, the $\mathrm{U}(6)$-type with highest weight $(2,2,2,0,0,0)$).

Notice that the highest derivative of $\pi$ is $ \left(\frac{\det}{|\det|}\right)_{G_1} \times \left(\frac{\det}{|\det|}\right)_{G_1} \times \left(\frac{\det}{|\det|}\right)_{G_1} $ of $G_3$, which does not satisfy condition (i) of Theorem \ref{main result_unitary}, and hence is not $G_3'$-distinguished. In Theorem \ref{deri_thm} of Section \ref{Necessary}, we will prove that the highest derivative of a $G_n'$-distinguished unitary representation $\pi$ is also $G_{n-d}'$-distinguished, where $d$ is the depth of $\pi$.
\end{example}

\subsection{Archimedean Bernstein-Zelevinsky derivative}
The Archimedean analog of the Bernstein-Zelevinsky derivative (which we will simply call the derivative) for Casselman-Wallach representations was introduced in \cite{AGS15a,AGS15b} and developed in~\cite{WZ25,WY26}. We will use the derivative as the main ingredient in Section~\ref{Necessary} to identify the non-distinguished unitary representations. Here, we restate some definitions and facts about derivatives, which we need frequently. 
\begin{definition}
    Let $\pi\in \Smod_{P_n}$. Define
    $$\Psi(\pi):=|\det|^{-1} \otimes \pi/\pi(\psi_n)$$
and
$$\Phi(\pi):= \varprojlim_{l} \pi /\Span\{\beta v\,|\,v  \in \pi ,\beta \in ({\mathfrak{v}}_{n})^{\otimes l} \}.$$
$\Psi(\pi)$ is a representation of $P_{n-1}$ and $\Phi(\pi)$ is a representation of $P_n$. The $k$-th derivative of $\pi$ is defined to be $D^k(\pi):=\Phi\Psi^{k-1}(\pi)$. The depth of representation $\pi$ is defined to be the maximal positive integer $k$ such that $D^k(\pi)\neq 0$, and in this case, we call $D^k(\pi)$ the highest derivative of $\pi$, and denote it simply by $\pi^{-}$.
\end{definition}
By~\cite[Theorem 1.2]{WY26}, $\Psi(\pi)$ is Hausdorff, hence lies in $\Smod_{P_{n-1}}$. Furthermore, by~\cite[Theorem 1.3]{WZ25}, when $\pi$ is restricted from a Casselman-Wallach representation of $G_n$, $D^k(\pi)\in \Smod_{P_{n-k+1}}$. The main results we shall need are as follows:
\begin{theorem}[\cite{AGS15a,AGS15b,WY26}]\label{1.1}
    Let $\pi$ be a Casselman-Wallach representation of $G_n$, and $\mathrm{Ann}(\pi)$ be the annihilator variety of $\pi$. Using non-degenerate invariant bilinear form on $\mathfrak{g}$ and $\mathfrak{g}^*$, $\mathrm{Ann}(\pi)$ is identified with a subvariety of $\mathfrak{g}_n=\mathcal{M}_{n\times n}(\mathbb{C})$. Then
    \begin{enumerate}
        \item The depth of $\pi$ is the minimal positive integer $d$ such that $X^d=0$ for all $ X\in \mathrm{Ann}(\pi)$.
        \item Suppose $\pi$ is of depth $d$. Then $D^d(\pi)=\Psi^{d-1}(\pi)$ is a Casselman-Wallach representation of $G_{n-d}$.
        \item $D^{k+1}:\mathrm{CW}_{G_n}\longrightarrow  \Smod_{P_{n-k}}$ is an exact functor for any non-negative integer $k$.
        \item  Let $n=n_1+\cdots +n_d$ and let $\chi_i$ be characters of $G_{n_i}$. Let $\pi= \chi_1 \times \cdots \times \chi_d $ denote the corresponding monomial representation. Then, by the assertion (1), $\pi$ is of depth $d$. Moreover
    $$D^d(\pi)\cong((\chi_1)|_{G_{n_1-1}} \times  \cdots  \times (\chi_d)|_{G_{n_d-1}}). $$
    \end{enumerate}
\end{theorem}
By (2), we will freely regard $\pi^-$ as a representation of both $P_{n-d+1}$ and $G_{n-d}$. We will also use the following Bernstein-Zelevinsky filtration for the restriction of Casselman-Wallach representations of $G_n$ to $P_n$.
\begin{theorem}[\cite{WY26}, Theorem 1.7]\label{BZ-fil}
    Let $\pi$ be a Casselman-Wallach representation of $G_n$. Then $\pi|_{P_n}$ has a decreasing filtration
    \[
    \pi|_{P_n}=\pi_0\supset \pi_1\supset \dots\supset \pi_n= 0
    \]
    such that
    \[
    \pi_k/\pi_{k+1}\simeq I^{k} ( D^{k+1}(\pi)) \text{ for}\quad 0\leq k\leq n-1,
    \]
    where $I^k$ is the iteration of the Mackey induction.
\end{theorem}

The main idea that the Bernstein-Zelevinsky filtration can be used to identify non-distinguished unitary representations through the highest derivative is the following proposition. Its $p$-adic analogy was proved in~\cite[section 7.3]{Ber}. Note that $D^{k+1}(\pi)$ has a decreasing filtration $F_k^l(\pi):=\mathfrak{v}_{n-k}^l \Psi^{k}(\pi)$ such that $D^{k+1}(\pi)=\varprojlim_l F_k^0(\pi)/F_k^l(\pi)$. By~\cite[Theorem 1.3]{WZ25}, each subquotient $F_k^0(\pi)/F_k^l(\pi)$ is a Casselman-Wallach representation of $G_{n-k-1}$.
\begin{proposition}[\cite{WZ25}, Theorem 1.7]\label{unitary characterization}
    Let $\pi$ be an irreducible unitary representation of $G_n$ of depth $d$. Then for any non-negative integer $k,l$ and any irreducible component $\tau$ in $F_k^0(\pi)/F_k^l(\pi)$, the real part of the central character of $\tau$ is positive
    \[
     \mathrm{Re}\,\omega_{\tau} >0
    \]
    unless $d=k+1$.
\end{proposition}

\subsection{Branching laws of $\mathrm{U}(n)$ and distinguished minimal $K$-type }\label{Non-van}
We show two facts: that the minimal $K$-types (in the sense of the norm defined in \cite{V81}) are likely to be test vectors, and that this phenomenon is closely related to distinction in some cases. 

The first fact appears in \cite{Sun}. If $\pi$ is a distinguished representation coming from cohomological induction, then the period does not vanish on the bottom layer under certain conditions (cf. \cite[Theorem 1.10]{Sun}). In particular, the bottom layer coincides with the minimal $K$-type for irreducible unitary representations with nonzero cohomology.

The second fact appears in \cite{Var}. Suppose $(G, H)$ is a symmetric pair, and both have compact Cartan subgroups. Choose a maximal compact subgroup $K$ of $G$ such that $K_H:=K\cap H$ is a maximal compact subgroup of $H$. Suppose $\pi$ is a discrete series of $G$ with minimal $K$-type $\tau$, and $\rho$ is a discrete series of $H$ with minimal $K_H$-type $\sigma$. Then $\Hom_{K_H}(\tau|_{K_H},\sigma) \neq 0$ implies $\Hom_{H}(\pi|_{H},\rho) \neq 0$ (cf. \cite[Theorem 1]{Var}). 

As Theorem~\ref{Non} states, we show that periods of $(G_n,G_n')$-distinguished representations are non-vanishing on distinguished minimal $\mathrm{U}(n)$-types. More generally, assume that $H$ is a spherical subgroup of the real reductive group $G$, and $K$ is a maximal compact subgroup of $G$ such that $K_H:=K\cap H$ is a maximal compact subgroup of $H$: 
\begin{definition} \label{distinguished minimal}
Let $\pi$ be a Casselman-Wallach representation of $G$. We define \textbf{distinguished minimal $K$-types} $\gamma$ of $\pi$ to be the minimal elements in the set 
\begin{equation} \{ \gamma \in \hat{K}\  | \ \Hom_{K}(\gamma,\pi)\neq 0\ \text{and}\ \Hom_{K_H} (\gamma|_{K_H}, \mathbb{C}) \neq 0 \}.
\end{equation}
\end{definition} 

The set of distinguished minimal $K$-types is non-empty for irreducible $H$-distinguished representations. In general, distinguished minimal $K$-types may not be unique. However, for the $(\mathrm{GL}_n(\mathbb{C}),\mathrm{GL}_n(\mathbb{R}))$ case, we will show that it is unique for both irreducible generic and irreducible unitary representations.

The standard maximal compact subgroups are $K=\mathrm{U}(n)$ and $K'=K_H=\mathrm{O}(n)$ in our case. This pair is a strong Gelfand pair. 
Consequently, the non-vanishing of period on a specific $K$-type is equivalent to specifying a test vector in that $K$-type.

We will need some classical results about branching laws for $\mathrm{U}(n+m)\downarrow \mathrm{U}(n)\times \mathrm{U}(m),$ and $ \mathrm{U}(n)\downarrow \mathrm{O}(n)$ in the proof of Theorem~\ref{induction}. Recall that an irreducible representation of $\mathrm{U}(n)$ is determined by its highest weight, which is a non-increasing $n$-tuple $\mu=(\mu_1,\dots,\mu_n)$ of integers
$$   \mu_1 \geq \mu_2 \geq \dots \geq \mu_n.
$$ From now on, we will not distinguish between an $n$-tuple and an irreducible representation of $\mathrm{U}(n)$. Assume that $\mu$ is an $n$-tuple, and $\gamma$ is an $m$-tuple. Then we use $(\mu,\gamma)$ to denote the $n+m$-tuple which is the non-increasing reordering of $(\mu_1,\dots,\mu_n,\gamma_1,\dots,\gamma_m)$. The branching law, $\mathrm{U}(n+m)\downarrow \mathrm{U}(n)\times \mathrm{U}(m)$, is described by Littlewood-Richardson theorem (cf. \cite[\S 9.6]{Kn96}).  The following lemma can be proved by using the highest weight theory, and we omit the proof.
\begin{lemma}\label{branching-unitary}
    Let $\mu$ be an irreducible representation of $\mathrm{U}(n)$, and $\gamma$ be an irreducible representation of $\mathrm{U}(m)$. Then
    $$
       \dim \Hom_{\mathrm{U}(n)\times \mathrm{U}(m)}((\mu,\gamma)|_{\mathrm{U}(n)\times \mathrm{U}(m)},\mu\boxtimes\gamma)=1.
    $$
\end{lemma}

We note that the category of finite-dimensional representations of $\mathrm{U}(n)$ is equivalent to the category of finite-dimensional holomorphic representations of $\mathrm{GL}_{n}(\mathbb{C})$ through the complexification $\mathrm{U}(n)\hookrightarrow \mathrm{GL}_{n}(\mathbb{C})$. Hence, for every finite-dimensional representation of $\mathrm{U}(n)$, we can associate a representation of $\mathrm{GL}_{n}(\mathbb{R})$ on the same space by restriction, such that:
\begin{itemize}
    \item The $\mathrm{O}(n)$-action is the same;
    \item The highest weight is the same with respect to diagonal Cartan subgroups.
\end{itemize}
Consequently, by Helgason's theorem (cf.  \cite[Theorem 8.49]{Kn96}), we have the following lemma:
\begin{lemma}\label{even}
    Let $\mu=(\mu_1 , \mu_2 , \dots , \mu_n)$ be an irreducible representation of $\mathrm{U}(n)$, then it is $\mathrm{O}(n)$-distinguished if and only if each $\mu_i$ is even for $1\leq i\leq n$.
\end{lemma}

\subsection{Schwartz analysis}
We introduce the Schwartz analysis involved in this article, namely, the Schwartz homology and generic computation of the distribution due to Kobayashi and Speh \cite{KS15}. By Theorem~\ref{main result_unitary}, 
we are mainly concerned with the representations of the form
\[\Pi=\left( (\frac{\det}{|\det|})^k\otimes |\det|^s \right)_{G_n}\times\left( (\frac{\det}{|\det|})^k\otimes |\det|^{-s}\right)_{G_n},\] and \[\pi=|\det|_{G_n'}^{-n/2}\times |\det|_{G_n'}^{n/2},\] for some $k \in \mathbb{Z}$ and $s \in \mathbb{C}$. Note that the trivial representation is a subrepresentation of $\pi$ and hence, $\Hom_{G_{2n}'}(\Pi,\mathbb{C})\hookrightarrow \Hom_{G_{2n}'}(\Pi,\pi)$.

Now we turn to a general setting, where we have a connected reductive Nash group $G$, with a Nash subgroup $G'$, and a Nash subgroup $P$ of $G$, with Nash subgroup $P'=P\cap G'$ of $G'$. We remark that the Nash subgroup is automatically closed. Let $X:=P\backslash G$ and $Y:=P'\backslash G'$.
Furthermore, we have tempered Fr\'echet bundle $\mathcal{V}$ over $X$ and tempered Fr\'echet bundle $\mathcal{W}$ over $Y$. From now on, we assume that:
\begin{center}
    $G'$ has finitely many orbits on $X$.
\end{center}
In particular, when all the groups involved are algebraic, $G'$ will have a unique open dense orbit on $X$. 

\subsubsection{Generic calculation of distribution kernel}
Similar to  \cite[Proposition 3.2]{KS15}, we have the following results:
\begin{proposition}
    Let $W$ be the fiber of $\mathcal{W}$ at $P'e$. It is a Fr\'echet representation of $P'$. 
    \begin{enumerate}
        \item There is a natural injective map:
    \begin{equation}\label{eqdis}
        \Hom_{G'}(\mathcal{S}(X,\mathcal{V}),\mathcal{S}(Y,\mathcal{W}))\hookrightarrow (\mathcal{G}(X,\mathcal{V}^*)\hat{\otimes} W)^{\Delta(P')}, \text{ denoted by } T\mapsto K_T,
    \end{equation}
where $P'$ acts diagonally via the action of $G\times P'$ on $\mathcal{G}(X,\mathcal{V}^*)\hat{\otimes} W$.
    \item If $P$ is cocompact in $G$, then \eqref{eqdis} is an isomorphism.
    \end{enumerate}
\end{proposition} 

\begin{proof}
    (1) By Schwartz kernel theorem \cite[chapter 51]{Tr}, $\Hom_{G'}(\mathcal{S}(X,\mathcal{V}),\mathcal{S}(Y,\mathcal{W}))\hookrightarrow \mathcal{G}(X\times Y,\mathcal{V}^*\boxtimes \mathcal{W})^{\Delta(G')}$. Consider the map 
    \begin{equation}\label{desc}
        G\times G' \longrightarrow G,  \text{ defined by } (g,g')\mapsto g g'^{-1}.
    \end{equation}
    By descending through this map, we get $\mathcal{G}(X\times Y,\mathcal{V}^*\boxtimes \mathcal{W})^{\Delta(G')}\simeq (\mathcal{G}(X,\mathcal{V}^*)\hat{\otimes} W)^{\Delta(P')}$.\\
    The assertion (2) follows directly from \cite[Proposition 3.2]{KS15}.
\end{proof}

Then for any symmetry breaking operator $T$, it has a coarse invariant, the support of distribution kernel $\Supp K_T$. It is a $P'$-invariant closed subset in $X$. Following Kobayashi and Speh, if $\Supp K_T=X$, we call $T$ \textbf{regular}. Until the end of this subsection, we further assume that $\mathcal{W}$ contains the trivial representation of $G'$ with multiplicity one. Then we have:
\begin{equation}\label{eq eb}
  \mathcal{S}(X,\mathcal{V}) '^{G'}= \Hom_{G'}(\mathcal{S}(X,\mathcal{V}),\mathbb{C})\hookrightarrow \Hom_{G'}(\mathcal{S}(X,\mathcal{V}),\mathcal{S}(Y,\mathcal{W}))
\end{equation}
Let the open orbit of $G'$ on $X$ be $U$. Then we have 
$ \mathcal{S}(X,\mathcal{V})' \rightarrow \mathcal{S}(U,\mathcal{V}|_U)' $
given by restriction. Applying $G'$-invariance, we get 
\begin{equation}\label{eq res}
     \mathcal{S}(X,\mathcal{V})'^{G'} \rightarrow \mathcal{S}(U,\mathcal{V}|_U)'^{G'} .
\end{equation} Then we have the  following observation by explicit calculation of Equation~\eqref{desc}:
\begin{lemma}\label{regular}
    If $T\in \mathcal{S}(X,\mathcal{V})'^{G'}$ maps to a nonzero element of $\mathcal{S}(U,\mathcal{V}|_U)'^{G'} $ under \eqref{eq res}, then viewing $T$ as a symmetry breaking operator through \eqref{eq eb}, $T$ is regular.
\end{lemma}
\

Now we specialize to the case that $G'$ is a symmetric subgroup of $G$ fixed by involution $\sigma$, and $P$ is a parabolic subgroup of $G$ with Langlands decomposition $P=MAN$ compatible with $\sigma$. Hence $P'=M'A'N'$, where $M'=M^{\sigma},A'=A^{\sigma},N'=N^{\sigma}$. Let $N_{-}$ be the opposite unipotent subgroup of $P$. The composition $N_{-}\hookrightarrow G\rightarrow X$ is an open dense embedding, such that $\mathcal{V}^*|_{N_{-}}$ has a natural trivialization $N_{-}\times V^*\stackrel{\simeq}{\longrightarrow}\mathcal{V}^*|_{N_{-}}$. Note $M'A'$ still acts on $N_{-}$, and $\mathfrak{n}'$ acts infinitesimally on $N_{-}$ since it is open in $X$. Consequently, by restriction, we have: 
\begin{equation} \label{res}
    (\mathcal{G}(X,\mathcal{V}^*)\hat{\otimes} W)^{\Delta(P')} \longrightarrow (\mathcal{G}(N_{-})\hat{\otimes} V^{\vee} \hat{\otimes} W \otimes \delta_P^2 )^{M'A',\mathfrak{n}'}.
\end{equation}
We denote the image still by $K_T$. Here, $K_T\neq 0$ if $T$ is regular.
\begin{remark}
    In general, particularly in our case, \eqref{res} is not an isomorphism. In some specific cases, like $PN_{-}P'=G$, \cite[Theorem 3.16]{KS15} shows that it is an isomorphism. We call the explicit description of right-hand side the generic calculation of the distribution kernel.
\end{remark}
\subsubsection{Meromorphic continuation}
Assume all the notation  as mentioned above. Now we fix two irreducible representations $\rho$ and $\rho'$ of $M$ and $M'$ respectively. And choose holomorphic families of characters of $A$ and $A'$ depending on parameter $s,t\in\mathbb{C}$, which are denoted by $\nu_s$ and $\nu'_t$. In this subsection, we set:
$$
    V:=\rho\otimes \nu_s ~\text{ and }~ W:=\rho'\otimes \nu'_t.
$$
Observe that for different $s$, $\mathcal{S}(X,\mathcal{V})\simeq \Ind_P^G(\rho\otimes \nu_s) $ is realized as the same space $\Ind_{K\cap M}^K \rho$ by restriction to $K$. The case $\mathcal{S}(Y,\mathcal{W})$ is similar. We will freely use the transformation from $K$-picture to $N_{-}$-picture, which means pulling back functions $i^*:\Ind_{K\cap M}^K \rho \simeq \Ind_P^G(\rho\otimes \nu_s)\rightarrow \mathcal{C}^{\infty}(N_{-})$ through $N_{-}\hookrightarrow G$.

\begin{definition}
        The family $T_{s,t}\in  \Hom_{G'}(\mathcal{S}(X,\mathcal{V}),\mathcal{S}(Y,\mathcal{W}))$ is called \textbf{holomorphically} (resp. \textbf{meromorphically}) dependent on $(s,t)\in\mathbb{C}^2$ if for any $f\in \Ind_{K\cap M}^K \rho$,~ $T_{s,t}(f)$ is \textbf{holomorphically} (resp. \textbf{meromorphically}) dependent on $(s,t)\in\mathbb{C}^2$.
\end{definition}

Although \eqref{res} is not an isomorphism, the following observations tell us that one can still  detect vanishing through it.
\begin{proposition}\label{mero}
   Suppose $\Ind_{P'}^{G'} (\rho'\otimes \nu'_t)$ contains the trivial representation. Let $T_{s}\in \Hom_{G'}(\mathcal{S}(X,\mathcal{V}),\mathbb{C})$ be a holomorphic family of regular periods. Let $f\in\Ind_{K\cap M}^K \rho$ and $\Omega\subset \mathbb{C}$ be a non-empty domain such that ${K_T}_{s}(i^*f)$ is holomorphic on $\Omega$. Suppose ${K_T}_{s}(i^*f)$ can be extended holomorphically to $\mathbb{C}$, then $T_{s}(f)={K_T}_{s}(i^*f)$ for all $s\in\mathbb{C}$.
\end{proposition}
\begin{proof}
    Since ${K_T}_{s}(i^*f)$ is holomorphic on $\Omega$ and $T_s$ is regular, we have $T_{s}(f)={K_T}_{s}(i^*f)$ on $\Omega$. Then the proposition follows by the identity theorem of holomorphic functions.
\end{proof}

\subsubsection{Description of double cosets $P\backslash G/H$ and Schwartz homology}\label{double coset}
We recall the results of Springer \cite{Spr} about the double coset decomposition. Given a connected reductive group $G$ over an algebraically closed field, let $\theta$ be an involution of $G$, and let $H=G^{\theta}$. Let $B$ be a $\theta$-stable Borel subgroup of $G$, and let $T$ be a $\theta$-stable maximal torus of $B$. Then, one has the following result for the double coset $B\setminus G/H$. 

\begin{lemma}\cite[\S 2 and \S 4]{Spr} 
Let $S=\{x\in G \ |\ x\theta(x)=1\}$, and the map from $G/H$ to $S$, sending $x$ to $x\theta(x)^{-1}$ is an embedding as a $G$-variety, where $G$ acts on $S$ by $g.x=gx\theta(g)^{-1}$. Every $B$-orbit in $S$ meets with $N:=N_G(T)$, and for each $w\in N\cap S$, the $B$-orbits in $wT$ are parameterized by $\mathrm{Ker}(w\theta+1)/\mathrm{Im}(w\theta-1)$, which is a finite abelian group of exponent $2$. 
\end{lemma}
Let $W_n = S_n$ be the Weyl group of $G_n$ and $W_{n,2}=\{w \in W_n \mid w^2=1\}$ the set of involutions in the Weyl group. Using the above result of Springer, one can prove the following result.
\begin{lemma}\cite[Lemma 3.1]{Kem15a}
   There is a bijective map between $B_n \backslash G_n \slash G_n'$ and the set of involutions $W_{n,2}$.  
\end{lemma}
Explicitly, $w=1$ corresponds to $g=1$, and the transposition $w=(k,l)$ $(1\leq k<l\leq n)$ corresponds to the 
\[g_w=\begin{pmatrix} I_{k-1} & & & &\\
& 1 & & \sqrt{-1} &\\ & & I_{l-k-1} & & \\ & \sqrt{-1} & & 1 &\\ & & & & I_{n-l}\end{pmatrix} \in G_n.\]
For the disjoint product of simple involutions $w=\prod_{j=1}^s(k_j,l_j)$, it corresponds to $g_w:=\prod_{j=1}^sg_{(k_j,l_j)}$. Let $P$ be the standard parabolic subgroup with Levi subgroup $L=\mathrm{GL}_{n_1}(\mathbb{C})\times\dots\times \mathrm{GL}_{n_k}(\mathbb{C})$, then $P\backslash G_n/G_n'$ is the set of equivalence classes of $W_{n,2}$ under equivalence relation:
\begin{equation*}
    w_1 \sim w_2 \text{ if and only if } w_1=\sigma w_2 \sigma'
\end{equation*}
for some $\sigma,\sigma' \in S_{n_1}\times \dots\times S_{n_k}$. Moreover, one can choose an ordered set of representatives $\{g_j\}_{0\leq j\leq l}$ such that $U_j:=\cup_{i\geq j} Pg_iH$ is an open subset of $G$ for all $0\leq j\leq l$.

We will use Schwartz homology to study the distinguished problem in Lemma \ref{red_int}, Lemma \ref{dist_comp} and Lemma \ref{mult<=1}.
We refer to \cite{CS} for the basic results, see also \cite{ST}.

Let $G=G_n$ and $H=G_n'$. Let $\sigma$ be a Casselman-Wallach representation of $L$, and let $\pi=\mathrm{Ind}_P^G(\sigma)$ be the normalized induction, equivalently, the Schwartz sections of the tempered bundle $\sigma\cdot \delta_P^{1/2}\times_P G$ over $P\backslash G$. Let $\mathrm{H}^{\mathcal{S}}_i(H,\pi)$ be the $i$-th Schwartz homology, so $\mathrm{H}^{\mathcal{S}}_0(H,\pi)'\simeq \mathrm{Hom}_{H}(\pi,\mathbb{C})$.

The restriction of $\pi$ to $H$ has the following filtration which will be used later. Let $\Pi_j$ be the Schwartz sections of the restriction of $\pi$ to the open set $U_j$ (defined by $\cup_{i\geq j} Pg_iH$ as above), then one has a filtration of $\pi$: 
\[\pi=\Pi_0\supset \Pi_1\supset \dots\supset \Pi_l\supset 0.\]
Moreover, by \cite[Proposition 8.3]{CS}, one has a filtration of $\rho_j:=\Pi_j/\Pi_{j+1}$: 
\[\rho_j\supset \rho_{j,1}\supset \rho_{j,2}\supset \dots\] 
such that $\rho_j\simeq \mathop{\varprojlim}\limits_{k}\rho_j/\rho_{j,k}$, and
\[\rho_{j,k}/\rho_{j,k+1}\simeq \mathcal{S}(H/Q_j,H\times_{Q_j}{^{g_j^{-1}}}(\sigma\cdot\delta_P^{1/2})\otimes \mathrm{Sym}^k(\mathfrak{g}/(\mathfrak{h}+\mathrm{Ad}(g_j)\mathfrak{p}))^{\vee}_{\mathbb{C}}).\]
where $Q_j=H\cap g_jPg_j^{-1}$, and ${}^{g_j^{-1}}\sigma(x):=\sigma(Ad(g_j^{-1})x)$ for all $x\in Q_j$. Moreover, we have the following Mittag-Leffler lemma. 
\begin{lemma}[see \cite{Gr}, Chapter 0, Proposition 13.2.3]\label{inverse lim lem}
    Let $\{\gamma_k,\alpha_k:\gamma_{k+1}\to \gamma_k\}_{k\geq 0}$ be an inverse system of smooth moderate growth Fr\'echet representation of $H$. Let $\gamma:=\varprojlim\limits_{k}\gamma_k$. Assume
    \begin{enumerate}
        \item $\{\gamma_k,\alpha_k:\gamma_{k+1}\to \gamma_k\}_{k\geq 0}$ is stationary;
        \item $\{\rmh^{\mathcal{S}}_{i+1}(H,\gamma_k),\alpha_k:\rmh^{\mathcal{S}}_{i+1}(H,\gamma_{k+1})\to \rmh^{\mathcal{S}}_{i+1}(H,\gamma_k)\}_{k\geq 0}$ is also stationary for some integer $i$.
    \end{enumerate}
    Then there is an isomorphism as vector spaces
    \[
   \varprojlim_k\rmh^{\mathcal{S}}_{i}(H,\gamma_k) \simeq \rmh^{\mathcal{S}}_{i}(H,\gamma).
    \]
\end{lemma}
In particular, the inverse system $\rho_j/\rho_{j,k}$ is surjective for any $j$, thus satisfies the stationary condition (1). If one can moreover show that $\rmh^{\mathcal{S}}_{i+1}(H,\rho_j/\rho_{j,k})=0$ for any integer $k$, then this lemma will imply
\begin{equation}\label{lim_hom}
     \varprojlim_k\rmh^{\mathcal{S}}_{i}(H,\rho_j/\rho_{j,k}) \simeq \rmh^{\mathcal{S}}_{i}(H,\rho_j).
\end{equation}
In practice, what concerns us is $\Hom_H(\rho_j,\mathbb{C})$. By the topological isomorphism, $\rho_j\simeq \mathop{\varprojlim}\limits_{k}\rho_j/\rho_{j,k}$, any continuous linear functional on $\rho_j$ must factor through some $\rho_j/\rho_{j,k}$. Thus, we also obtain
\begin{equation}\label{continous eq}
    \Hom_H(\rho_j,\mathbb{C})\simeq \varinjlim_k\Hom_H(\rho_j/\rho_{j,k},\mathbb{C}).
\end{equation}

By Shapiro's lemma, see \cite[Proposition 6.7]{CS},
\begin{equation}\label{shap_lem}\begin{aligned}\mathrm{H}^{\mathcal{S}}_l(H,\rho_{j,k}/\rho_{j,k+1})\simeq \ & \mathrm{H}^{\mathcal{S}}_l(Q_j,{^{g_j^{-1}}}(\sigma\cdot\delta_P^{1/2})\otimes \mathrm{Sym}^k(\mathfrak{g}/(\mathfrak{h}+\mathrm{Ad}(g_j)\mathfrak{p}))^{\vee}_{\mathbb{C}}\cdot \delta_H\cdot \delta_{Q_j}^{-1})\\
\simeq\  & \mathrm{H}^{\mathcal{S}}_l(P^j, \sigma\otimes \mathrm{Sym}^k(\mathfrak{g}/(\mathrm{Ad}(g_j^{-1})\mathfrak{h}+\mathfrak{p}))^{\vee}_{\mathbb{C}}\cdot \delta_{H^j}\cdot \delta_P^{1/2}\cdot\delta_{P^j}^{-1}),\end{aligned}\end{equation}
where $H^j\simeq g_j^{-1}Hg_j$ and $P^j=P\cap H^j$.
We will need some lemmas to compute Schwartz homology.

\begin{lemma}\label{scalar}
    Let $G$ be a Nash group and $\pi$ be a moderate-growth smooth Fr\'echet representation of $G$. Suppose there exists a central element $g \in G$ such that $g$ acts on $\pi$ by a scalar $a \neq 1$. Then $\mathrm{H}^{\mathcal{S}}_i(G, \pi) = 0$ for all $i \in \mathbb{Z}$. 
\end{lemma}

The proof is straightforward and left to the reader.

\begin{lemma}\label{van_hom}
    Let $G$ be a connected real reductive group. Let $\pi$ be a non-trivial irreducible finite dimensional representation of $G$. Then $\mathrm{H}^{\mathcal{S}}_i(G, \pi) = 0$ for all $i \in \mathbb{Z}$.
\end{lemma}
\begin{proof}
    By~\cite[Theorem 7.7]{CS}, $\mathrm{H}^{\mathcal{S}}_i(G, \pi)\simeq \mathrm{H}_i(\mathfrak{g}, K; \pi)$ for any integer $i$. Since $G$ is connected, $\mathrm{H}_i(\mathfrak{g}, K; \pi)\simeq \mathrm{H}_i(\mathfrak{g}, \mathfrak{k}; \pi)$ for any integer $i$. Then the result follows from~\cite[Chapter 1, Theorem 5.3]{BW}.
\end{proof}

In the following two lemmas, let $G = L \ltimes Q$ be a Nash group, where $L$ and $Q$ are Nash subgroups and $Q$ is normal. The next lemma is the Hochschild-Serre spectral sequence for normal subgroups.

\begin{lemma}[\cite{Geng}, Corollary A.10]\label{spectral sequence}
    Let $\pi$ be a moderate-growth smooth Fr\'echet representation of $G$. Assume that $\mathrm{H}^{\mathcal{S}}_i(Q, \pi)$ is Hausdorff. Then there exist convergent first-quadrant spectral sequences:
    \[
    E^2_{p,q} = \mathrm{H}^{\mathcal{S}}_p(L, \mathrm{H}^{\mathcal{S}}_q(Q, \pi)) \Rightarrow \mathrm{H}^{\mathcal{S}}_{p+q}(G, \pi).
    \]
\end{lemma}
We will also need a lemma concerning invariant distributions.
\begin{lemma}\label{invariant dis}
    Let $V$ be a Fr\'echet space, and $D\in \mathcal{S}(G,V)'$. Then $D$ is invariant under right translation of $Q$ if and only if $D=S\otimes dq$ for some $S\in\mathcal{S}(L,V)'$. Here $dq$ is the right invariant measure of $Q$.
\end{lemma}
\begin{proof}
    Under $Q$-action, we have $\mathcal{S}(G,V)\simeq \mathcal{S}(G)\hat{\otimes}V\simeq (\mathcal{S}(L)\hat{\otimes}V)\hat{\otimes} \mathcal{S}(Q) $. Hence the coinvariant
    \begin{equation*}
        \mathcal{S}(G,V)_Q\simeq (\mathcal{S}(L)\hat{\otimes}V)\hat{\otimes} (\mathcal{S}(Q)_Q)\simeq \mathcal{S}(L,V),
    \end{equation*}
    where the second isomorphism follows from \cite[Theorem 1.2]{CS}.
\end{proof}



\section{Necessary condition}\label{Necessary}

\subsection{Reduction to the integral infinitesimal character case}
We begin with a lemma that reduces the proof of necessary part of Theorem \ref{equi_thm} to the case where $\pi$ is a product of unitary characters of the form $\left(\frac{\det}{|\det|}\right)^k_{G_m}$. Assume $n=n_1 + n_2$.

\begin{lemma}\label{red_int}
Let $\pi = \pi_1 \times \pi_2$, where $\pi_1$ is a unitary representation of $G_{n_1}$ given by a product of characters of the form $(\frac{\det}{|\det|})^k$, and $\pi_2$ is a unitary representation of $G_{n_2}$ given by a product of unitary characters $\chi_i:=|\det|^{u_i}(\frac{\det}{|\det|})^{k_i}$ with $u_i \in \sqrt{-1}\mathbb{R}\setminus\{0\}$ and some complementary series. Assume that $\pi_1$ and $\pi_2$ satisfy the condition (i) in Theorem \ref{main result_unitary}. Then 
\[
\dim \mathrm{Hom}_{G_n'}(\pi, \mathbb{C}) \leq \dim\mathrm{Hom}_{G_{n_1}' \times G_{n_2}'}(\pi_1 \boxtimes \pi_2, \mathbb{C}).
\]
\end{lemma}

\begin{proof}
Write $\pi_1 = \mathrm{Ind}_{P_1}^{G_{n_1}}(\mathop{\boxtimes}\limits_{i=1}^r \alpha_i)$ for unitary characters $\alpha_i$ of $G_{a_i}$ and a parabolic subgroup $P_1 \subset G_{n_1}$ with Levi subgroup $\prod_i G_{a_i}$. Similarly, write $\pi_2 = \mathrm{Ind}_{P_2}^{G_{n_2}}(\mathop{\boxtimes}\limits_{j=1}^{s} \beta_j)$ for characters $\beta_j$ of $G_{b_j}$ and a parabolic subgroup $P_2 \subset G_{n_2}$ with Levi subgroup $\prod_j G_{b_j}$. Let $P$ be the parabolic subgroup of $G_n$ with Levi subgroup $\prod_i G_{a_i}\times \prod_j G_{b_j}$. Then $\pi = \mathrm{Ind}_{P}^{G_n} \gamma$ with 
\[
\gamma = \left(\mathop{\boxtimes}\limits_{i=1}^r \alpha_i\right) \boxtimes \left(\mathop{\boxtimes}\limits_{j=1}^s \beta_j\right).
\]
Let $P_0$ be the parabolic subgroup of $G_n$ with Levi subgroup $G_{n_1} \times G_{n_2}$; then $\pi = \mathrm{Ind}_{P_0}^{G_n} (\pi_1 \boxtimes \pi_2)$.  

By the results in Section \ref{double coset}, we set $G=G_n$, $H=G_n'$, and let $\{g_j\}_{1 \leq j \leq m}$ be the representatives of $P \backslash G / H$. Assume that $g_1, \dots, g_l$ satisfy $\bigcup_{1 \leq j \leq l} P g_j H = P_0 H$. Define $U_k$ as the open subset $\bigcup_{j \geq k+l} P g_j H$ for $k>0$, and $U_0 = G$. This gives a filtration:
\[
G = U_0 \supset U_1 \supset \cdots \supset U_{m-l}.
\]
Let $\mathcal{V}$ be the tempered bundle $\gamma \cdot \delta_P^{1/2} \times_P G$. Define $\Pi_j := \mathcal{S}(U_j, \mathcal{V}|_{U_j})$ as the subspace of Schwartz sections of $\mathcal{V}$ restricted to $U_j$, and $\Pi_{m-l+1}=0$.

Let $\rho_j := \Pi_{j}/\Pi_{j+1}$, and set $Q_j := H \cap g_{j+l} P g_{j+l}^{-1}$. The $H$-module $\rho_j$ admits a filtration $\rho_j = \rho_{j,0} \supset \rho_{j,1} \supset \cdots$ such that:
\[
\rho_{j,k}/\rho_{j,k+1} \cong \mathcal{S}\left(H/Q_j, H \times_{Q_j} {^{g_{j+l}^{-1}}}(\gamma \cdot \delta_P^{1/2}) \otimes \mathrm{Sym}^k\left(\mathfrak{g}/(\mathfrak{h} + \mathrm{Ad}(g_{j+l})\mathfrak{p})\right)_{\mathbb{C}}^{\vee}\right), \quad \forall k \in \mathbb{Z}_{\geq 0}.
\]
By \eqref{shap_lem},
\[
\mathrm{Hom}_H(\rho_{j,k}/\rho_{j,k+1}, \mathbb{C}) \cong \mathrm{Hom}_{P^j}\left(\gamma \otimes \mathrm{Sym}^k\left(\mathfrak{g}/(\mathfrak{h}^j + \mathfrak{p})\right)_{\mathbb{C}}^{\vee} \otimes \delta_{H^j} \cdot \delta_P^{1/2} \cdot \delta_{P^j}^{-1}, \mathbb{C}\right),
\]
where $P^j = P \cap H^j$ and $H^j = g_{j+l}^{-1} H g_{j+l}$.

Using the exact sequences:
\[
0 \to \Pi_{j+1} \to \Pi_j \to \rho_j \to 0, \quad 0 \leq j \leq m-l,
\]
to show $\dim \mathrm{Hom}_H(\Pi, \mathbb{C}) \leq \dim\mathrm{Hom}_{G_{n_1}' \times G_{n_2}'}(\pi_1 \boxtimes \pi_2, \mathbb{C})$, it suffices by \eqref{lim_hom} to prove $\mathrm{Hom}_H(\rho_{j,k}, \mathbb{C}) = 0$ for all $j \geq 1,k\in\mathbb{Z}_{\geq 0}$, and  
\[
\dim \mathrm{Hom}_{H}(\rho_0, \mathbb{C}) \leq \dim\mathrm{Hom}_{G_{n_1}' \times G_{n_2}'}(\pi_1 \boxtimes \pi_2, \mathbb{C}).
\]

\textbf{Step 1:}  $\mathrm{Hom}_H(\rho_{j,k}, \mathbb{C}) = 0$ for $j \geq 1,k\in\mathbb{Z}_{\geq 0}$.  
Assume that $g_{j+l}$ corresponds to $w \in W_{n,2}$. Then there exist indices $1 \leq k_1 \leq n_1$ and $1 \leq k_2 \leq n_2$ with $w(k_1) = n_1 + k_2$. The subgroup
\[
\left\{t(x) := \begin{pmatrix} 
I_{k_1-1} & & & & \\ 
& x &  & & \\ 
& & I_{n_1 - k_1 + k_2 - 1} & & \\ 
& & & x & \\ 
& & & & I_{n_2 - k_2}
\end{pmatrix} \ \bigg\vert\ x \in \mathbb{R}_{>0} \right\}
\]
lies in $P^j$ and acts diagonally on $\gamma \otimes \mathrm{Sym}^k\left(\mathfrak{g}/(\mathfrak{h}^j + \mathfrak{p})\right)_{\mathbb{C}}^{\vee} \otimes \delta_{H^j} \cdot \delta_P^{1/2} \cdot \delta_{P^j}^{-1}$. The eigenvalues of $t(x)$ on $\mathrm{Sym}^k(\cdots)$ are of the form $x^m$ for $m \in \mathbb{Z}$. By assumption on $\gamma$, $t(x)$ acts on $\gamma$ via $x^s$ with $s \notin \mathbb{Z}$, so there is no non-zero fixed element under the action of $t(x)$. Hence, $\mathrm{Hom}_H(\rho_{j,k}, \mathbb{C}) = 0$.

\textbf{Step 2:}  Bound $\dim \mathrm{Hom}_H(\rho_0, \mathbb{C})$.  
The space $\Pi_1$ consists of Schwartz sections of $\gamma \cdot \delta_P^{1/2} \times_{P} G$ over $U_1$. Let $L_0$ be the Levi subgroup of $P_0$. The bundle map
\[
\gamma \cdot \delta_{P}^{1/2} \times_{P} G \to \left(\gamma \cdot \delta_{P \cap L_0}^{1/2} \times_{P \cap L_0} L_0\right) \cdot \delta_{P_0}^{1/2} \times_{P_0} G \to (\pi_1 \boxtimes \pi_2) \cdot \delta_{P_0}^{1/2} \times_{P_0} G
\]
identifies $\Pi_1$ with Schwartz sections of $(\pi_1 \boxtimes \pi_2) \times_{P_0} G$ over $U_1$. The quotient $\rho_0 := \Pi / \Pi_1$ has a filtration $\rho_0 \supset \rho_{0,1} \supset \cdots$ with:
\[
\rho_{0,k}/\rho_{0,k+1} \simeq \mathcal{S}\left(H/(H \cap P_0), H \times_{H \cap P_0} \left((\pi_1\boxtimes \pi_2) \cdot \delta_{P_0}^{1/2}\right) \otimes \mathrm{Sym}^k\left(\mathfrak{g}/(\mathfrak{h} + \mathfrak{p}_0)\right)_{\mathbb{C}}^{\vee}\right).
\]
Thus,
\[
\mathrm{Hom}_H(\rho_{0,k}/\rho_{0,k+1}, \mathbb{C}) \simeq \mathrm{Hom}_{G_{n_1}' \times G_{n_2}'}\left((\pi_1 \boxtimes \pi_2) \otimes \mathrm{Sym}^k(\cdots) \otimes \delta_{P_0}^{1/2} \cdot \delta_H \cdot \delta_{P_0 \cap H}^{-1}, \mathbb{C}\right).
\]
The subgroup $\left\{\begin{pmatrix} a_1 I_{n_1} & \\ & a_2 I_{n_2} \end{pmatrix} \ \big\vert\ a_1, a_2 \in \mathbb{R}^\times \right\} \subset P_0 \cap H$ acts trivially on $\pi_1 \boxtimes \pi_2$ since $\pi_1,\pi_2$ satisfy condition (i) of Theorem \ref{main result_unitary}. As $\delta_{P_0}^{1/2} = \delta_{P_0 \cap H}$, $\delta_H = 1$, and the action of the subgroup on $\mathrm{Sym}^k(\cdots)$ is semisimple with no non-zero fixed element for $k > 0$, we conclude $\mathrm{Hom}_H(\rho_{0,k}/\rho_{0,k+1}, \mathbb{C}) \neq 0$ only when $k = 0$. In this case,
\[
\mathrm{Hom}_H(\rho_{0,0}/\rho_{0,1}, \mathbb{C}) \simeq \mathrm{Hom}_{G_{n_1}' \times G_{n_2}'}(\pi_1 \boxtimes \pi_2, \mathbb{C}).
\]
This completes the proof.
\end{proof}

\subsection{Distinction under derivative}\label{derivative section}
First, we establish a lemma that will be utilized in the proof of Theorem~\ref{deri_thm}.
\begin{lemma}\label{key lemma}
    Let $\pi$ be an irreducible unitary representation of $G_{n}$ with depth $d \geq 2$. Let $\sigma$ denote the transpose of the associated partition of the associated variety of $\mathrm{Ann}(\pi)$. Suppose $\pi$ satisfies condition (i) in Theorem~\ref{equi_thm}, while $\pi^-$ is not distinguished. Then we have:
    \begin{equation*}
        \mathrm{H}^{\mathcal{S}}_0(H_{n,d}', (\pi^- \boxtimes \psi_{1}^{\sigma}) \otimes |\det|^{-1} \otimes \mathrm{Sym}^k(\mathfrak{{p}}_n / (\mathfrak{h}_{n,d} + \mathfrak{p}_n'))^{\vee}_{\mathbb{C}} \cdot \delta_{P_n'} \cdot \delta_{H_{n,d}}^{1/2} \cdot \delta_{H_{n,d}'}^{-1}) =0
\end{equation*}
for all $k\in\mathbb{Z}_{\geq 0}$.
\end{lemma}

\begin{proof}
    By direct computation, as $H'_{n,d}$-representations, the terms $\psi_1^{\sigma}$ and $|\det|^{-1} \cdot \delta_{P_n'} \cdot \delta_{H_{n,d}}^{1/2} \cdot \delta_{H_{n,d}'}^{-1}$ are trivial. Therefore, for all $i \in \mathbb{Z}$ and $k \in \mathbb{Z}_{\geq 0}$,
    \begin{equation*}
        \begin{aligned}
          &  \mathrm{H}^{\mathcal{S}}_i(H_{n,d}', (\pi^- \boxtimes \psi_{1}^{\sigma}) \otimes |\det|^{-1} \otimes \mathrm{Sym}^k(\mathfrak{{p}}_n / (\mathfrak{h}_{n,d} + \mathfrak{p}_n'))^{\vee}_{\mathbb{C}} \cdot \delta_{P_n'} \cdot \delta_{H_{n,d}}^{1/2} \cdot \delta_{H_{n,d}'}^{-1}) \\
         & \simeq  \mathrm{H}^{\mathcal{S}}_i(H_{n,d}', \pi^- \otimes \mathrm{Sym}^k(\mathfrak{{p}}_n / (\mathfrak{h}_{n,d} + \mathfrak{p}_n'))^{\vee}_{\mathbb{C}}).
         \end{aligned}
    \end{equation*}
    Consider the Levi decomposition $H_{n,d}' = G_{n-d}' \ltimes Q$, where $Q$ is the unipotent radical of $H_{n,d}'$. Recall that by the definition of the Schwartz homology, for arbitrary representation $\gamma$ of $H_{n,d}'$, 
    \[   \mathrm{H}^{\mathcal{S}}_0 (H_{n,d}' , \gamma)\simeq \gamma_{H_{n,d}'}=  (\gamma_{Q})_{G_{n-d}'}\simeq  \mathrm{H}^{\mathcal{S}}_0(G_{n-d}',  \mathrm{H}^{\mathcal{S}}_0(Q,\gamma)).
    \]
    We first observe:
    \begin{equation*}
        \mathrm{H}^{\mathcal{S}}_0(Q, \pi^- \otimes \mathrm{Sym}^k(\mathfrak{{p}}_n / (\mathfrak{h}_{n,d} + \mathfrak{p}_n'))^{\vee}_{\mathbb{C}}) \simeq \mathrm{H}^{\mathcal{S}}_0(Q,  \mathrm{Sym}^k(\mathfrak{{p}}_n / (\mathfrak{h}_{n,d} + \mathfrak{p}_n'))^{\vee}_{\mathbb{C}}) \otimes \pi^-.
    \end{equation*}

    \textit{Case (a).} When $k = 0$, we have $\mathrm{H}^{\mathcal{S}}_0(H_{n,d}',\pi^{-})\simeq \mathrm{H}^{\mathcal{S}}_0(G_{n-d}',\pi^{-})=0$, since $\pi^{-}$ is not distinguished.\\

    \textit{Case (b).} When $k \in \mathbb{Z}_{>0}$, let $\mathfrak{u}_d$ denote the lower triangular subspace of $\mathfrak{g}_{d-1}$ with purely imaginary entries. Embed $\mathfrak{g}_{d-1}$ into the bottom-right block of $\mathfrak{g}_{n-1}$, so that $\mathfrak{u}_d$ embeds as a subspace of $\mathfrak{g}_n$ via the top-left embedding $\mathfrak{g}_{n-1}\hookrightarrow\mathfrak{g}_n$. Note that $\mathfrak{u}_d$ is a nilpotent $H_{n,d}'$-subrepresentation of $\mathfrak{{p}}_n / (\mathfrak{h}_{n,d} + \mathfrak{p}_n')$. For simplicity, we assume that $\mathfrak{u}_d$ is trivial; alternatively, one may consider a filtration where the successive quotients are trivial. Consider the filtration:
    \begin{equation*}
       \mathrm{Sym}^k(\mathfrak{{p}}_n / (\mathfrak{h}_{n,d} + \mathfrak{p}_n'))_{\mathbb{C}}=M_0 \supset M_1 \supset  \dots\supset M_k \supset 0=M_{k+1},
    \end{equation*} 
    where $M_j$ consists of the symmetric product with at least $j$ elements in $\mathfrak{u}_d$. Let $\rho_j:=M_{j}/M_{j+1}$ denote the successive quotient. 

    To prove case (b), it suffices to show the following for $0\leq j\leq k-1$:
    \begin{equation*}
      \text{(I). }   \mathrm{H}^{\mathcal{S}}_0(H_{n,d}',\pi^{-}\otimes\rho_k^{\vee})\simeq \mathrm{H}^{\mathcal{S}}_0(H_{n,d}',\pi^{-}\otimes \mathrm{Sym}^k(\mathfrak{u}_d)_{\mathbb{C}}^{\vee})=0, \text{ and } \text{(II). }\mathrm{H}^{\mathcal{S}}_0(H_{n,d}',\pi^{-}\otimes \rho_j)=0.
    \end{equation*}
    
    \textit{Part (I).} $\mathrm{H}^{\mathcal{S}}_0(H_{n,d}',\pi^{-}\otimes\rho_k^{\vee})\simeq\mathrm{H}^{\mathcal{S}}_0(G_{n-d}',\pi^{-})\otimes\rho_k^{\vee}=0$
    
    \textit{Part (II).} Consider the central diagonal element $g$ of $G_{n-d}'$ with entries $a\in\mathbb{R}_{>1}$. Note that $g$ acts on $\pi^{-}$ by the scalar one, since $\pi$ satisfies condition (i). On the other hand, $g$ acts on $\rho_j$ by the scalar $a^{k-j}$, and thus on $\mathrm{H}^{\mathcal{S}}_0(Q,\rho_j)$ by the scalar $a^{k-j}\neq 1$ as well. Therefore, by Lemma~\ref{scalar}, 
 \[
 \mathrm{H}^{\mathcal{S}}_0(H_{n,d}',\pi^{-}\otimes \rho_j)\simeq \mathrm{H}^{\mathcal{S}}_0(G_{n-d}',\pi^{-}\otimes \mathrm{H}^{\mathcal{S}}_0(Q,\rho_j))=0.
 \]
\end{proof}

The main result of this subsection is the following, which is inspired by \cite{Ka}:
\begin{theorem}\label{deri_thm}
    Let $\pi$ be an irreducible unitary representation of $G_{n}$ with depth $d$. Suppose $\pi$ satisfies condition (i) in Theorem~\ref{equi_thm}, while $\pi^-$ is not distinguished. Then $\pi$ is not distinguished.
\end{theorem}

\begin{proof}
    If $d = 1$, then $\pi$ is a unitary character, and the statement is trivial. For simplicity, we assume $d = 2$. The case $d > 2$ is analogous. Suppose, on the contrary, that $\pi^-$ is not distinguished but $\pi$ is distinguished. Let $\kappa \in \Hom_{P_{n}'}(\pi, \mathbb{C})$ be a non-zero functional. By Theorem~\ref{BZ-fil}, there is an exact sequence:
    \[
   \Hom_{P_n'} (D^1(\pi),\mathbb{C}) \longrightarrow \Hom_{P_n'} (\pi|_{P_n},\mathbb{C}) \stackrel{\vartheta}{\longrightarrow} \Hom_{P_n'} (\SInd_{H_n}^{P_{n}}(\Psi(\pi) \boxtimes \psi_{n}),\mathbb{C}) \longrightarrow0
    \]
    By Proposition~\ref{unitary characterization}, we obtain
    \[
    \dim \Hom_{P_n'} (D^1(\pi),\mathbb{C} )\leq \dim \Hom_{G_{n-1}'} (D^1(\pi),\mathbb{C})=0
    \]
    since $\Hom_{G_{n-1}'} (F_1^0(\pi)/F_1^l(\pi),\mathbb{C})=0$ for any non-negative integer $l$. Thus, $\vartheta$ is surjective. Consider the composition $D_{\kappa}$:
    \[
    \mathcal{S}(P_{n}, \Psi(\pi)) \stackrel{av}{\longrightarrow} \SInd_{H_n}^{P_{n}}(\Psi(\pi) \boxtimes \psi_{n}) \stackrel{\kappa \circ \varphi}{\longrightarrow} \mathbb{C},
    \]
    where $av$ is defined by
    \begin{equation*}
        av(f)(g) := \int_{H_n} |\det h|^{-1} \sigma(h)^{-1} f(hg) \, dh
    \end{equation*}
    for $h \in H_n$, $g \in P_{n}$, and $\sigma := \Psi(\pi) \cdot \psi_{n}$. Here, $av$ is a surjection; see \cite[\S 6.3]{CS}. Note that $D_{\kappa}$ is non-zero since $\vartheta$ is surjective. The $\sigma$-tempered distribution $D_{\kappa}$ satisfies the following:
    \begin{align}
        & D_{\kappa} \circ R(g_0) = D_{\kappa} \quad \text{for all } g_0 \in P_{n}', \label{eq5} \\
        & D_{\kappa} \circ L(h) = |\det h|^3 D_{\kappa} \circ \sigma(h)^{-1} \quad \text{for all } h \in H_n, \label{eq6}
    \end{align}
    where $R$ and $L$ denote the right and left translations of $P_{n}$ on $\mathcal{S}(P_{n}, \Psi(\pi))$, respectively.

    \begin{claim}
        $D_{\kappa}$ is supported on $H_n \cdot G_{n-1}'$.
    \end{claim}

    We view $\psi_{n}$ as a smooth function on $P_{n}$ via the canonical projection $P_{n} \rightarrow V_{n}$. Then, by Equation~\eqref{eq6}, $\psi_{n} \cdot D_{\kappa}$ is left invariant under $V_{n}$. By Lemma~\ref{invariant dis},
    \begin{equation}\label{tensor_inv}
        \psi_{n} \cdot D_{\kappa} = dv \otimes S,
    \end{equation}
    for some $S$, a $\sigma$-tempered distribution on $G_{n-1}$. Now, suppose $g \in \Supp(S)$. Since the right-hand side of \eqref{tensor_inv} is $V_{n}'$-right invariant, by Equation~\eqref{eq5}, we have $\psi_{n}(gv) = \psi_{n}(g)$ for all $v \in V_{n}'$. Hence, $g \in P_{n-1} \cdot G_{n-1}'$ by the definition of $\psi_{n}$. Therefore, $D_{\kappa}$ is supported on $P_{n-1} \cdot G_{n-1}' \cdot V_{n} = H_n \cdot G_{n-1}'$.

    By the theory in the Preliminary Section~\ref{double coset}, $\SInd_{H_n}^{P_{n}}(\Psi(\pi) \boxtimes \psi_{n})$ admits a $P_n'$-subrepresentation consisting of Schwartz sections restricted to the open set $P_n \setminus H_{n} \cdot P_n'$. Denote the quotient by $\rho$, which possesses a decreasing filtration $\rho = \rho_0 \supset \rho_1 \supset \dots$ by Borel's lemma. Note that $D_{\kappa}$ factors through to a non-zero linear functional on $\rho$ since it is supported on $H_n \cdot G_{n-1}'$. However, by Lemma~\ref{key lemma}, we have $\mathrm{H}^{\mathcal{S}}_0(P_n', \rho_k / \rho_{k+1}) = 0$ for all $k \in \mathbb{Z}_{\geq 0}$. Consequently, $\mathrm{H}^{\mathcal{S}}_0(P_n', \rho / \rho_{k}) = 0$ for all $k \in \mathbb{Z}_{\geq 0}$, which implies
    \begin{equation*}
        \Hom_{P_n'}(\rho,\mathbb{C})=0.
    \end{equation*}
    This leads to a contradiction.
\end{proof}

\begin{conjecture}\label{natural embedding}
    Suppose $\pi$ is an irreducible $G_n'$-distinguished representation of $G_n$ with depth $d$, then there is a natural embedding:
    \begin{equation*}
        \Hom_{G_{n}'}(\pi,\mathbb{C})\hookrightarrow \Hom_{G_{n-d}'}(\pi^{-},\mathbb{C}).
    \end{equation*}
    This implies that if $\pi$ is distinguished, then $\pi^-$ is also distinguished.
\end{conjecture}

\subsection{Proof of condition (ii) of Theorem \ref{equi_thm}}

Let us prove the necessary condition (ii) of Theorem \ref{equi_thm}, that is, all irreducible unitary $G'_n$-distinguished representations satisfy condition (ii) of Theorem \ref{equi_thm}. 

\begin{proof}
Let $\pi$ be $G'_n$-distinguished. If possible, we assume that $\pi$ does not satisfy the condition (ii). By Lemma \ref{red_int}, one can reduce the problem to the case where $\pi=\prod_{i=1}^m(\frac{\det}{|\det|})^{k_i}_{G_{n_i}}$. Assume that $n_{i_0}$ is the largest integer such that $k_{i_0}\in 2 \mathbb{Z}+1$ and the multiplicity of  $(\frac{\det}{|\det|})^{k_{i_0}}_{G_{n_{i_0}}}$ is odd. Applying the highest derivative $n_{i_0} - 1$ times to $\pi$, we obtain an irreducible unitary representation 
\[
\begin{aligned}
\pi^{(n_{i_0}-1)} := & \left(\frac{\det}{|\det|}\right)^{k_1} \times \dots \times \left(\frac{\det}{|\det|}\right)^{k_1} \times \dots \times \left(\frac{\det}{|\det|}\right)^{k_{i_0}} \times \dots \times \left(\frac{\det}{|\det|}\right)^{k_{i_0}} \times \dots \\
& \times \left(\frac{\det}{|\det|}\right)^{k_m} \times \dots \times \left(\frac{\det}{|\det|}\right)^{k_m}
\end{aligned}
\]
of $G_{r}$ for some $r$, where $(\frac{\det}{|\det|})^{k_i}$ is a character of $G_{n_i-n_{i_0}+1}$ when $n_i\geq n_{i_0}$, and is omitted when $n_i<n_{i_0}$. Note that if $k_j\in 2\mathbb{Z}+1$, then $n_j-n_{i_0}+1\leq 1$ by the assumption on $n_{i_0}$. The representation $\pi^{(n_{i_0}-1)}$ does not satisfy the condition (i) of Theorem \ref{equi_thm} as there does not exist $w \in W_{r,2}$ satisfying the condition (i) of Theorem \ref{equi_thm}, since the Langlands parameter of $\pi^{(n_{i_0}-1)}$ contains $z\mapsto (\frac{z}{|z|})^{k_{i_0}}$ ($k_{i_0}\in 2\mathbb{Z}+1$) with odd multiplicity. So, $\pi^{(n_{i_0}-1)}$ is not $G_{r}'$-distinguished. This contradicts Theorem \ref{deri_thm}, and hence, the statement follows.

\end{proof}

\subsection{Period integral through degenerate Whittaker functions}
Assume that $\pi$ is an irreducible unitary distinguished representation of $G_n$ with depth $d$. Let $\sigma=(\sigma_1,\dots,\sigma_r)$ denote the transpose of the associated partition of the associated variety of $\mathrm{Ann}(\pi)$. We will define an inverse map $\Hom_{G_{n-d}'}(\pi^{-},\mathbb{C})\rightarrow \Hom_{P_n'}(\pi,\mathbb{C})$ through integral of \textbf{degenerate Whittaker functions}. This fact also provides strong support for Conjecture~\ref{natural embedding}.

We know $\pi$ embeds as a unique subrepresentation of $ \Ind_{H_{n,d}}^{G_n}((\pi^{-}\boxtimes \psi^{\sigma}_1)\otimes|\det|^{-1})$. We fix such an embedding. For $v\in V_{\pi}$, we denote its image by $W_v$. It is a smooth function from $G_n$ to $V_{\pi^{-}}$, satisfying
\begin{equation}\label{eq12}
    W_v(h(a,n)g)=|\det a|^{d-1}\psi^{\sigma}_1(n)\pi^{-}(a)\cdot W_v(g)
\end{equation}
for $a\in G_{n-d}$ and $n\in N_{d}$. The function space of such $W_v$ is denoted by $\mathcal{W}_1(\pi,\psi)$. Similarly, applying the highest derivative to $\pi^{-}$ and so on, we get a family of function spaces $\mathcal{W}_j(\pi,\psi),1\leq j\leq r$. Such function spaces are all unique, since the degenerate Whittaker functional of $\pi$ is unique by Theorem~\ref{1.1} (4).
Let $\beta\in \Hom_{G_{n-d}'}(\pi^{-},\mathbb{C})$, we define the integral formally as
\begin{equation}\label{degenerate}
    \Lambda(W_v):=\int\limits_{H_{n,d}'\backslash P_n' } \beta( W_v(p))dp.
\end{equation}
We show that it is absolutely convergent using the following asymptotic behavior of the degenerate Whittaker function. The result for functions in $\mathcal{W}_j(\pi,\psi)$ follows similarly.
\begin{theorem}\label{thm7.2}
    Let $(\pi,V_{\pi})$ be a distinguished irreducible unitary representation of $G_n$ for $n\geq 2$. Assume that $\pi$ is of depth $d\geq 2$. Then there exists a finite set $C$ of characters of $(\mathbb{C}^{\times})^{d-1}$ with positive real parts, a positive integer $n_{\chi}$ for each $\chi\in C$, and for each $v\in V_{\pi}$, there exists a finite set of Schwartz functions $\{\phi_x\in \mathcal{S}(\mathbb{C}^{d-1}, V_{\pi^{-}})\}_{x\in X}$ indexed by the finite set $X$ of functions on $(\mathbb{C}^{\times})^{d-1}$ of the form $\chi(a)\cdot (\log|a|)^n$, where $\chi\in C$ and $0\leq n\leq n_{\chi}$ such that
    \begin{equation}
        W_v(a)=\delta^{1/2}_{d-1}(a) \sum_{x\in X} x(a_{n-d+1},\dots,a_{n-1})\phi_x(a_{n-d+1},\dots,a_{n-1}),
    \end{equation}
    where $a=m(a_{n-d+1},\dots, a_{n-1})$.
\end{theorem}
\begin{proof}
Since $\pi$ is unitary, we have:
\begin{equation*}
    \int\limits_{H_{n-1,d}\backslash G_{n-1}} \Vert W_v(g)\Vert^2 dg 
     =\int\limits_{H_{n,d}\backslash P_n} \Vert W_v(p)\Vert^2 dp < +\infty.
\end{equation*}
Thanks to Iwasawa decomposition, we have:
\begin{equation*}
    \int\limits_{A_{n-d+1}\dots A_{n-1}} \Vert W_v(a)\Vert^2 \delta_{d-1}^{-1}(a)d^{\times}a_{n-d+1}\dots d^{\times}a_{n-1} < +\infty,
\end{equation*}
where $a=m(a_{n-d+1},\dots, a_{n-1})$.
Hence the result follows from Proposition~\ref{prop8.2}, Proposition~\ref{prop8.3} and Proposition~\ref{prop8.4}.
\end{proof}
\begin{remark}
    When $\pi$ is generic (not necessarily distinguished), the result was established by \cite{JS1}.
\end{remark}
\subsubsection{The absolute convergence of the integral~\eqref{degenerate}:}
    First of all, by Equation~\eqref{eq12}, the integral is well-defined. If $d=1$, then the result is obvious. For $d\geq 2$, we have
\begin{equation*}
    \Lambda(W_v)=\int\limits_{H_{n-1,d}'\backslash G_{n-1}'} \beta( W_v(g))dg.
\end{equation*}
Thanks to Iwasawa decomposition, it is sufficient to show the following integral is absolutely convergent
\begin{equation*}
    \int\limits_{A'_{n-d+1}\dots A'_{n-1}}\beta( W_v(a))\delta^{-1/2}_{d-1}(a) d^{\times}a_{n-d+1}\dots d^{\times}a_{n-1},
\end{equation*}
where $a=m(a_{n-d+1},\dots, a_{n-1})$. Then it follows from Theorem~\ref{thm7.2} directly.
\quad  \quad \quad \qedsymbol{}

\begin{remark}
The integral~\eqref{degenerate} is also non-zero. This is because $\pi|_{P_n'}$ contains a subrepresentation $c$-$\Ind_{H_{n,d}'}^{P_n'}(\pi^{-}\boxtimes \mathbb{1})$.     
\end{remark}

\subsubsection{ }Now we use integral~\eqref{degenerate} to define some invariants of $\pi$. Let $w_n$ be the representative of the longest Weyl element in $\mathrm{GL}_n$
given by the matrix whose entries are 1 on the anti-diagonal and zeros otherwise. Let $w_{n,d}:=\begin{bmatrix}
    I_{n-d} & \\
         & w_{d}
\end{bmatrix}$. For $W \in \mathcal{W}_j(\pi, \psi)$, define $W^{\sharp}(g) := W(w_{n,d_j}(g^{-t}))$, where $d_j = \sigma_1 + \dots + \sigma_j$. Using the MVW involution, we realize $\pi^{\vee}$ as a representation on the space $V_{\pi}$ with the action given by $\pi^{\vee}(g) = \pi(g^{-t})$. Hence, $W^{\sharp} \in \mathcal{W}_j(\pi^{\vee}, \psi^{-1})$. Applying the same convergence argument to $\pi^{\vee}$, we observe that the map $v \mapsto \Lambda(W_v^{\sharp})$ is $P_n'$-invariant, and thus $G_n'$-invariant for all $v \in \pi^{\vee}$. Consequently, the map $v \mapsto \Lambda(W_v^{\sharp})$ is also $G_n'$-invariant for all $v \in \pi$. We obtain
\[
\Lambda(W_v^{\sharp}) = c \cdot \Lambda(W_v)
\]
for all $v \in V$. Note that $c$ is independent of the choice of $\beta$ and is, in fact, an invariant of $\pi$. We denote $c$ by $c_j(\pi)$. Since $(W^{\sharp})^{\sharp} = W$, we have $c_j(\pi) \cdot c_j(\pi^{\vee}) = 1$. Moreover, by condition (i) of Theorem 3.1, we have $\overline{\pi} \simeq \pi^{\vee}$, which implies $c_j(\pi^{\vee}) = c_j(\overline{\pi}) = c_j(\pi)$. Consequently, $c_j(\pi) = \pm 1$. 

\section{Sufficient condition}\label{sufficient}
\subsection{Distinction under parabolic induction}\label{parabolic induced}
Let $\pi_1$ (resp.\ $\pi_2$)  be representations of $G_n$ (resp.\ $G_m$), we have a bilinear map:
\[
\Hom_{G_n'}(\pi_1, \mathbb{C}) \times \Hom_{G_m'}(\pi_2, \mathbb{C}) \longrightarrow \Hom_{G_{n+m}'}(\pi_1 \times \pi_2, \mathbb{C})
\]
mapping $(\mathcal{P}_1,\mathcal{P}_2)$, the $G_n'$ and $G_m'$-periods of $\pi_1$ and $\pi_2$ respectively, to the composition of the following maps:
\begin{equation*}
    \pi_1 \times \pi_2 \xrightarrow{\mathrm{res}} \left(\pi_1 \otimes |\det|^{\frac{m}{2}}\right)_{G_n'} \times \left(\pi_2 \otimes |\det|^{\frac{-n}{2}}\right)_{G_m'} \xrightarrow{\mathcal{P}_1 \boxtimes \mathcal{P}_2} \left(|\det|^{\frac{m}{2}}\right)_{G_n'} \times \left(|\det|^{\frac{-n}{2}}\right)_{G_m'} \xrightarrow{\mathrm{pr}} \mathbb{C},
\end{equation*}
where the first map restricts functions on $G_{n+m}$ to $G_{n+m}'$, the second map applies $\mathcal{P}_1 \boxtimes \mathcal{P}_2$, and the third map is the canonical quotient.
This induces a linear map:
\[
I: \Hom_{G_n'}(\pi_1, \mathbb{C}) \otimes \Hom_{G_m'}(\pi_2, \mathbb{C}) \longrightarrow \Hom_{G_{n+m}'}(\pi_1 \times \pi_2, \mathbb{C}).
\]
 We denote $I(\mathcal{P}_1 \otimes \mathcal{P}_2)$ by $\mathcal{P}_1 \times \mathcal{P}_2$ and call it the \textbf{parabolic induced period} of $\mathcal{P}_1$ and $\mathcal{P}_2$. Let $N$ be the opposite unipotent radical of the standard parabolic subgroup $P$ with Levi $G_n \times G_m$. Then, via restriction to the trivialization $N \hookrightarrow P \backslash G$, we conclude that $I$ is injective.
\begin{theorem}\label{induction}
    Let $(\pi_1,V_1)$ (resp.\ $(\pi_2,V_2)$) be $G'_n$-distinguished (resp. $G'_m$-distinguished) representation of $G_n$ (resp.\ $G_m$) with period $\mathcal{P}_1$ (resp.\ $\mathcal{P}_2$), which is non-vanishing on the $K$-type $\gamma_1$ (resp.\ $\gamma_2$). Then the period $\mathcal{P}_1 \times \mathcal{P}_2$ is non-vanishing on the $K$-type $(\gamma_1, \gamma_2)$.
\end{theorem}

\begin{proof}
Let $\gamma = (\gamma_1, \gamma_2)$. By definition, the $\gamma$-isotypic space in $\pi_1 \times \pi_2$ consists of functions $f \in \mathcal{C}^{\infty}(G_{n+m}, V_1 \otimes V_2)$ such that
\[
f|_{\mathrm{U}(n+m)}(k) = u(k \cdot v), \quad k \in \mathrm{U}(n+m),
\]
for some $u \in \mathrm{Hom}_{\mathrm{U}(n) \times \mathrm{U}(m)}(\gamma, \gamma_1 \boxtimes \gamma_2)$ and $v \in \gamma$. Take $v$ to be the $\mathrm{O}(n+m)$-fixed vector in $\gamma$. Then the image of $f$ in $\left(|\det|^{\frac{m}{2}}\right)_{G_n'} \times \left(|\det|^{-\frac{n}{2}}\right)_{G_m'}$ under $(\mathcal{P}_1 \times \mathcal{P}_2) \circ \mathrm{res}$ is a function $h \in \mathcal{C}^{\infty}(G_{n+m}', \mathbb{C})$ such that
\[
h|_{\mathrm{O}(n+m)}(k') = \mathcal{P}_1 \boxtimes \mathcal{P}_2(u(k' \cdot v)), \quad k' \in \mathrm{O}(n+m).
\]
Notice that $u(v)$ is fixed by $\mathrm{O}(n) \times \mathrm{O}(m)$. Hence, if $u(v) \neq 0$, then $u(v)$ has a non-zero image under $\mathcal{P}_1 \boxtimes \mathcal{P}_2$ by the assumption on $\mathcal{P}_1$ and $\mathcal{P}_2$. Meanwhile, if $u(v) \neq 0$, the function $u(k' \cdot v)$, $k' \in \mathrm{O}(n+m)$, lies in the trivial $\mathrm{O}(n+m)$-type since $v$ is an $\mathrm{O}(n+m)$-fixed vector. Hence, its image in $\left(|\det|^{\frac{m}{2}}\right)_{G_n'} \times \left(|\det|^{-\frac{n}{2}}\right)_{G_m'}$ is also an $\mathrm{O}(n+m)$-fixed vector and thus has a non-zero image under the morphism ``$\mathrm{pr}$". Therefore, it remains to show that $u(v) \neq 0$.

By the proof of Helgason's theorem in \cite[Theorem 8.49]{Kn96}, the $\mathrm{O}(n+m)$-fixed vector can be obtained by
\[
\int_{\mathrm{O}(n+m)} k \cdot v' \, dk \neq 0,
\]
where $v'$ is the highest weight vector of $\gamma$. We may assume that $v'$ lies in $\gamma_1 \boxtimes \gamma_2$ by choosing a positive order properly. Notice that the $\mathrm{O}(n) \times \mathrm{O}(m)$-fixed vector
\[
z := \int_{\mathrm{O}(n) \times \mathrm{O}(m)} k \cdot v' \, dk ,
\]
is non-zero since
\[
\int_{\mathrm{O}(n+m)} k \cdot v' \, dk = \int_{\mathrm{O}(n+m)/\mathrm{O}(n) \times \mathrm{O}(m)} k_1 \cdot \left(\int_{\mathrm{O}(n) \times \mathrm{O}(m)} k \cdot v' \, dk\right) \, dk_1.
\]
Let $\langle \cdot, \cdot \rangle$ denote the invariant inner product on $\gamma$. Notice that $z \in \gamma_1 \boxtimes \gamma_2$ since $v' \in \gamma_1 \boxtimes \gamma_2$. Moreover,
\[
\langle v, z \rangle = c \left\langle  v, \int_{\mathrm{O}(n+m)/\mathrm{O}(n) \times \mathrm{O}(m)} k \cdot z\ dk\right\rangle = c \cdot \langle v, v \rangle \neq 0,
\]
where $c$ is a positive constant. Hence, $u(v) \neq 0$ since $\langle v, z \rangle \neq 0$. This completes the proof.
\end{proof}

\subsection{Building blocks of distinguished representations}
This section proves that a particular induced representation of $G_{2n}$ is $G'_{2n}$-distinguished. Let 
\[P=\left\{\begin{pmatrix}
A & B \\ 0 & C
\end{pmatrix}\ |\ A,B,C\in M_{n\times n}(\mathbb{C}), ~~  \det(A)\det(C)\neq 0\right\}\] 
be the parabolic subgroup of $G=G_{2n}$, whose Levi subgroup is $L\simeq G_n\times G_n$. Let $\gamma=(\frac{\det}{|\det|}\otimes |\det|^s)_{G_n}\boxtimes (\frac{\det}{|\det|}\otimes |\det|^{-s})_{G_n}$ be a character of $L$ for some $s \in \mathbb{C}$. Let $\Pi=\mathrm{Ind}_P^G(\gamma)$ be a representation of $G$, and assume that $\Pi$ is irreducible, which is equivalent to $s\notin \mathbb{Z}\backslash\{0\}$. The following is the main result of this section.
\begin{lemma}\label{dist_comp}
Assume that $s\in  \mathbb{C}$ and $s\notin \mathbb{Z}\backslash\{0\}$. The representation $\Pi=\mathrm{Ind}_P^G(\gamma)$, as defined above, is $H$-distinguished for $H=G'_{2n}$. 
\end{lemma}

\begin{proof}
As $\mathrm{Ind}_P^G(\gamma)$ is irreducible, it is isomorphic to $(\frac{\det}{|\det|}\otimes |\det|^{-s})_{G_n}\times (\frac{\det}{|\det|}\otimes |\det|^{s})_{G_n}$. Without loss of generality, we may assume that the real part $\mathrm{Re}(s)\geq 0$ by replacing $s$ with $-s$ if necessary.

Following Section \ref{double coset}, the double cosets of $P\setminus G/H$ are represented by the following $n+1$ elements:
\[
g_j = \begin{pmatrix}
I_{n-j} & 0_{(n-j)\times j} & 0_{(n-j)\times j} & 0_{(n-j)\times (n-j)} \\
0_{j\times (n-j)} & I_{j} & \sqrt{-1}I_j & 0_{j\times (n-j)} \\
0_{j\times (n-j)} & \sqrt{-1}I_j & I_j & 0_{j\times (n-j)} \\
0_{(n-j)\times (n-j)} & 0_{(n-j)\times j} & 0_{(n-j)\times j} & I_{n-j}
\end{pmatrix}, \quad j=0,1,\dots,n.
\]
Let $U_j$ denote the open subset $\cup_{i\geq j}Pg_iH$, so that $G=U_0\supset U_1\supset \dots\supset U_n$. Let $\mathcal{V}$ be the tempered bundle $\gamma\cdot\delta_P^{1/2}\times_P G$, and let $\Pi$ be realized as the space of Schwartz sections of $\mathcal{V}$. Define $\Pi_j:=\mathcal{S}(U_j,\mathcal{V}_{U_j})$ as the subspace of Schwartz sections of the restriction of $\mathcal{V}$ to $U_j$.

Let $\rho_j:=\Pi_j/\Pi_{j+1}$, and let $Q_j:=H\cap g_j P g_j^{-1}$. Then $\rho_j$ admits a decreasing filtration $\rho_j\supset \rho_{j,1}\supset\rho_{j,2}\supset \dots$ such that
\[
\rho_{j,k}/\rho_{j,k+1} \cong \mathcal{S}\left(H/Q_j, H\times_{Q_j}{^{g_j^{-1}}}(\gamma\cdot\delta_P^{1/2}) \otimes \mathrm{Sym}^k(\mathfrak{g}/(\mathfrak{h}+Ad(g_j)\mathfrak{p}))_{\mathbb{C}}^{\vee}\right)
\]
as $H$-modules. Hence, by \eqref{shap_lem},
\begin{equation}\label{h_iH=h_iP}
\mathrm{H}^{\mathcal{S}}_i(H,\rho_{j,k}/\rho_{j,k+1}) \cong \mathrm{H}^{\mathcal{S}}_i\left(P^j; \gamma\otimes \mathrm{Sym}^k(\mathfrak{g}/(\mathfrak{h}^j+\mathfrak{p}))_{\mathbb{C}}^{\vee} \otimes \delta_{H^j} \cdot \delta_P^{1/2} \cdot \delta_{P^j}^{-1}\right),
\end{equation}
where $H^j = g_j^{-1}Hg_j$ and $P^j = P\cap H^j$.

By the following exact sequences
\begin{equation}\label{filtration}
    0 \to \Pi_{j+1} \to \Pi_j \to \rho_j \to 0, \quad 0 \leq j \leq n-1,
\end{equation}
to show that $\mathrm{Hom}_H(\Pi,\mathbb{C})\neq 0$, i.e., $\mathrm{H}^{\mathcal{S}}_0(H,\Pi) \simeq \mathrm{Hom}_H(\Pi,\mathbb{C})^{\vee} \neq 0$, it suffices to prove that 
$\mathrm{H}^{\mathcal{S}}_0(H,\rho_j)=0$, $\mathrm{H}^{\mathcal{S}}_1(H,\rho_j)=0$ for all $0 \leq j \leq n-1$, and
$\mathrm{H}^{\mathcal{S}}_0(H,\Pi_n) \neq 0$. Therefore, by \eqref{lim_hom}, it suffices to show that \begin{itemize}
\item $\mathrm{H}^{\mathcal{S}}_i(H,\rho_{j,k}/\rho_{j,k+1})=0$ for $i=0,1$, $\forall\ 0 \leq j \leq n-1$, and $\forall \ k \in \mathbb{Z}_{\geq 0}$,
\item $\mathrm{H}^{\mathcal{S}}_0(H,\Pi_n) \neq 0$.
\end{itemize}

By definition, one has
\begin{equation}\label{pj}
\begin{aligned}
P^j = P \cap H^j = \left\{\begin{pmatrix}
a & b & -\sqrt{-1}\cdot\overline{b} & c \\
0 & d & 0 & e \\
0 & 0 & \overline{d} & \sqrt{-1}\cdot\overline{e} \\
0 & 0 & 0 & f
\end{pmatrix} \;\middle|\; 
\begin{aligned}
& a,f \in G_{(n-j)}'; \; b \in M_{(n-j)\times j}(\mathbb{C}),\\
&d \in G_j,\; c \in M_{(n-j)\times j}(\mathbb{R}),\\  
& e \in M_{j\times (n-j)}(\mathbb{C})
\end{aligned}
\right\},
\end{aligned}
\end{equation}
and 
\begin{equation}\label{g/hj+p}
\mathfrak{g}/\mathfrak{h}^j+\mathfrak{p} \simeq \left\{\begin{pmatrix}
0 & 0 & 0 & 0 \\
0 & 0 & 0 & 0 \\
0 & 0 & 0 & 0 \\
x & 0 & 0 & 0 
\end{pmatrix} \;\middle|\; x \in \sqrt{-1}\cdot M_{(n-j)\times (n-j)}(\mathbb{R})\right\},
\end{equation}
where the isomorphism is as $K\cap L^j$-modules, and $L^j$ denotes the standard Levi subgroup of $P^j$.

Notice that $\delta_{H^j}=1$ since $H^j$ is reductive, and one has $\delta_{P}^{1/2}\cdot\delta_{P^j}^{-1}=1$ by direct computation using \eqref{pj}.

For $\mathrm{H}^{\mathcal{S}}_i(H,\rho_{j,k}/\rho_{j,k+1})=0$ ($i\in\mathbb{Z},0 \leq j \leq n-1$): 

Let $L^j$ denote the standard Levi subgroup of $P^j$, and let $\mathfrak{n}^j$ denote the nilradical of $\mathfrak{p}^j$. Let $T$ be the diagonal subgroup $\left\{\mathrm{diag}(\underbrace{t,\dots,t}_{n-j},1,\dots,1,\underbrace{t^{-1},\dots,t^{-1}}_{n-j})\ \middle|\ t>0\right\}$. Observe that $T$ acts on both $(\mathfrak{g}/(\mathfrak{h}^j+\mathfrak{p}))^{\vee}$ and $\mathfrak{n}^j$  with positive weight $t^{\mathbb{Z}_{\geq 0}}$. Moreover, the actions of $\mathfrak{n}^j$ on $\gamma$ and $\mathfrak{g}/(\mathfrak{h}^j+\mathfrak{p})$ are both trivial. 
Hence, when $k>0$, the action of $T$ on $\mathfrak{n}^j$-homology of $\gamma\otimes \mathrm{Sym}^k(\mathfrak{g}/(\mathfrak{h}^j+\mathfrak{p}))^{\vee}$ has no trivial composition factor; When $k=0$, similarly, the $\mathfrak{n}^j$-homology $\gamma$ has no trivial composition factor under the action of $T$. By the Lemma \ref{spectral sequence},  the Schwartz homology $\mathrm{H}^{\mathcal{S}}_{i'}(L^j,\mathrm{H}^{\mathcal{S}}_{i''}(\mathfrak{n}^j,\bullet)\Rightarrow \mathrm{H}^{\mathcal{S}}_{i'+i''}(P^j,\bullet)$. Combining \eqref{h_iH=h_iP}, Lemma \ref{van_hom} implies that $\mathrm{H}^{\mathcal{S}}_i(H,\rho_{j,k}/\rho_{j,k+1})$ vanishes for all $i$.

For $\mathrm{H}^{\mathcal{S}}_0(H,\Pi_n) \neq 0$: We have
\[
\mathrm{H}^{\mathcal{S}}_0(H,\Pi_n) \simeq \mathrm{H}^{\mathcal{S}}_0(P^n, \gamma\otimes \delta_{H^n} \cdot \delta_P^{1/2} \cdot \delta_{P^n}^{-1}) =\mathbb{C},
\]
since $\gamma$ is a trivial representation of $P^n$. This finishes the proof. 
\end{proof}

\begin{lemma}\label{dist_comp_ev}
The representations $|\det|_{G_m}^{s} \times |\det|_{G_m}^{-s}$ are $G_{2m}'$-distinguished when $s\in \mathbb{C}\backslash \mathbb{Z}$.
\end{lemma}
\begin{proof}
The result follows from a similar argument as the above lemma. We omit the proof.
\end{proof}

\subsection{Proof of the sufficient part of Theorem \ref{main result_generic}}\label{sec:proof_gen}
We first remark that an irreducible representation is generic if and only if it is a standard representation (see \cite[Appendix A]{Kem15a} for details). Thus, for an irreducible generic representation $\pi$ satisfying the condition in Theorem \ref{main result_generic}, it can be written as
\[
\pi\simeq \chi_1\times \dots\times \chi_r\times \sigma_1\times \dots\times \sigma_s\quad (r+2s=n)
\]
where $\chi_i=\left(\frac{\det}{|\det|}\right)^{k_i}$ for some even integer $k_i$, and 
\[
\sigma_j\simeq \left(\frac{\det}{|\det|}\right)^{k_j'}\otimes \left(|\det|^{t_j}\times |\det|^{-t_j}\right)
\]
for some integer $k_j'$ and complex numbers $t_j\notin \mathbb{Z}\setminus\{0\}$. Then the sufficient part of Theorem \ref{main result_generic} follows from Lemma~\ref{dist_comp} and Lemma~\ref{dist_comp_ev}.

\subsection{Proof of the sufficient part of Theorem \ref{equi_thm}}\label{sec:proof_uni}
It remains to show that any irreducible unitary representation of $G_n$ satisfying conditions (i) and (ii) in Theorem~\ref{equi_thm} is $G_n'$-distinguished. By Section~\ref{parabolic induced}, it suffices to show that the following representations of $G_n$ are $G_n'$-distinguished:
\begin{itemize}
\item[(i)] $(\frac{\det}{|\det|})^{k},$ for $ k\in 2\mathbb{Z}$,  
\item[(ii)] $(|\det|_{G_{m}}^{u}\otimes (\frac{\det}{|\det|})_{G_{m}}^{k})\times (|\det|_{G_{m}}^{-u}\otimes (\frac{\det}{|\det|})_{G_{m}}^{k})$, for $u\in \sqrt{-1}\mathbb{R}\backslash\{0\}$, $k\in \mathbb{Z}$, 
\item[(iii)] $\tau_{k,t,m}:=(\frac{\det}{|\det|})_{G_{2m}}^{k}\otimes (|\det|_{G_{m}}^{t}\times |\det|_{G_{m}}^{-t})$, for $k\in\mathbb{Z}$, $0<t<1$,
\item[(iv)]  $\left(|\det|_{G_{2m}}^{u}\otimes \tau_{k,t,m} \right)  \times  \left(|\det|_{G_{2m}}^{-u}\otimes \tau_{k,t,m} \right)$, for $u\in \sqrt{-1}\mathbb{R}\backslash\{0\}$, $k\in\mathbb{Z}$, $0<t<1$,
\item[(v)] $(\frac{\det}{|\det|})_{G_m}^{k}\times (\frac{\det}{|\det|})_{G_m}^{k},$ for $ k\in 2\mathbb{Z}+1$. 
\end{itemize}

For (i), the representation is trivial on $G_n'$. 

For (ii), when $k\in 2\mathbb{Z}+1$, notice that $(\frac{\det}{|\det|})^{k-1}$ is trivial on $G_{2m}'$, by Lemma \ref{dist_comp}, we conclude that (ii) is also $G_{2m}'$-distinguished. When $k\in 2\mathbb{Z}$, it suffices to show $|\det|_{G_{m}}^{u} \times |\det|_{G_{m}}^{-u}$  is $G_{2m}'$-distinguished, which follows from Lemma \ref{dist_comp_ev}.

For (iii), the argument is similar to that of (ii).

For (iv), one can write it as the product of 
\[((\frac{\det}{|\det|})_{G_{m}}^{k}\otimes |\det|_{G_{m}}^{u}\otimes |\det|_{G_{m}}^{t})\times((\frac{\det}{|\det|})_{G_{m}}^{k}\otimes |\det|_{G_{m}}^{-u}\otimes |\det|_{G_{m}}^{-t})\] and \[((\frac{\det}{|\det|})_{G_{m}}^{k}\otimes |\det|_{G_{m}}^{u}\otimes |\det|_{G_{m}}^{-t})\times((\frac{\det}{|\det|})_{G_{m}}^{k}\otimes |\det|_{G_{m}}^{-u}\otimes |\det|_{G_{m}}^{t}).\] These two representations are $G_{2m}'$-distinguished by the same argument as (ii), so (iv) is $G_{4m}'$-distinguished by Theorem \ref{induction}. 

For (v), the result follows from Lemma~\ref{dist_comp}. \qed

\section{Non-vanishing on the distinguished minimal $K$-type}\label{non-vanishing}
By Theorem~\ref{induction}, Theorem~\ref{main result_generic} and Theorem~\ref{main result_unitary}, to prove Theorem~\ref{Non}, we first deal with $G'_2$-distinguished representations of $G_2$, showing that the period of $\Hom_{G_2'}(\Pi,\mathbb{C})$ does not vanish on the distinguished minimal $K$-type.

\subsection{Distinguished minimal $K$-type non-vanishing: generic case}\label{generic case non-vanishing}
 Following Section~\ref{Non-van}, let $V^*=\left( (\frac{\det}{|\det|})^k\otimes |\det|^s \right)_{G_n}\times\left( (\frac{\det}{|\det|})^k\otimes |\det|^{-s}\right)_{G_n}$ and $W$ be the trivial representation of $P'$. Let $\mathcal{P}$ be a period of $\Hom_{G_n'}(\Pi,\mathbb{C})$. We do the generic calculation in this case. Firstly, note that by the calculations in the above sections, the $\mathcal{P}$ has nonzero image in $\mathcal{S}(U,\mathcal{V}|_U)'^{G_n'}$, thus by Lemma~\ref{regular}, $\mathcal{P}$ is regular. 
 
 Let $m=\begin{bmatrix}
     m_1 & \\
         & m_2\\
 \end{bmatrix}\in M'A'$. The $M'A'$-invariance of the distribution kernel amounts to
 \begin{equation}\label{inv}
     f(u)=\chi(m)f(m_2\cdot u\cdot m_1^{-1} ),
 \end{equation} where $\chi(m)=(\frac{\det m_1}{|\det m_1|})^k|\det m_1|^{-n}\cdot (\frac{\det m_2}{|\det m_2|})^k|\det m_2|^{n}$, $u\in M_{n\times n}(\mathbb{C})$, and $f\in\mathcal{G}(N_{-})\simeq \mathcal{G}(M_{n\times n}(\mathbb{C}))$. The critical point that this method is hard to use in $n\geq 2$ is that killing by $\mathfrak{n}'$ is hard to calculate. Hence, from now on, we restrict to $n=1$. Then $\mathfrak{n}'$ is one dimensional, and we need only to look at $X=\begin{bmatrix}
     0 & 1\\
     0 & 0\\
 \end{bmatrix}\in\mathfrak{n}'$. Note that:
 $$
     \begin{bmatrix}
         1 & \\
         u & 1\\
     \end{bmatrix}\begin{bmatrix}
         1 & t\\
           & 1
     \end{bmatrix}=\begin{bmatrix}
         1 & t\\
         u & tu+1
     \end{bmatrix}=\begin{bmatrix}
         (tu+1)^{-1} & t\\
               & tu+1\\
     \end{bmatrix}\begin{bmatrix}
         1 & \\
         \frac{u}{tu+1}  & 1  \\
     \end{bmatrix}.
 $$
 Consequently, we find the action of $\exp(tX)$ on $f$ is given by:
 \begin{equation*}
     (\exp(tX)\cdot f)(u)=f\left(\frac{u}{1+tu}\right)\left|\frac{1}{1+tu}\right|^{2+2s}.
 \end{equation*}
We first use the Euclidean coordinate that $u=x+iy$. Then we get
\begin{equation}\label{kill}
   \frac{d(\exp(tX)\cdot f)}{dt}\bigg|_{t=0}=\frac{\partial f}{\partial x}(x^2-y^2)+\frac{\partial f}{\partial y}(2xy)=0.
\end{equation}
 
Now we turn to polar coordinates, and Equation~\eqref{inv}, $f(tu)=\sign(t)^{k}(\frac{1}{|t|})^{1+s}f(u)$, implies that $f(r,\theta)=(\frac{1}{r})^{1+s}g(\theta)$ if $f$ is supported on $N_{-}$. On the other hand, Equation~\eqref{kill} simplifies to
 \begin{equation*}
     r\cos\theta\frac{\partial f}{\partial r}+\sin\theta \frac{\partial f}{\partial \theta}=0,
 \end{equation*} which shows $g(\theta)=C\sign(\sin\theta)^k|\sin\theta|^{1+s}$ for some nonzero $C$. $C$ will not influence the conclusion and we assume to be 1. Hence, generically, the distribution kernel is given by $f=\sign(\sin\theta)^k(\frac{1}{r})^{1+s}|\sin\theta|^{1+s}$. For the use of meromorphic extension, we observe that
 \begin{equation}\label{ext}
     P(a,b)=\frac{1}{\Gamma(a)\Gamma(b)}\int\limits_{-1}^{1} (1-s)^{a-1}(1+s)^{b-1}ds
 \end{equation}
 can be extended holomorphically to $(a,b)\in\mathbb{C}^2$. More explicitly, after extension, $P(a,b)=\frac{2^{a+b-1}}{\Gamma(a+b)}$. 
 
 Calculation of $i^*$: note the coordinate change from $N_{-}$ to $\mathrm{SU}(2)$
 \begin{equation}
     \begin{bmatrix}
         1 & \\
     u & 1\\
     \end{bmatrix}=\begin{bmatrix}
         \frac{1}{(1+|u|^2)^{1/2}} & *\\
              & (1+|u|^2)^{1/2}
     \end{bmatrix}\begin{bmatrix}
     \frac{1}{(1+|u|^2)^{1/2} } &  \frac{-\overline{u}}{(1+|u|^2)^{1/2}}\\
      \frac{u}{(1+|u|^2)^{1/2}} &  \frac{1}{(1+|u|^2)^{1/2} }\\
     \end{bmatrix}
 \end{equation}
 \begin{claim} Let ${K_T}_{s}=\sign(\sin\theta)^k(\frac{1}{r})^{1+s} |\sin\theta|^{1+s}$ be the holomorphic family of distribution kernels.
    Then ${K_T}_{s}(i^*h)$ equals $\mathcal{P}(h)$ up to some non-zero scalar for $\mathrm{Re}(s)\geq 0$ and $s\notin \mathbb{Z}\backslash\{0\}$, where $h$ is a non-zero $\mathrm{O}(2)$-invariant function in the distinguished minimal $K$-type.
 \end{claim}
 First of all, we note that $\Gamma(-\frac{s-1}{2})\Gamma(\frac{3s+1}{2})$ is holomorphic and non-zero at $\mathrm{Re}(s)\geq 0$ and $s\notin \mathbb{Z}\backslash\{0\}$.
 \begin{Case}
     $k$ is even. Then the distinguished minimal $K$-type is $(k,k)$, and is realized as constant functions on $\mathrm{SU}(2)$. Under this circumstance, suppose $h=1$,
     \begin{align*}
           \left(  \sign(\sin\theta)^k(\frac{1}{r})^{1+s} |\sin\theta|^{1+s}\right)(i^*h)  &\  = \int\limits_0^{2\pi}\int\limits_0^{+\infty} \left( \frac{1}{1+r^2}\right)^{1+s} (\frac{1}{r})^{1+s} |\sin\theta|^{1+s} rdr d\theta\\
        = \frac{1}{2^{s+1}}\int\limits_{-1}^1 (1-t)^{-\frac{s-1}{2}-1}(1+t)^{\frac{3s+1}{2}-1}dt \int\limits_0^{2\pi}&|\sin\theta|^{1+s}d\theta =\frac{\Gamma(-\frac{s-1}{2})\Gamma(\frac{3s+1}{2})}{2\Gamma(s+1)}\int\limits_0^{2\pi}|\sin\theta|^{1+s}d\theta .
     \end{align*}
    Since the result has no poles when $\mathrm{Re}(s)\geq 0$ and $s\notin \mathbb{Z}\backslash\{0\}$, we observe ${K_T}_s(i^*h)$ is holomorphic and non-zero on $\mathrm{Re}(s)\geq 0$ and $s\notin \mathbb{Z}\backslash\{0\}$. Thus, Claim follows from Proposition~\ref{mero}.
 \end{Case}
 
 \begin{Case}
      $k$ is odd. Then the distinguished minimal $K$-type is $(k+1,k-1)$, and on $\mathrm{SU}(2)=\left\{\begin{bmatrix}
          z_1 & z_2\\
          -\overline{z_2} & \overline{z_1}\\
      \end{bmatrix}\Bigg||z_1|^2+|z_2|^2=1\right\}$, this $K$-type is spanned by $\overline{z_1}z_2,|z_1|^2-|z_2|^2,z_1\overline{z_2}$. The unique $\mathrm{O}(2)$-invariant subspace is generated by $h=\overline{z_1}z_2-z_1\overline{z_2}$. Calculation shows that
      $$ i^*(h)=2\sqrt{-1}\left(\frac{1}{1+r^2}\right)^{2+s}r\sin\theta,
      $$ hence,
      \begin{align*}
          & \left( \sign(\sin\theta)^k(\frac{1}{r})^{1+s} |\sin\theta|^{1+s}\right)(i^*h)  = 2\sqrt{-1}\int\limits_0^{2\pi}\int\limits_0^{+\infty} \left( \frac{1}{1+r^2}\right)^{2+s} (\frac{1}{r})^{1+s} |\sin\theta|^{2+s} r^2dr d\theta\\
           &= \frac{2\sqrt{-1}}{2^{s+1}}\int\limits_{-1}^1 (1-t)^{-\frac{s-2}{2}-1}(1+t)^{\frac{3s+2}{2}-1}dt \int\limits_0^{2\pi}|\sin\theta|^{2+s}d\theta\\
           &=\frac{\sqrt{-1}\cdot \Gamma(-\frac{s-2}{2})\Gamma(\frac{3s+2}{2})}{\Gamma(s+2)}\int\limits_0^{2\pi}|\sin\theta|^{2+s}d\theta.
      \end{align*}
      Similarly, we observe ${K_T}_s(i^*h)$ is holomorphic and non-zero on $\mathrm{Re}(s)\geq 0$ and $s\notin \mathbb{Z}\backslash\{0\}$ since $\frac{\Gamma(-\frac{s-2}{2})\Gamma(\frac{3s+2}{2})}{\Gamma(s+2)}$ is holomorphic and non-zero on this domain. This completes the proof of Claim.
 \end{Case}

For other $K$-type, similar computation shows after meromorphic continuation, they are holomorphic and non-zero on $\mathrm{Re}(s)\geq 0$ and $s\notin \mathbb{Z}\backslash\{0\}$.

\begin{proof}[Proof of the generic case of Theorem~\ref{Non}]
Assume that $\pi$ is an irreducible generic representation of $G_n$. By Theorem~\ref{Lang_cls}, $\pi$ can be written as $\mathrm{Ind}^{G_n}_{B_n} (\lambda_1 \boxtimes \cdots \boxtimes \lambda_n)$, where $\lambda_i(z) = \left(\frac{z}{|z|}\right)^{m_i}|z|^{s_i}$. Moreover, we may rearrange the order of the $\lambda_i$'s such that $m_1 \geq \dots \geq m_n$. In this case, the lowest $K$-type of $\pi$ is $(m_1, \dots, m_n)$. 

Assume that $\pi$ satisfies condition (i) of Theorem~\ref{main result_generic}. Then, for each odd $m_i$, it appears an even number of times in $(m_1, \dots, m_n)$. Write $(m_1, \dots, m_n)$ as $(k_1, \dots, k_1; \dots; k_r, \dots, k_r)$, where $k_1 > k_2 > \dots > k_r$, and each $k_i$ has multiplicity $n_i$. For each odd $k_i$ (where $n_i$ is even, since $\pi$ is $G_n'$-distinguished), by replacing $(\underbrace{k_i, \dots, k_i}_{n_i})$ with 
\[
(\underbrace{k_i+1, \dots, k_i+1}_{n_i/2}, \underbrace{k_i-1, \dots, k_i-1}_{n_i/2}),
\]
one gets $(\widetilde{m_1}, \dots, \widetilde{m_n})$. By Lemma~\ref{even}, $(\widetilde{m_1}, \dots, \widetilde{m_n})$ is the unique distinguished minimal $K$-type of $\pi$.

From the result for $G_2$, the $G_2'$-period of any $G_2'$-distinguished representation does not vanish on the unique distinguished minimal $K$-type. Write $\pi = \pi_1 \times \dots \times \pi_l$, where each $\pi_i$ is either a $G_2'$-distinguished representation of $G_2$ or a $G_1'$-distinguished representation of $G_1$. By Theorem~\ref{induction}, the $G_n'$-morphism from $\pi$ to $\mathbb{C}$ is non-vanishing at the $K$-type $(\gamma_1, \dots, \gamma_l)$, where $\gamma_i$ is the distinguished minimal $K$-type of $\pi_i$ for $1 \leq i \leq l$. By the result for $G_2$, $(\gamma_1, \dots, \gamma_l)=(\widetilde{m_1}, \dots, \widetilde{m_n})$. This completes the generic part of the theorem.
\end{proof} 

\begin{example}
Consider the $G_4'$-distinguished  generic unitary representation of $G_4$:
\[\begin{aligned}\pi &=(\frac{\det}{|\det|})\times (\frac{\det}{|\det|})\times (\frac{\det}{|\det|}\otimes |\det|^{1/2})\times (\frac{\det}{|\det|}\otimes |\det|^{-1/2})\\ &\simeq  \left((\frac{\det}{|\det|})\times (\frac{\det}{|\det|})\right)\times \left((\frac{\det}{|\det|}\otimes |\det|^{1/2})\times (\frac{\det}{|\det|}\otimes |\det|^{-1/2})\right).\end{aligned}\]
Let $\tau_{s}:=\left((\frac{\det}{|\det|}\otimes |\det|^{s})\times (\frac{\det}{|\det|}\otimes |\det|^{-s})\right)$. Then $\pi$ has lowest $K$-type $(1,1,1,1)$. By changing from $(1,1,1,1)$ to $(2,2,0,0)$ as the proof above, one can get its distinguished minimal $K$-type  $\sigma=(2,2,0,0)$. Notice that this $K$-type has multiplicity two in $\pi$. One can decompose the $\sigma$-isotypic space of $\pi$ as $\sigma_1\oplus \sigma_2$, where $\sigma_1$ (resp. $\sigma_2$) is the irreducible $\sigma$-subspace in $\pi$ which comes from induction from the $K$-type $(1,1)\boxtimes (1,1)$ (resp. $(2,0)\boxtimes (2,0)$) of $\tau_0\times \tau_{1/2}$ of $G_2\times G_2$.  The above proof shows that the $G_4'$-period of $\pi$ does not vanish on $\sigma_2$. Actually, by the proof and the result of $G_2$, the $G_4'$-period of $\pi$ vanishes on $\sigma_1$.
\end{example}

\subsection{Distinguished minimal $K$-type non-vanishing: unitary case}
As remarked in Section~\ref{generic case non-vanishing}, the method becomes computationally challenging for $n \geq 2$. In this section, we will construct the period for distinguished unitary representations  and demonstrate that the periods do not vanish on the distinguished minimal $K$-types of these representations, except when the Langlands parameters contain certain exceptional cases:
\begin{equation}\label{exception}
    \lambda_i(z) = \left(\frac{z}{|z|}\right)^{m_i}|z|^{s_i} \quad \text{with} \quad m_i \in 2\mathbb{Z}+1 \quad \text{and} \quad s_i \in \mathbb{Z}\setminus\{0\}.
\end{equation}
If we write a distinguished unitary representation $\pi$ as the form in the Theorem \ref{equi_thm}, then ``the Langlands parameters of $\pi$ contain no parameters of the form \eqref{exception}'' is equivalent to $(\frac{\det}{|\det|})^{2k+1}_{G_m}$ for $k\in \mathbb{Z},m\geq 2$ does not show up in the product expression of $\pi$. 

The building block for this construction is the concrete realization of the period when $n=1$ and its meromorphic continuation. This follows directly from the discussion in Section~\ref{generic case non-vanishing}. More directly, the period also arises from integrating functions over the open double coset $BwH$:
\begin{equation*}
    \mathcal{P}:f \mapsto \int_{H^w\backslash H} f(wh) \, dh, \quad H^w = w^{-1}Bw \cap H.
\end{equation*}
It is straightforward to show that for $\text{Re}(s) > 0$, the integral is absolutely convergent and admits a meromorphic continuation to $s \in \mathbb{C}$ by \eqref{ext}, which is holomorphic on $s\in\mathbb{C}\backslash \mathbb{Z}$. The details are left to interested readers.

To construct the period for an irreducible unitary representation $\pi$, we first construct the period for a principal series representation that realizes $\pi$ as a quotient. We then show that the period of such a principal series is multiplicity-free, which implies that the period must factor through $\pi$.

\subsubsection{Construction of the periods}\label{construction of period}
To construct the $G_n'$-period map for all irreducible unitary distinguished representations of $G_n$, by parabolic induction \ref{induction}, and the classification of distinguished representations, one only needs to construct the period for the representations 
\[(|\det|^s)_{G_n}\times (|\det|^{-s})_{G_n} \]
for $s\in (\mathbb{C}\setminus \mathbb{Z})$, and the representations 
\[(\frac{\det}{|\det|}\otimes |\det|^s)_{G_n}\times (\frac{\det}{|\det|}\otimes |\det|^{-s})_{G_n} \]
for $s\in (\mathbb{C}\setminus \mathbb{Z})$.

For $s\in \mathbb{C}$ and $\epsilon\in \{0,1\}$, let 
\begin{align*}
    \pi_{\epsilon,s,n}&:=\left((\frac{\det}{|\det|})^{\epsilon}\otimes |\det|^s\right)_{G_n}\times \left((\frac{\det}{|\det|})^{\epsilon}\otimes |\det|^{-s}\right)_{G_n}, \\
    \Pi_{\epsilon,s,n}&:=\mathrm{Ind}_{B_{2n}}^{G_{2n}}(\mathop{\boxtimes}\limits_{i=1}^{2n}\lambda_i),
\end{align*}
where $\lambda_i$ is the character of $\mathbb{C}^{\times}$ defined by $\lambda_i(z):=(\frac{z}{|z|})^{\epsilon}|z|^{n+1+s-2i}$ for $1\leq i\leq n$, and $\lambda_{n+i}(z):=(\frac{z}{|z|})^{\epsilon}|z|^{n+1-s-2i}$ for $1\leq i\leq n$. Notice that $\pi_{\epsilon,s,n}$ is a quotient of  $\Pi_{\epsilon,s,n}$ when $s\in \mathbb{C}\setminus \mathbb{Z}_{<0}$.

We will construct the $G_{2n}'$-period of $\Pi_{\epsilon,s,n}$. To show that the period factors through $\pi_{\epsilon,s,n}$, we need the following multiplicity-one lemma on $\mathrm{Hom}_{G_{2n}'}(\Pi_{\epsilon,s,n},\mathbb{C})$. 

\begin{lemma}\label{mult<=1}
Let $\Pi=\mathrm{Ind}_{B_{2n}}^{G_{2n}}(\mathop{\boxtimes}\limits_{i=1}^{2n}\lambda_i)$, where $\lambda_i(z)=(\frac{z}{|z|})^{m_i}|z|^{s_i}$, $s_i=n+1+s-2i$ with $s\in \mathbb{C}\setminus \frac{1}{2}\mathbb{Z}$,  and $s_{n+i}=-s_{n+1-i}$,  $m_{n+i}=-m_{n+1-i}$ for $1\leq i\leq n$. 
Then 
$$\dim_{\mathbb{C}}\mathrm{Hom}_{G_{2n}'}(\Pi,\mathbb{C})= 1 .$$ 
\end{lemma}

\begin{proof}
Similar to the proof of Lemma~\ref{dist_comp}, we compute the Schwartz homology of $\Pi|_{G_{2n}'}$ by the double coset decomposition $B_{2n}\backslash G_{2n}/G_{2n}'$. We give a sketch of the proof. Following the notation in \ref{double coset}, we set $G=G_{2n}$, $H=G_{2n}'$, and $P=B_{2n}$. By Section~\ref{double coset}, the double cosets are parameterized by $W_{2n,2}$. We partition this set into two cases.

\textbf{Case 1:} $w\in W_{2n,2}, w(l)\neq 2n+1-l$ for some $l$. \\
Let $l'$ be the smallest index satisfying this property, so $1\leq l'\leq n$. Let $l''=w(l')$. Assume that such $w$ corresponds to $g_j$ as in \ref{double coset}. Then 
\[P^j=P\cap g_j^{-1}\cdot H\cdot g_j=\left\{\begin{pmatrix}
X & 0 & 0 & 0 & 0 \\ 0 & z & * &  0 & * \\ 0 & 0 & Y & 0 & *\\ 0 &  0 & 0 & \overline{z} & *\\ 0 & 0 & 0 & 0 & *
\end{pmatrix}\ \middle| \begin{array}{l} X\ \text{is an}\ (l'-1)\times (l'-1) \ \text{diagonal matrix},\\ z \in \mathbb{C}^{\times}; Y\in \mathrm{GL}_{l''-l'-1}\ \text{and}\ ``*" \ \text{are the sub-matrix}\\ \text{satisfying certain conditions} \end{array}
\right\},\]
where $\begin{pmatrix}
z & * & 0 \\ 0 & Y & 0 \\ 0 & 0 & \overline{z}\end{pmatrix}$ means a non-zero real number if $l'=l''$.

Let
\[T=\left\{\begin{pmatrix}
I_{l'-1} &  &  & & \\ & a &  & &\\ & & I_{l''-l'-1} & &\\ & & & a & \\ & & & & I_{2n-l''}
\end{pmatrix}\ \middle| \ a\in \mathbb{R}_{>0} \right\}\]
be the subgroup of $P^j$. Then the weights of $T$ on the nilpotent radical of $\mathfrak{p}^j$ are $|a|^{\mathbb{Z}_{>0}}$. Moreover, the action of $T$ on `` $\sigma\otimes \mathrm{Sym}^k(\dots)\dots$" is as follows: for $1\leq l'' \leq n$(resp. $n+1\leq l''\leq 2n$), the action of $T$ has weights $|a|^{2s+\mathbb{Z}}$(resp. $|a|^{\mathbb{Z}_{>0}}$). Thus, by the discussion as Lemma \ref{dist_comp} using spectral sequence, one concludes  $\mathrm{H}_i^{\mathcal{S}}(H,\rho_{j,k}/\rho_{j,k+1})=0$ for all $i$ in \eqref{shap_lem}.

\textbf{Case 2:} $w\in W_{2n,2}$ such that $w(l)= 2n+1-l$ for all $l$.  \\
Let $g_{j'}$ be the corresponding element as in \ref{double coset}. Note that this corresponds to the unique open set $Pg_{j'}H$. In this case, the term $\mathrm{Sym}^k(\dots)$ in $\sigma\otimes \mathrm{Sym}^k(\dots)\dots$ in \eqref{shap_lem} does not appear, and $\sigma\otimes\dots$ is the one-dimensional trivial representation of $P^{j'}$.  

Combining the results from the two cases, we obtain the lemma by applying the exact sequence analogous to \eqref{filtration}.
\end{proof}

Recall that for any $w \in S_{2n}$, one can formally define the intertwining operator between the induced modules $\mathrm{Ind}_B^G(\lambda)$ and $\mathrm{Ind}_B^G(w\lambda)$ for any character $\lambda$, as described in \cite[III. 3]{D75}:
\begin{equation}\label{intertw}
    \iota_{w} : \mathrm{Ind}_B^G(\lambda) \longrightarrow  \mathrm{Ind}_B^G(w\lambda), \qquad  f  \mapsto  \int_{w^{-1}Nw \cap \overline{N}} f(\overline{n}w^{-1}g) \, d\overline{n},
\end{equation}
where $\overline{N}$ is the unipotent subgroup corresponding to the negative roots. 

For $s \in \mathbb{C}$ and $\epsilon \in \{0,1\}$, define $\lambda_{\epsilon,s} := \lambda_1 \boxtimes \dots \boxtimes \lambda_{2n}$ with
\[
\lambda_i(z) := \left(\frac{z}{|z|}\right)^{\epsilon} |z|^{n+1+s-2i} \quad \text{for } 1 \leq i \leq n,
\]
and
\[
\lambda_{n+i}(z) := \left(\frac{z}{|z|}\right)^{\epsilon} |z|^{n+1-s-2i} \quad \text{for } 1 \leq i \leq n.
\]

Let $\widetilde{w}_1$ be the longest element in $S_n$, the group of permutations of the set $\{1,\dots,n\}$, and let $\widetilde{w}_2 = (n+1, n, \dots, 1)$ under the cycle notation. Consider the intertwining operator
\begin{equation}\label{int_w}
\begin{aligned}
\mathrm{Ind}_{B}^G(\lambda_{\epsilon,s}) & \xrightarrow{\iota_{\widetilde{w}_1,\epsilon,s}}  \mathrm{Ind}_{B}^G(\lambda_n \boxtimes \lambda_{n-1} \boxtimes \dots \boxtimes \lambda_1 \boxtimes \lambda_{n+1} \boxtimes \dots \boxtimes \lambda_{2n}) \\
& \xrightarrow{\iota_{\widetilde{w}_2,\epsilon,s}}  \mathrm{Ind}_{B}^G(\lambda_{n+1} \boxtimes \lambda_n \boxtimes \lambda_{n-1} \boxtimes \dots \boxtimes \lambda_1 \boxtimes \lambda_{n+2} \boxtimes \dots \boxtimes \lambda_{2n}).
\end{aligned}
\end{equation}
Notice that $\iota_{\widetilde{w}_j,\epsilon,s}$ can be defined over the space $\mathrm{Ind}_{K \cap B}^K \left(\prod_{i=1}^{2n} \left(\frac{z}{|z|}\right)^{\epsilon}\right)$. Although the integral defining $\iota_{\widetilde{w}_2,\epsilon,s}$ is not convergent, we can normalize it using meromorphic continuation and requiring it to be the identity on the lowest $K$-type (see \cite[III. 4]{D75}). After normalization, the family of intertwining operators, still denoted by $\iota_{\widetilde{w}_2,\epsilon,s}$, depends meromorphically on $s \in \mathbb{C}$ and is holomorphic and invertible on $s \in \mathbb{C} \setminus \mathbb{Z}$.

Applying similar intertwining operators as $\iota_{\widetilde{w}_2,\epsilon,s}$, we obtain
\[
\begin{aligned}
& \mathrm{Ind}_{B}^G(\lambda_{n+1} \boxtimes \lambda_n \boxtimes \lambda_{n-1} \boxtimes \dots \boxtimes \lambda_1 \boxtimes \lambda_{n+2} \boxtimes \dots \boxtimes \lambda_{2n}) \\
\longrightarrow \  
 &\mathrm{Ind}_{B}^G(\lambda_{n+1} \boxtimes \lambda_n \boxtimes \lambda_{n+2} \boxtimes \lambda_{n-1} \boxtimes \dots \boxtimes \lambda_{2n} \boxtimes \lambda_1).
\end{aligned}
\]
Combining these intertwining operators, we obtain the \textbf{period}
\begin{equation}\label{period def}
    \begin{aligned}
    \psi_{\epsilon,s}: &\mathrm{Ind}_{B}^G(\lambda_{\epsilon,s}) \stackrel{\iota_{\widetilde{w}_3,\epsilon,s}}{\longrightarrow} \mathrm{Ind}_{B}^G\left((\lambda_{n+1} \boxtimes \lambda_{n}) \boxtimes (\lambda_{n+2} \boxtimes \lambda_{n-1}) \boxtimes \dots \boxtimes (\lambda_{2n} \boxtimes \lambda_1)\right) \\
    & \stackrel{\simeq}{\longrightarrow} \pi_{\epsilon,n-s-1,1} \times \pi_{\epsilon,n-s-3,1} \times \dots \times \pi_{\epsilon,-n-s+1,1} \stackrel{\mathcal{P} \times \dots \times \mathcal{P}}{\longrightarrow} \mathbb{C},
\end{aligned}
\end{equation}
where $\mathcal{P} \times \dots \times \mathcal{P}$ is the parabolic induced period from $G_2$, as described in Section~\ref{parabolic induced}. The map $\psi_{\epsilon,s}$ is holomorphic for $s \in \mathbb{C} \setminus \mathbb{Z}$.

By Lemma~\ref{mult<=1}, the period of $\Pi_{\epsilon,s,n}$ factors through $\pi_{\epsilon,s,n}$ when $s \notin \frac{1}{2}\mathbb{Z}$. For $s \in (\frac{1}{2}\mathbb{Z})\setminus \mathbb{Z}$, we can also define the period, which may be zero, on $\Pi_{\epsilon,s,n}$ since $\psi_{\epsilon,s}$ is holomorphic on $\mathbb{C} \setminus \mathbb{Z}$, and the period also factors through $\pi_{\epsilon,s,n}$ by continuity of $\psi_{\epsilon,s}$ on $s$.   We will show the period does not vanish at the distinguished minimal $K$-type below, and so is a non-zero period.

\subsubsection{Non-vanishing on the distinguished minimal $K$-type}\label{nonvanish-unitary}

\begin{proof}[Proof of Theorem \ref{Non} in the unitary case excluding \eqref{exception}:]

As in the generic case, given any $G_n'$-distinguished irreducible unitary representation $\pi$ of $G_n$, $\pi$ contains a unique distinguished minimal $K$-type.

By parabolic induction, we only need to show that in the following cases (i) and (ii), the period does not vanish on the distinguished minimal $K$-type.

\begin{itemize}
    \item[(i)] For $\pi_{0,s,n}$ with $s \in \mathbb{C} \setminus \mathbb{Z}$, the distinguished minimal $K$-type is $(\underbrace{0,\dots,0}_{2n})$.
\end{itemize}

When $n=1$, for $\pi_{0,s,1}$ with $s \in \mathbb{C} \setminus \mathbb{Z}$, by the generic case Section~\ref{generic case non-vanishing}, the period map is a holomorphic map and does not vanish at the distinguished minimal $K$-type, which is the trivial $K$-type.

For $n > 1$, the intertwining operator $\iota_{\widetilde{w}_3,0,s}$ does not vanish at the distinguished minimal $K$-type. By the parabolic induction theorem~\ref{induction}, the period of $\Pi_{0,s,n}$ also does not vanish at the distinguished minimal $K$-type.

\begin{itemize}
    \item[(ii)] For $\pi_{1,s,n}$ with $s \in \mathbb{C} \setminus \mathbb{Z}$, the distinguished minimal $K$-type is $(\underbrace{2,\dots,2}_{n}, \underbrace{0,\dots,0}_{n})$.
\end{itemize}

When $n=1$, for $\pi_{1,s,1}$, by the generic case Section~\ref{generic case non-vanishing}, the period map is a holomorphic map and is identity at the distinguished minimal $K$-type $(2,0)$. To proceed by induction on $n$, we define the period in another way different from the $\psi_{\epsilon,s}$ above, but these two periods coincide with each other by Lemma \ref{mult<=1} and continuity on $s$. Let $\widetilde{w}_4$ be the longest element in $S_{n-1}$, the group of permutations of the set $\{n+2,\dots,2n\}$. Consider 
\begin{align*}
    \widetilde{\phi}_{1,s}: \mathrm{Ind}_{B}^G(\lambda_{1,s}) &\xrightarrow{\iota_{\widetilde{w}_2,1,s}\circ \iota_{\widetilde{w}_1,1,s}}  \mathrm{Ind}_{B}^G(\lambda_{n+1} \boxtimes \lambda_n \boxtimes \lambda_{n-1} \boxtimes \dots \boxtimes \lambda_1 \boxtimes \lambda_{n+2} \boxtimes \dots \boxtimes \lambda_{2n})\\& \xrightarrow{\iota_{\widetilde{w}_4,1,s}} \mathrm{Ind}_{B}^G\left((\lambda_{n+1} \boxtimes \lambda_{n}) \boxtimes (\lambda_{n-1} \boxtimes \dots\boxtimes \lambda_1) \boxtimes  (\lambda_{2n} \boxtimes\dots\boxtimes \lambda_{n+2})\right).
\end{align*} 
Since the composition of intertwining operators is independent of the order, one can realize $\widetilde{\phi}_{1,s}$ in another way: 
\begin{align*}
    \mathrm{Ind}_{B}^G(\lambda_{1,s}) &\to \mathrm{Ind}_{B}^G(\lambda_{n+1} \boxtimes \lambda_1\boxtimes \dots \boxtimes \lambda_n \boxtimes \lambda_{n+2} \boxtimes \dots \boxtimes \lambda_{2n})\\&\to \mathrm{Ind}_{B}^G(\lambda_{n+1} \boxtimes \lambda_n \boxtimes \lambda_1\boxtimes \dots \boxtimes  \lambda_{n-1}\boxtimes \lambda_{n+2} \boxtimes \dots \boxtimes \lambda_{2n})\\&\to \mathrm{Ind}_{B}^G(\lambda_{n+1} \boxtimes \lambda_n \boxtimes \lambda_{n-1}\boxtimes\dots \boxtimes \lambda_1\boxtimes \lambda_{2n} \boxtimes \dots \boxtimes \lambda_{n+2}),
\end{align*}
where the last step is change the subscripts of $\lambda$: $(1,2,\dots,n-1)$ (resp. $(n+2,n+3,\dots,2n)$) to  $(n-1,n-2,\dots,1)$ (resp. $(2n,2n-1,\dots,n+2)$), and the image of the last step is contained in $\pi_{1,n-1-s,1}\times \pi_{1,s+1,n-1}$. So is the image of $\widetilde{\phi}_{1,s}$.

By the results of generic cases, the $G_2'$-period,  denoted by $\mathcal{P}_1$, of $\pi_{1,n-1-s,1}$ does not vanish at the distinguished minimal $K$-type $(2,0)$. By induction on $n$, the $G_{2n-2}'$-period, denoted by $\mathcal{P}_2$, of $\pi_{1,s+1,n-1}$ does not vanish at the distinguished minimal $K$-type $(\underbrace{2,\dots,2}_{n-1},\underbrace{0,\dots,0}_{n-1})$.  By parabolic induction, see Theorem \ref{induction}, we get a $G_{2n}'$-period of $\Pi_{\epsilon,s,n}$, that is
\[ (\mathcal{P}_1\times \mathcal{P}_2)\circ \widetilde{\phi}_{1,s}.\]
 
To show this $G_{2n}'$-period does not vanish on the $K$-type $\sigma =(\underbrace{2,\dots,2}_n,\underbrace{0,\dots,0}_n)$, it suffices to show that the image of the $\sigma$-isotypic space of  $\Pi_{1,s,n}$ under $\widetilde{\phi}_{1,s}$ has non-zero component in the $K$-type induced from $\mathrm{U}(2)\times \mathrm{U}(2n-2)$-type $(2,0)\boxtimes (2,\dots,2,0,\dots,0)$ of $\pi_{1,n-1-s,1}\times \pi_{1,s+1,n-1}$. This will be proved in Appendix B, and the proof is now finished.
\end{proof}

Here, we provide an example to illustrate the proof above.
\begin{example}
Consider the $G_4'$-distinguished representation $\pi=(\frac{\det}{|\det|}\cdot |\det|^s)_{G_2}\times (\frac{\det}{|\det|}\cdot |\det|^{-s})_{G_2}$ ($s\notin \mathbb{Z}$) of $G_4$. One can realize $\pi$ as a submodule of 
\[\Pi:=(\frac{\det}{|\det|}\cdot |\det|^{-1+s})_{G_1}\times (\frac{\det}{|\det|}\cdot |\det|^{1+s})_{G_1}\times (\frac{\det}{|\det|}\cdot |\det|^{-1-s})_{G_1}\times  (\frac{\det}{|\det|}\cdot |\det|^{1-s})_{G_1}.\]
 There exists an obvious period of $\Pi$, which arises from the parabolic induction of the period $\mathcal{P}_0$ of $(\frac{\det}{|\det|}\cdot |\det|^{-1})_{G_1}\times  (\frac{\det}{|\det|}\cdot |\det|)_{G_1}$. However, the period $\mathcal{P}_0\times \mathcal{P}_0$ vanishes on $\pi$, since $\mathcal{P}_0$ vanishes on the $G_2$-submodule $(\frac{\det}{|\det|})_{G_2}$. Thus, one cannot obtain the period in this manner. The above usage of the intertwining operator resolves this issue in an elegant way.
\end{example}

\begin{remark}
The proof of Theorem~\ref{Non} in the unitary case actually shows slightly more: for $(\frac{\det}{|\det|})_{G_m}\times (\frac{\det}{|\det|})_{G_m}$ with $m\leq 4$, the period also does not vanish on the distinguished minimal $K$-type. However, when $m\geq 5$, the argument fails because the period in case (ii) of Section~\ref{generic case non-vanishing} vanishes for $s \in -4 - 2\mathbb{Z}_{\geq 0}$. These exceptional cases are already included in \eqref{exception}.
\end{remark}

\section{Local factors of distinguished representations}\label{local factor}
To each finite-dimensional semisimple complex representation $\lambda$ of the Weil group of a local field, one can associate a local $L$-factor as well as a local $\varepsilon$-factor with certain nice properties (see \cite{Jac} for more details). Let $\psi^{\circ}$ be the additive character of $\mathbb{C}$ defined by
\[z \mapsto \text{exp} \left( 2 \pi \sqrt{-1} (z + \overline{z}) \right),\]
where $\pi$ denotes the standard irrational number. Then for any non-trivial additive unitary character of $\mathbb{C}$, there exists $b \in \mathbb{C}^\times$ such that it is of the form
\[\psi_{b}: ~ z \mapsto \psi^{\circ}(bz) \text{ for } z \in \mathbb{C}.\]
We have $\psi_{b}\big|_{\mathbb{R}}=1 \Longleftrightarrow b \in \sqrt{-1} \mathbb{R}^{\times}$. We assume such standard characters and denote them also by $\psi$ in this section.

Let $\mathbb{F}$ be a local field and $\psi$ be a non-trivial character of $\mathbb{F}$. To each pair of irreducible smooth representations (Casselman-Wallach representations in the archimedean place) $\pi$ and $\pi^\prime$ of $\mathrm{GL}_n(\mathbb{F})$  and $\mathrm{GL}_m(\mathbb{F})$ respectively, one can attach local Rankin–Selberg L-factor $L(s, \pi \times \pi^\prime)$, $\varepsilon$-factor $\varepsilon(s, \pi \times \pi^\prime,\psi)$ and $\gamma$-factor $\gamma(s, \pi \times \pi^\prime,\psi)$ (see \cite{JPSS} for more details) and we have 
\begin{equation}\label{eq:A}
\gamma(s, \pi \times \pi^\prime,\psi)= \varepsilon(s, \pi \times \pi^\prime,\psi) \frac{L(1-s, \pi^{\vee} \times \pi^{\prime \vee})}{L(s, \pi \times \pi^\prime)}
\end{equation}
where ${\pi}^{\vee}$ denotes the contragredient of $\pi$. 

The local $\varepsilon$-factor corresponding to a one dimensional representation $\varkappa_{m,t}: z \mapsto \left(\frac{z}{|z|}\right)^m |z|^{2t}$ of $\mathbb{C}^\times$ for some $m \in \mathbb{Z} , t \in \mathbb{C}$ and the additive character $\psi^{\circ}$ is given by 
\begin{equation*}
    \varepsilon(s,\varkappa_{m,t}, \psi^{\circ})= (\sqrt{-1})^{|m|},
\end{equation*}
for each $s \in \mathbb{C}$. The \textbf{local $\varepsilon$-factor} corresponding to the character $\varkappa_{m,t}$ of $\mathbb{C}^\times$ and the additive character $\psi_{b}$ is given by (cf. \cite[Remark 4.3.6]{Jac})
\begin{align}
    \varepsilon(s,\varkappa_{m,t}, \psi_b) &= \varkappa_{m,t}(b) |b|_{\mathbb{C}}^{s-\frac{1}{2}} ~ \varepsilon(s,\varkappa_{m,t}, \psi^{\circ}) \nonumber\\
    &= (\sqrt{-1})^{|m|} ~ b^m ~|b|^{(2t-m+2s-1)}, \label{new:A}
\end{align}
where $|b|_{\mathbb{C}}=|b|^2$ is the normalized absolute value on $\mathbb{C}$, for each $s \in \mathbb{C}$. Therefore, if $\psi_b \big|_\mathbb{R}=1$, we have $b\in \sqrt{-1}\mathbb{R}^{\times}$ and by \eqref{new:A},
\begin{align}\label{new:B}
    \varepsilon\left(\frac{1}{2},\varkappa_{m,t}, \psi_b\right)
    =(\sqrt{-1})^{|m|}\left(\frac{b}{|b|}\right)^m |b|^{2t}.
\end{align}
In general, if $\pi=Q(\lambda_1,\dots ,\lambda_n) $ is the unique irreducible quotient of $ \mathrm{Ind}^{G_n}_{B_n}
(\lambda_1 \otimes \cdots\otimes \lambda_n)$, then
\begin{equation}\label{eq:B}
 \varepsilon(s,\pi, \psi)=\varepsilon\left(s,\mathrm{Ind}^{G_n}_{B_n}
(\lambda_1 \otimes \cdots\otimes \lambda_n), \psi\right)=\prod\limits_{j=1}^n \varepsilon\left(s,\lambda_j, \psi\right).   
\end{equation}
We now have the following $\varepsilon$-factor result of distinguished representation  similar to \cite[Theorem 1.1(2)]{ST} for the Archimedean case and \cite[Theorem 0.1]{Off},\cite[Theorem 1.1]{MO} for the non-archimedean case:
\begin{theorem}\label{thm:epsilon}
 Let $\pi$ be an irreducible $\mathrm{GL}_n(\mathbb{R})$-distinguished representation of $\mathrm{GL}_n(\mathbb{C})$,  and $\psi$ be an additive character of $\mathbb{C}$ such that $\psi\big|_\mathbb{R}=1$. Then
 \[\varepsilon\left(\frac{1}{2}, \pi , \psi\right)=1.\]  
\end{theorem}
\begin{proof}
 Let $\pi=Q(\lambda_1,\dots ,\lambda_n)$ for some character $\lambda_j=\varkappa_{m_j,t_j}$ with $m_j \in \mathbb{Z}$ and $t_j \in \mathbb{C}$. Since $\pi$ is $\mathrm{GL}_n(\mathbb{R})$-distinguished, there exists an involution $w \in S_n$  such that $\lambda_{w(j)}=\overline{\lambda}_j^{-1}$ and when $w(j)=j$, we have $\lambda_j(-1)=1$. Define $A=\{j \mid j < w(j) \}$ and $B=\{j \mid j = w(j) \}$. If $w(j)=j$, $\lambda_{w(j)}=\overline{\lambda}_j^{-1} \implies t_j=0$ and $\lambda_j(-1)=1 \implies m_j=2k_j $ an even integer. Since $b\in \sqrt{-1}\mathbb{R}^{\times}$, we have $\left(\frac{b}{|b|}\right)^2=-1$, and hence, by \eqref{new:B}
 \begin{align*}
 \varepsilon\left(\frac{1}{2},\varkappa_{2k_j,0},\psi_b\right)
 &=(\sqrt{-1})^{2|k_j|}\left(\frac{b}{|b|}\right)^{2k_j}\\
 &=(-1)^{|k_j|+k_j}=1.
 \end{align*}
 Therefore, by relation (\ref{eq:B}) 
 \begin{align*}
     \varepsilon\left(\frac{1}{2},\pi, \psi\right)&=\prod\limits_{j\in A} \varepsilon\left(\frac{1}{2},\lambda_j, \psi \right)\varepsilon\left(\frac{1}{2},\lambda_{w(j)}, \psi \right)  \times \prod\limits_{j\in B} \varepsilon\left(\frac{1}{2},\lambda_j, \psi \right)\\
     &= \prod\limits_{j\in A} \varepsilon\left(\frac{1}{2},\varkappa_{m_j,t_j}, \psi \right)\varepsilon\left(\frac{1}{2},\overline{\varkappa_{m_j,t_j}}^{-1}, \psi \right)  \times \prod\limits_{j\in B} \varepsilon\left(\frac{1}{2},\varkappa_{2k_j,0}, \psi \right)\\
     &= \prod\limits_{j\in A} \varepsilon\left(\frac{1}{2},\varkappa_{m_j,t_j}, \psi \right)\varepsilon\left(\frac{1}{2},\varkappa_{m_j,-t_j}, \psi \right), \text{ as } \varepsilon\left(\frac{1}{2},\varkappa_{2k_j,0}, \psi \right)=1\\
     &= \prod\limits_{j\in A} (\sqrt{-1})^{2|m_j|}\left(\frac{b}{|b|}\right)^{2m_j}|b|^{2t_j}|b|^{-2t_j}\\
     &= \prod\limits_{j\in A} (-1)^{|m_j|+m_j}\\
     &= 1.
 \end{align*}
\end{proof}

\begin{corollary}(cf. \cite[Theorem 6.3]{MO})\label{triviality of ep}
 Let $\pi$ be an irreducible $\mathrm{GL}_n(\mathbb{R})$-distinguished representation of $\mathrm{GL}_n(\mathbb{C})$, $\pi'$ be an irreducible $\mathrm{GL}_m(\mathbb{R})$-distinguished representation of $\mathrm{GL}_m(\mathbb{C})$. Let $\psi$ be an additive character of $\mathbb{C}$ with a trivial restriction to $\mathbb{R}$. Then
 \[\varepsilon\left(\frac{1}{2}, \pi \times \pi' , \psi\right)=1.\]
\end{corollary}
\begin{proof}
    Let $\pi=Q(\lambda_1,\ldots,\lambda_n)$ and $\pi'=Q(\mu_1,\ldots,\mu_m)$ for characters $\lambda_i, \mu_j$ of $\mathbb{C}^\times$. As $\pi$ and $\pi'$ are both distinguished, by \cite[Theorem 1.3]{Kem15a}, there exist involutions $w\in S_n$ and $w^\prime\in S_m$ such that
    \[
    \lambda_{w(i)}=\overline{\lambda}_i^{-1},\qquad
    \mu_{w^\prime(j)}=\overline{\mu}_j^{-1},
    \]
    and for every fixed point of $w$ or $w^\prime$, the corresponding character takes the value $1$ at $-1$. By the multiplicativity of the $\varepsilon$-factor, we have
    \[
    \varepsilon\left(\frac12,\pi\times\pi',\psi\right)
    =\prod_{i=1}^n\prod_{j=1}^m
    \varepsilon\left(\frac12,\lambda_i\mu_j,\psi\right).
    \]
    Define an involution $\sigma$ on the set $\{1,\ldots,n\}\times\{1,\ldots,m\}$ as follows: for $1 \le i \le n$ and $1 \le j \le m$
    \[
    \sigma(i,j)=(w(i),w^\prime(j)).
    \]
    Then
    \[
    \lambda_{w(i)}\mu_{w^\prime(j)}
    =\overline{\lambda_i}^{-1}\overline{\mu_j}^{-1}
    =\overline{\lambda_i\mu_j}^{-1}.
    \]
    Moreover, if $\sigma(i,j)=(i,j)$, then $w(i)=i$ and $w^\prime(j)=j$, and hence
    \[
    (\lambda_i\mu_j)(-1)=\lambda_i(-1)\mu_j(-1)=1.
    \]
    Thus, the characters $\lambda_i\mu_j$ satisfy the same involution and fixed-point conditions as in the proof of Theorem \ref{thm:epsilon}. Applying the same computation to the above product produces
    \[
    \varepsilon\left(\frac12,\pi \times \pi',\psi\right)=1.
    \]
\end{proof}
\begin{remark}
After our first preprint, Nadir Matringe pointed out that the above result already occurred in their paper \cite[Theorem 6.3]{MO}. We still keep the result here as our approach is a little different.   
\end{remark}
\begin{remark}
We want to mention the following $\gamma$-factor characterization due to Kemarsky. Let $\pi$ be an irreducible, generic, unitary representation $\mathrm{GL}_n(\mathbb{C})$. Then 
\[\pi\text{ is }\mathrm{GL}_n(\mathbb{R})\text{-distinguished if and only if }\gamma\left(\frac{1}{2}, \pi \times \chi, \psi\right)=1\]
for every unitary character $\chi:\mathbb{C}^\times \longrightarrow \mathbb{C}^\times$ such that $\chi \big|_{\mathbb{R}^\times}=1$. The forward implication follows from \cite[Theorem 1.1]{Kem15a}, and the converse follows from \cite[Theorem B.3]{Kem15a} together with Theorem \ref{main result_generic}.
\end{remark}

\section{Further application and problem}\label{application}
\subsection{Transferring branching laws: Theta correspondence}
In this subsection, $\mathbb{F}$ is a local field of characteristic 0. It is well-known to experts that theta correspondence can transfer branching laws. Here we state a direct corollary of the main result implying some branching laws $\mathrm{GL}_{2m}(\mathbb{R})\downarrow\mathrm{GL}_m(\mathbb{C})$. We state the idea in a general setting. We only consider type II dual pair. Let $D$ be a central division algebra with center $\mathbb{F}$. Let $V,W$ be left $D$-vector spaces with dimension $a,n$ respectively. Set:
$$
    \mathbb{X}:=V\otimes W \qquad \mathbb{W}:=\mathbb{X}\oplus \mathbb{X}^*.
$$
Then $\mathbb{W}$ is equipped with a natural symplectic form and
\begin{equation*}
    \mathrm{GL}_D(V)\times \mathrm{GL}_D(W) \hookrightarrow \mathrm{M}(\mathbb{X})\hookrightarrow \mathrm{Sp}(\mathbb{W}),
\end{equation*}
where $\mathrm{M}(\mathbb{X})$ is the Levi subgroup of Siegel parabolic that stabilizes $\mathbb{X}$. Since $\mathrm{M}(\mathbb{X})$ splits in the metaplectic group $\mathrm{Mp}(\mathbb{W})$, we have the splitting
\begin{equation}
    \mathrm{GL}_D(V)\times \mathrm{GL}_D(W) \hookrightarrow \mathrm{Mp}(\mathbb{W}).
\end{equation}
Now let $\mathbb{K}$ be a quadratic extension over $\mathbb{F}$, $D_{\mathbb{K}}:=D\otimes_{\mathbb{F}}\mathbb{K}$, and let $V'$ be a left $D_{\mathbb{K}}$-vector space with dimension $m$. Moreover, let $W_{\mathbb{K}}:=D_{\mathbb{K}}\otimes_D W$. Then we have \textbf{seesaw pair}:
\[
\begin{tikzpicture}
  \matrix[matrix of math nodes,row sep=2em,column sep=4em,minimum width=2em] (m) {
     \mathrm{GL}_{D}(V') & \mathrm{GL}_{D_\mathbb{K}}(W_\mathbb{K}) \\
     \mathrm{GL}_{D_\mathbb{K}}(V')   & \mathrm{GL}_D(W) \\};
  \draw[-] (m-1-1) -- (m-2-1); 
  \draw[-] (m-1-2) -- (m-2-2); 
  \draw[-] (m-1-1) -- (m-2-2); 
  \draw[-] (m-1-2) -- (m-2-1); 
\end{tikzpicture}
\]
and assume that $\pi_1$ (resp. $\pi_2$) is a moderate growth smooth Fr\'echet representation (or smooth representation in non-archimedean case) of $\mathrm{GL}_{D_\mathbb{K}}(V')$ (resp. $\mathrm{GL}_D(W)$). Consequently, we have the following isomorphism\cite{Pra96}
	\begin{equation}\label{seesaw}
	\Hom_{\mathrm{GL}_{D_\mathbb{K}}(V')}(\Theta(\pi_2),\pi_1)\cong \Hom_{\mathrm{GL}_D(W)}(\Theta(\pi_1),\pi_2). 
	\end{equation} 
In our case, $D=\mathbb{F}$, $D_{\mathbb{K}}=\mathbb{K}$ and $\pi_2$ is the trivial representation. Hence we get the following corollary:
\begin{corollary}\label{theta}
   For an irreducible representation $\pi$ of $\mathrm{GL}_m(\mathbb{K})$, $\Hom_{\mathrm{GL}_{m}(\mathbb{K})}(\Theta(\mathbb{C}),\pi)\neq 0$ if and only if $\Theta(\pi)$ is distinguished.
\end{corollary}
For concrete application of Corollary~\ref{theta}, it is favorable to calculate the big theta lift. Unfortunately, compared to the small theta lift, we know much less about the big theta lift. Hence, we specialize to some cases, where $\theta=\Theta$. Assume that $m=n$. By relating problem with poles of Godement-Jacquet $L$-function, \cite{FSX} shows that for generic representation $\pi$ of $\mathrm{GL}_n(\mathbb{K})$, $\Theta(\pi)=\theta(\pi)=\pi^{\vee}$. On the other hand, we need to calculate the big theta lift of the trivial representation of $\mathrm{GL}_{n}(\mathbb{F})$. This follows from a more general result:
\begin{theorem}\label{Theta}
    Let $\mathbb{F}$ be a local field of characteristic 0. And let $(G,G')=(\mathrm{GL}_n(\mathbb{F}),\mathrm{GL}_m(\mathbb{F}))$ be a reductive dual pair in the stable range($G$ is smaller), which means $m\geq 2n$. Then for irreducible unitary representation $\pi$ of $G$, $\Theta(\pi)=\theta(\pi)$ is an irreducible unitary representation of $G'$.
\end{theorem}
\begin{proof}
    If $\mathbb{F}$ is a non-archimedean local field, it is proved by~\cite{CLLTZ25} (actually, the strategy in \cite{CZ24} can already prove it). For $\mathbb{F}$ an Archimedean field, it is explained in the Appendix B.
\end{proof}
When $\pi$ is the trivial representation and $m=2n$, we have $\theta(\pi)=\Ind_{P(\mathbb{F})}^{G(\mathbb{F})}(\mathbb{1}\boxtimes \mathbb{1})$, where $P(\mathbb{F})$ is the parabolic subgroup 
\[P(\mathbb{F})=\left\{\begin{pmatrix}
A & B \\ 0 & C
\end{pmatrix}\ |\ A,B,C\in M_{n\times n}(\mathbb{F}), \det(A),\det(C)\neq 0\right\}.\] 
Consequently, we get Corollary~\ref{new-branching}. Similarly, by Theorem~\ref{Theta}, we can obtain a result like the above corollary in the $p$-adic case using the classification given by \cite{Mat14}. 

\subsection{Harmonic analysis on symmetric space $\mathrm{GL}_n(\mathbb{R})\backslash \mathrm{GL}_n(\mathbb{C})$}
In this subsection, let $X=\mathrm{GL}_n(\mathbb{R})\backslash \mathrm{GL}_n(\mathbb{C})$. One of the most important questions in harmonic analysis on $X$ is to study the Plancherel decomposition of $L^2(X)$. In particular, we want to determine the support of the Plancherel measure of $L^2(X)$. Some results were obtained by \cite{BS1} and \cite{BS2}. However, we want to pursue a more involved answer, describing the support of the Plancherel measure in terms of $L$-parameters or $A$-parameters of $\mathrm{GL}_n(\mathbb{C})$. We note that 
Sakellaridis and Venkatesh have proposed conjectures in terms of the dual group of $X$ (cf. \cite[Conjecture 1.3.1]{SV}). In our case, it is easy to see $\rho_{\mathfrak{g}_n'}\leq \rho_{\mathfrak{g}_n/\mathfrak{g}_n'}$, where $\rho$ is the moment as a representation of $\mathfrak{g}_n'$. Consequently, by  \cite[Theorem 4.14]{Ko}, we have that $L^2(X)$ is a tempered representation of $\mathrm{GL}_n(\mathbb{C})$. Hence, the support of the Plancherel measure is contained in the set of tempered distinguished representations.
\begin{question}
According to Theorem~\ref{main result_unitary}, we have the description of the set of tempered distinguished representations in terms of $L$-parameters. Can we figure out the subset of the support of the Plancherel measure of $X$?     
\end{question}
There is another question closely related to the harmonic analysis on symmetric spaces. The notation used here are the same as in Section~\ref{derivative section}. Let $\pi$ be an irreducible unitary distinguished representation. Denote the image of $\pi \hookrightarrow \mathcal{C}^{\infty}(X)$ by $\mathcal{C}_{\pi}$. Then we can formally define a model transfer, which is important in the relative trace formula:
\begin{align*}
    &\mathcal{C}_{\pi} \longrightarrow \mathcal{W}_r(\pi,\psi) \\
    f\mapsto (g&\mapsto \int_{N_n'\backslash N_n} f(ng)(\psi_r^{\sigma})^{-1}(n)dn).
\end{align*}
\begin{question}
    Can we specify the $L$-parameter of $\pi$ such that the model transfer integral is convergent?
\end{question}
This question has a partial answer in the generic case, see \cite{Kem15a}.

\appendix
\section{Mellin Transform}

In this appendix, we generalize the method of the Mellin transform from \cite{JS1} to apply to the proof of Theorem~\ref{thm7.2}. Most of the generalizations follow the approach in \cite{JS1}, and we refer to \textit{loc. cit.} for the proof.

The following notation are fixed in this appendix. Let $\mathbb{F}$ be an Archimedean local field, i.e., $\mathbb{F}=\mathbb{R}$ or $\mathbb{C}$. Define $|x|_\mathbb{F}:=|x|$ if $\mathbb{F}=\mathbb{R}$, and $|x|_\mathbb{F}:=|x|^2$ if $\mathbb{F}=\mathbb{C}$. Let $V$ be a Fréchet space. Let $T$ be a locally compact topological group. Our prototypes of $T$ include $(\mathbb{R}^{+})^n$, $(\mathbb{R}^{\times})^n$, or $(\mathbb{C}^{\times})^n$ (we call these three cases the $*$-class). We use $\chi$ to denote a character of $T$. If $T$ belongs to the $*$-class, we decompose $\chi$ as $\chi:=\chi_0 |t|^s$, where $\chi_0=(\frac{t}{|t|})^m$, $m\in\mathbb{Z}$, is a unitary character, $t\in T$, and $s$ is a real number. We also embed the $*$-class into the group $X=\mathbb{F}^n$.

Now we formally define the Mellin transform for $f\in \mathcal{C}(T, V)$: \begin{equation*} \hat{f}(\chi)=\int_T f(t)\chi(t) dt. \end{equation*} If $T$ belongs to the $*$-class, we also define: \begin{equation}\label{eq 18} \hat{f}(\chi, n)=\int_T f(t)\chi(t)(\log|t|_{\mathbb{F}})^n dt, \end{equation} where $n$ is a positive integer. Note that when $f\in \mathcal{S}(X, V)$, the above integral is absolutely convergent for $\mathrm{Re}(\chi)>0$.

\begin{example} The most important example is the Gamma function, which appears in the explicit formula for Whittaker functions on $\mathrm{GL}_2(\mathbb{C})$: 
\begin{equation*} 
\Gamma(s):=\int\limits_{0}^{+\infty} e^{-t} t^s \frac{dt}{t}. 
\end{equation*} It satisfies the functional equation $\Gamma(s+1)=s\Gamma(s)$, and we use this functional equation to provide the meromorphic extension of the Gamma function to the complex plane. \end{example}
From now on, we assume that $T$ is in $*$-class and $f\in \mathcal{S}(X,V)$. Then $\hat{f}(\chi,n)$ also has a functional equation, hence a meromorphic extension. We note that $\hat{f}(s,m;n)$ has poles of order(at most) $kn$ at points $-k-m$, where $k$ is a natural number and $m=0,1$, if $\mathbb{F}=\mathbb{R}$; $\hat{f}(s,m;n)$ has poles of order(at most) $kn$ at points $-k-\frac{|m|}{2}$, where $k$ is a natural number and $m\in\mathbb{Z}$, if $\mathbb{F}=\mathbb{C}$.\\
Following the idea in \cite{JS1}, Section {\bf3.4}, we define an appropriate function space on which the Mellin transform of the degenerate Whittaker functions will fall. Recall that for finite-dimensional semi-simple representation $\sigma$ of $\mathbb{F}^{\times}$, we can define $L$-function $L(\chi,\sigma)$ on the character group, which are pieces of the complex plane. 
\begin{definition}
    Let $\{\sigma_i\}_{i=1}^{n}$ be a set of finite-dimensional semi-simple representations of $\mathbb{F}^{\times}$. Let $\mathcal{M}(\sigma_1,\dots,\sigma_n)$ be the subspace of meromorphic functions generated by the functions of the form
    \begin{equation}
        M(\chi_1,\dots,\chi_n)=\prod_{i=1}^nL(\chi_i,\sigma_i)h(\chi_1,\dots,\chi_n)
    \end{equation}
    where $h$ is a holomorphic function on the character group of $T$ and valued in $V$. Moreover, $M(\chi)$ should satisfy some rapid decay condition along the vertical strip, which we refer to \textit{loc.cit.} for precise definition.
\end{definition}

\begin{proposition}
    Let $\{\sigma_i\}_{i=1}^{n}$ be a set of finite-dimensional semi-simple representations of $\mathbb{F}^{\times}$. Then any element of $\mathcal{M}(\sigma_1,\dots,\sigma_n)$ is the Mellin transform of a function $f$ of the form:
    \begin{equation}\label{eq21}
        f(t)=\sum_{j=1}^r \phi_j(t) \xi_j(t),
    \end{equation}
where $\phi_j\in\mathcal{S}(X,V)$ and $\xi_j$ is a left-translate finite function,i.e., $\xi(t)=\chi(t)(\log|t|_{\mathbb{F}})^q$ for some character $\chi$ and natural number $q$.
\end{proposition}
\begin{proof}
    When $V=\mathbb{C}$, it is Proposition 4 in \textit{loc.cit.} The idea is seeking some function $f$ in the form of Equation~\eqref{eq21} such that after the Mellin transform it will have the prescribed principal part at some point. This is achieved by solving ODEs. Now, suppose the principal part at $(a_1,\dots ,a_n)$ is given by $\sum_{0< (l_1,\dots,l_n)\leq N }\frac{v_{(l_1\dots l_n)}}{(z_i-a_i)^{l_i}}$ for some finite $n$-tuple $N$ and $v_{(l_1\dots l_n)}\in V$. Then consider the finite-dimensional subspace of $V$ generated by $v_{(l_1,\dots,l_n)}$. The ODEs can be solved within this subspace by the Borel Lemma.
\end{proof}
 We will need some estimates of $\chi$ in Equation~\eqref{eq21}.
\begin{proposition}\label{prop8.2}
   Let $\{\sigma_i\}_{i=1}^{n}$ be a set of finite-dimensional semi-simple representations of $\mathbb{F}^{\times}$. Let $f$ be a function whose Mellin transform $\hat{f}$ is in $\mathcal{M}(\sigma_1,\dots,\sigma_n)$. Suppose $\hat{f}$ is holomorphic in the product of the multi-half-planes $\mathrm{Re} \chi_i\geq 0$ for $1\leq i\leq n$. Then $f$ can be written in the form:
   \begin{equation}\label{eq22}
        f(t)=\sum_{j=1}^r \phi_j(t) \xi_j(t),
    \end{equation}
where $\phi_j\in\mathcal{S}(X,V)$ and $\xi_j$ is a left-translate finite function with $\mathrm{Re} \xi_j>0$ for $1\leq j\leq r$.
\end{proposition}
\begin{proof}
    When $V=\mathbb{C}$, it is Proposition 2 in \textit{loc.cit.} In general, let $\Omega$ be a complex domain (may not connected). We endow $H(\Omega,V)$, the holomorphic functions on it, with the subspace topology of $\mathcal{C}^{\infty}(\Omega,V)$. Then $H(\Omega,V)$ is a nuclear space and we have $H(\Omega,V)\simeq H(\Omega)\widehat{\otimes} V$, see \cite{Tr} chapter 44. Hence the Cauchy integral formula and maximum principle also hold for $H(\Omega,V)$.
\end{proof}

Now suppose $V$ has a continuous inner product.
\begin{proposition}\label{prop8.3}
    Let $\{\sigma_i\}_{i=1}^{n}$ be a set of finite-dimensional semi-simple representations of $\mathbb{F}^{\times}$. Let $f$ be a function whose Mellin transform $\hat{f}$ is in $\mathcal{M}(\sigma_1,\dots,\sigma_n)$. Suppose further that $f$ is square-integrable on $T$. Then $\hat{f}$ is holomorphic in the multi-half-planes $\mathrm{Re} \chi_i\geq 0$ for $1\leq i\leq n$.
\end{proposition}
\begin{proof}
     When $V=\mathbb{C}$, it is Proposition 5 in \textit{loc.cit.}, and the proof here applies directly to general $V$.
\end{proof}

\begin{proposition}\label{prop8.4}
    Let $(\pi,V_{\pi})$ be a distinguished irreducible unitary representation of $G_n$. For $v\in V_{\pi}$, $W_v$ is defined as in Section~\ref{Necessary}. Then the Mellin transform of $W_v|_{A_{n-d+1}\dots A_{n-1}}$ falls in some $\mathcal{M}(\sigma_1,\dots,\sigma_{d-1})$.
\end{proposition}
\begin{proof}
    For generic $\pi$, it is proved in Chapter 4 of \cite{JS1}. And the proof here applies directly to general $V_{\pi}$.
\end{proof}

\section{Proof of the claim in subsection \ref{nonvanish-unitary}}\label{appendix B}

We refer the detailed computation of intertwining operators to \cite[Section 5]{WZ24}.

Recall what we need to prove is the following. 

\begin{lemma}
Consider the intertwining operator as the proof of Subsection \ref{nonvanish-unitary}:
\begin{align*}
    \widetilde{\phi}_{1,s}: \mathrm{Ind}_{B}^G(\lambda_{1,s})     &\xrightarrow{ \iota_{\widetilde{w}_1,1,s}}  \mathrm{Ind}_{B}^G(\lambda_n \boxtimes  \dots \boxtimes \lambda_1 \boxtimes  \lambda_{n+1} \boxtimes \dots \boxtimes \lambda_{2n})\\    &\xrightarrow{\iota_{\widetilde{w}_2,1,s}}  \mathrm{Ind}_{B}^G(\lambda_{n+1} \boxtimes \lambda_n \boxtimes \lambda_{n-1} \boxtimes \dots \boxtimes \lambda_1 \boxtimes \lambda_{n+2} \boxtimes \dots \boxtimes \lambda_{2n})\\& \xrightarrow{\iota_{\widetilde{w}_4,1,s}} \mathrm{Ind}_{B}^G\left((\lambda_{n+1} \boxtimes \lambda_{n}) \boxtimes (\lambda_{n-1} \boxtimes \dots\boxtimes \lambda_1) \boxtimes  (\lambda_{2n} \boxtimes\dots\boxtimes \lambda_{n+2})\right),
\end{align*} 
then the image of the $\sigma = (\underbrace{2,\dots,2}_n, \underbrace{0,\dots,0}_n)$-isotypic subspace of $\Pi_{1,s,n}=\mathrm{Ind}_{B}^G(\lambda_{1,s})$ has non-zero component on the $\sigma$-subspace which is contained in the parabolic induction of the $\mathrm{U}(2)\times \mathrm{U}(2n-2)$-type $(2,0)\boxtimes (\underbrace{2,\dots,2}_{n-1},\underbrace{0,\dots,0}_{n-1})$.
\end{lemma}
\begin{proof}
Let $\langle \cdot ,\cdot\rangle$ denote the $K$-invariant inner product on $\sigma$. The $\sigma$-isotypic space of $\Pi_{1,s,n}=\mathrm{Ind}_{B}^G(\lambda_{1,s})$ is spanned by matrix coefficients
\[
\{k \mapsto \langle k \cdot u, v \rangle \mid u \in \sigma, v \text{ in the } (1,\dots,1)\text{-weight space of } \sigma\}.
\]
We will simply denote the matrix coefficient by $\langle k \cdot u, v \rangle$. Let $v_0$ be the non-zero vector in the $\mathrm{U}(n) \times \mathrm{U}(n)$-submodule with highest weight $(1,\dots,1) \boxtimes (1,\dots,1)$. Let us show that the image of $\langle k \cdot u, v_0 \rangle$ under $\widetilde{\phi}_{1,s}$ satisfies the lemma.

Notice that the image of $\langle k \cdot u, v_0 \rangle$ under $\iota_{\widetilde{w}_1,1,s}$ is a non-zero multiple of $\langle k \cdot u, v_0 \rangle$ by \cite[Lemma 5.2]{WZ24}.
Let $\langle k \cdot u, v_1 \rangle$ be the image of $\langle k \cdot u, v_0 \rangle$ under $\iota_{\widetilde{w}_2,1,s}$. Then $v_1$ is in the $\mathrm{U}(1)\times \mathrm{U}(n)\times \mathrm{U}(n-1)$-type $(1)\boxtimes (\underbrace{1,\dots,1}_n)\boxtimes (\underbrace{1,\dots,1}_{n-1})$  by \cite[Lemma 5.2]{WZ24}. By this lemma cited again, the image of $\langle k \cdot u, v_1 \rangle$ under $\iota_{\widetilde{w}_4,1,s}$
is a non-zero multiple of $\langle k \cdot u, v_1 \rangle$. 

Suppose that $v_1$ decomposes as $\sum_{i=0}^{n-1} \mu_i + \sum_{i=0}^{n-1} \delta_i$, where $\mu_i$ (resp. $\delta_i$) lies in the $\mathrm{U}(2) \times \mathrm{U}(2n-2)$-type
\[
(1,1) \boxtimes (\underbrace{2,\dots,2}_i, \underbrace{1,\dots,1}_{2n-2-2i}, \underbrace{0,\dots,0}_i) \quad \text{(resp. } (2,0) \boxtimes (\underbrace{2,\dots,2}_i, \underbrace{1,\dots,1}_{2n-2-2i}, \underbrace{0,\dots,0}_i))
\]
in $\sigma$. To show the statement, it suffices to prove that $\delta_{n-1}$ is non-zero.

Let us first show that $\sum \delta_i\neq 0$.
Assume that it is not true, then this image lies in the sum of the $\mathrm{U}(2) \times \mathrm{U}(2n-2)$-type $(1,1)\boxtimes (\dots)$. Since $v_1$ is in the $\mathrm{U}(1)\times \mathrm{U}(n)\times \mathrm{U}(n-1)$-type $(1)\boxtimes (1,\dots,1)\boxtimes (1,\dots,1)$, one gets that $v_1$ must belong to the $\mathrm{U}(n+1) \times \mathrm{U}(n-1)$-type $(1,\dots,1)\boxtimes (1,\dots,1)$ in $\sigma$. However, this is impossible since $\sigma$ contains no non-zero element of such a $\mathrm{U}(n+1) \times \mathrm{U}(n-1)$-type.

Notice that the $\mathrm{U}(2) \times \mathrm{U}(2n-2)$-type $(2,0)\boxtimes \gamma$ appears in $\sigma$ if and only if $\gamma = (\underbrace{2,\dots,2}_{n-1}, \underbrace{0,\dots,0}_{n-1})$, which means that $\delta_{n-1}\neq 0$. This completes the proof. 

\end{proof}

\section{Proof of Theorem~\ref{Theta} in archimedean case}

In the Archimedean case, experts have long believed that Theorem~\ref{Theta} is true, but the proof has never appeared in the literature. Therefore, this appendix does not contribute new results but rather explains this issue. It is worth noting that Loke and Ma proved the type I case in \cite{LM15}. Since the type II case also involves the Fock model and relevant results about compact dual pairs, their argument works equally well in the type II case if Jian-Shu Li's construction (cf. \cite{Li89}) holds for the type II case. Li has remarked in \emph{loc. cit.} that the type II case is easier. Thus, we address this gap by showing that Li's construction works in the type II case.

In this appendix, $\mathbb{F}$ is an Archimedean local field, and $(G, G')$ is the reductive dual pair $(\mathrm{GL}_n(\mathbb{F}), \mathrm{GL}_m(\mathbb{F}))$ in the stable range, that is, $m \geq 2n$. The unitary Weil representation $\omega$ of $G \times G'$ is realized as $L^2(M_{n \times m}(\mathbb{F}))$, with the action given by: \begin{equation} ((g, h) \cdot f)(x) = |\det g|^{-\frac{m}{2}}|\det h|^{\frac{n}{2}} f(g^{-1} \cdot x \cdot h), 
\end{equation} 
where $g \in G$, $h \in G'$, and when $\mathbb{F} = \mathbb{C}$, $|\cdot|$ denotes the square norm. Moreover, $\omega^{\infty} = \mathcal{S}(M_{n \times m}(\mathbb{F}))$. Now, given an irreducible unitary $SAF$ representation $(\sigma, V_{\sigma})$ of $G$, we can formally define a Hermitian form on $\omega \hat{\otimes} V_{\sigma}$ by: \begin{equation}\label{matrix coe} (\Phi, \Phi') := \int_G ((\omega^{\infty} \otimes \sigma)(g) \Phi, \Phi')  dg, 
\end{equation}
where $(\cdot, \cdot) := (\cdot, \cdot)_{\omega} \otimes (\cdot, \cdot)_{\sigma}$. 
Consider the radical of this form: 
\begin{equation*} \Rad((\cdot, \cdot)) := \{\Phi \in \omega^{\infty} \hat{\otimes} V_{\sigma} \mid (\Phi, \Phi') = 0, \forall \Phi' \in \omega^{\infty} \hat{\otimes} V_{\sigma} \}. 
\end{equation*} 
We aim to show: \begin{itemize} \item The integral in Equation~\eqref{matrix coe} is actually convergent and positive in the stable range; \item $H(\sigma) := (\omega^{\infty} \hat{\otimes} V_{\sigma^{\vee}}) / \Rad((\cdot, \cdot)) = \theta(\sigma)$. \end{itemize} The first point follows from the $L^p$-estimate of the matrix coefficient.
\begin{proposition}
    For $\phi,\phi'\in \mathcal{S}(M_{n\times m}(\mathbb{F}))$, we have 
    \begin{equation*}
        \int_G |(\omega(g)\phi,\phi')|^p dg<+\infty
    \end{equation*}
    for some $p<1$.
\end{proposition}
\begin{proof}
    For simplicity, we prove $\mathbb{F}=\mathbb{R}$ case. By $KAK$-decomposition, we only need to prove 
    \begin{equation*}
        \int\limits_{A^+} \Delta(a)|(\omega(a)\phi,\phi')|^p da=\int\limits_{A^+} \Delta(a)\prod_{1\leq i\leq n}\left(\frac{1}{a_i}\right)^{\frac{mp}{2}}\left|\int\limits_{M_{n\times m}(\mathbb{F})} \phi(a^{-1}x)\overline{\phi'(x)}dx\right|^p da< +\infty
    \end{equation*}
    where \[A^+:=\{(a_1,\dots,a_n)\in\mathbb{R}^n|a_1\geq a_2 \geq \dots \geq a_n>0\},~\Delta(a)=\prod_{1\leq i<j\leq n}\left(\frac{a_j}{a_i}-\frac{a_i}{a_j}\right)\prod_{1\leq i\leq n}\frac{1}{a_i},\] and $da$ is the Euclidean measure. Since $\phi$ is a Schwartz function, we only need to show for $A^+_{\epsilon}=\{(a_1,\dots,a_n)\in\mathbb{R}^n|a_1\geq a_2 \geq \dots \geq a_n>\epsilon\}$ and some $\epsilon>0$, the following integral
    \begin{equation*}
         \int\limits_{A^+_{\epsilon}} \prod_{1\leq i<j\leq n}\left(\frac{a_j}{a_i}-\frac{a_i}{a_j}\right)\prod_{1\leq i\leq n}\left(\frac{1}{a_i}\right)^{\frac{mp}{2}+1} da \quad \text{ is convergent}. 
    \end{equation*}
    The power of each $a_i$ is at most $-(\frac{mp}{2}+1)+(n-1)$ in each term. Thus, we can take suitable $p<1$, such that the power of each $a_i$ in each term is at most $-2+\eta$ for arbitrary $\eta>0$. The convergence follows from integration on each $a_i$ separately. 
\end{proof}
The positivity and the second point are proved by the characterization of the Hermitian form $(\cdot, \cdot)$ in the mixed model of the Weil representation. Let $V$ be an $n$-dimensional $\mathbb{F}$-vector space, and $W$ be an $m$-dimensional $\mathbb{F}$-vector space. We choose specific bases of $V$ and $W$ to make the identifications: \begin{equation} G \simeq \mathrm{GL}(V), \qquad G' \simeq \mathrm{GL}(W), \qquad M_{n \times m} \simeq \Hom_{\mathbb{F}}(V, W). \end{equation} We choose a direct sum decomposition of $W = W_1 \oplus W_2$, such that $\dim W_1 = \dim V$. Let $X = \Hom_{\mathbb{F}}(V, W_1)$. Note that $(G, \mathrm{GL}(W_2))$ is still a reductive dual pair, and its unitary Weil representation is denoted by $\omega_0$. Then $\mathcal{S}(X, \omega_0^{\infty} \hat{\otimes} V_{\sigma})$, as a $G$-representation with action given by \begin{equation*} (g \cdot \Phi)(x) := |\det g|^{-\frac{n}{2}} (\omega_0^{\infty} \hat{\otimes} \sigma)(g) \Phi(g^{-1} \cdot x), \end{equation*} is isomorphic to $\omega^{\infty} \hat{\otimes} V_{\sigma}$. Moreover, for $\Phi \in \mathcal{S}(X, \omega_0^{\infty} \hat{\otimes} V_{\sigma})$, we can define $\Phi_{\sigma}(g) := |\det g|^{n/2} \int_{G} (s \cdot \Phi)(g) ds$ as a constant function in $\mathcal{S}(X, \omega_0^{\infty} \hat{\otimes} V_{\sigma})$. We then observe: \begin{equation*} (\Phi_{\sigma}(x), \Phi'_{\sigma}(x)) = (\Phi, \Phi'), \end{equation*} where the left-hand side is the inner product of $\omega_0^{\infty} \hat{\otimes} V_{\sigma}$, and the right-hand side is the Hermitian form we defined above. As a result, the positivity of the Hermitian form follows. The irreducibility of $H(\sigma)$ and the identity $H(\sigma) = \theta(\sigma)$ follow formally from the argument in \cite{Li89}.

\end{document}